\documentclass[format=acmsmall, screen=true, review=false]{acmart}

\usepackage{booktabs}
\usepackage[mathscr]{euscript}
\usepackage{color}
\usepackage{graphicx}
\usepackage[linesnumbered,boxed,ruled,vlined]{algorithm2e}
\usepackage{tikz,pgfplots}
\usepackage{upgreek}
\usepackage{multirow}
\usepackage{yfonts}
\usepackage{mathtools}
\usepackage[normalem]{ulem}
\usepackage{latexsym}
\usepackage{url}
\usepackage{amsmath,amsbsy,amsfonts}
\usepackage{enumitem}
\usepackage{todonotes}
\usepackage{adjustbox}
\usepackage{arydshln}
\usepackage{listings}
\usetikzlibrary{math} 
\usetikzlibrary{spy,fit,matrix, positioning, intersections, math}
\usepgfplotslibrary{groupplots, statistics, colorbrewer}

\pgfplotsset{compat=1.16}

\definecolor{utorange}{rgb}{0.8,0.33,0.}
\definecolor{themec}{RGB}{51,108,121}
\definecolor{darkred}{rgb}{.6,.1,.1}
\definecolor{darkblue}{rgb}{.1,.1,.9}
\definecolor{greenback}{rgb}{.19,.94,.13}
\definecolor{orange}{rgb}{.76,.39,.13}
\definecolor{grass}{rgb}{.19,.64,.13}
\definecolor{sierp}{RGB}{209,28,209}
\definecolor{bgorange}{rgb}{1.,.95,.78}
\definecolor{grassgreen}{RGB}{92,135,39}
\definecolor{thinbox}{rgb}{.7,.8,1.}
\definecolor{darkgreen}{rgb}{0.25,0.65,0.10}
\definecolor{darkblue}{rgb}{0,0,0.7}
\definecolor{MitRed}{RGB}{163, 31, 52}
\definecolor{lakeblue}{RGB}{0,84,135}
\definecolor{chain}{RGB}{179,205,227}
\definecolor{kernel}{RGB}{254,217,166}
\definecolor{proposal}{RGB}{204,235,197}

\newcommand{\defeq}{\vcentcolon=}
\renewcommand{\vec}[1]{{\mathchoice
                     {\mbox{\boldmath$\displaystyle{#1}$}}
                     {\mbox{\boldmath$\textstyle{#1}$}}
                     {\mbox{\boldmath$\scriptstyle{#1}$}}
                     {\mbox{\boldmath$\scriptscriptstyle{#1}$}}}}

\newcommand{\mat}[1]{\mathbf{{#1}}}
\newcommand\restr[2]{{
  \left.\kern-\nulldelimiterspace 
  {#1}\vphantom{\big|} \right|_{#2}}}

\DeclareMathOperator*{\argmin}{argmin}
\DeclareMathOperator{\diag}{diag}

\newcommand{\R}{\mathbb{R}}

\newcommand{\D}{\mathcal{D}}

\newcommand{\bit}{\begin{itemize}}
\newcommand{\eit}{\end{itemize}}
\newcommand{\bdm} {\begin{displaymath}}
\newcommand{\edm} {\end{displaymath}}

\newcommand{\GM}[2]{\mathcal{N}\!\left( {#1}, {#2}\right)}

\newcommand{\iparspace}{\mathcal{M}}

\newcommand{\obs}{\vec{d}}
\newcommand{\ipar}{m}

\newcommand{\iparpr}{m_{\text{pr}}}
\newcommand{\iparmap}{m_{\scriptscriptstyle\text{MAP}}}

\newcommand{\ff}{\mathcal{F}} 
\newcommand{\qoi}{\mathcal{G}}
\newcommand{\tk}{K} 

\newcommand{\FF}{{\ensuremath{\mat{F}}}}

\newcommand{\prior} {\pi_{\mbox{\tiny prior}}}
\newcommand{\like}{\pi_{\mbox{\tiny like}}}

\newcommand{\post}{\pi_{\mbox{\tiny post}}}

\newcommand{\mupost}{\mu_{\mbox{\tiny post}}}
\newcommand{\muprior}{\mu_{\mbox{\tiny prior}}}
\newcommand{\mulaplace}{\hat{\mu}_{\mbox{\tiny post}}}

\newcommand{\ncov} {\bs{\Gamma}_{\!\mbox{\tiny noise}} }
\newcommand{\prcov} {\bs{\Gamma}_{\!\mbox{\tiny prior}} }
\newcommand{\postcov} {\bs{\Gamma}_{\!\mbox{\tiny post}} }

\newcommand{\Cprior}{\mathcal{C}_{\text{\tiny{prior}}}}

\newcommand{\map} {{m}_{\mbox{\tiny MAP}} }

\newcommand{\mpr} {{\vec{m}}_{\mbox{\tiny pr}} }

\newcommand{\Hmisfit}{ \mat{H}_{\mbox{\tiny misfit}} }

\newcommand{\mc}[1]{\mathcal{#1}}

\newcommand{\gbf}[1]{\boldsymbol{#1}}
\newcommand{\bs}[1]{\ensuremath{\boldsymbol{#1}}}

\renewcommand{\vec}[1]{\gbf{#1}}

\renewcommand{\matrix}[1] {\ensuremath{\boldsymbol{#1}}}

\usepackage{xspace}

\newcommand{\del}{\partial}
\newcommand{\Acal}{\mc{A}}

\newcommand{\ipart}{m_{\scriptscriptstyle\text{true}}}

\newcommand{\LRp}[1]{\left( #1 \right)}

\newcommand{\LRc}[1]{\left\{ #1 \right\}}

\newcommand{\half} {\ensuremath{\frac{1}{2}}}

\usepackage{pifont} %

\usepackage{hyperref}
\usepackage{url}

\pgfplotsset{every axis/.append style={
    tick label style={font=\Large},
    label style={font=\Large},
    legend style={font=\large}},
    every axis plot/.append style={ultra thick}
}

\begin{document}

\acmJournal{TOMS}

\title[hIPPYlib-MUQ: Integrating Data with Complex Predictive
  Models]{hIPPYlib-MUQ: A Bayesian Inference Software Framework for
  Integration of Data with Complex Predictive Models under
  Uncertainty}

\author{Ki-Tae Kim}
\orcid{TBD}
\affiliation{
\institution{University of California, Merced}
\department{Applied Mathematics, School of Natural Sciences}
\city{Merced}
\state{CA}
\country{USA}
}
\email{kkim107@ucmerced.edu}

\author{Umberto Villa}
\orcid{0000-0002-5142-2559}
\affiliation{
\institution{The University of Texas at Austin}
\department{Oden Institute for Computational Engineering \& Sciences}
\city{Austin}
\state{TX}
\country{USA}
}
\email{uvilla@austin.utexas.edu}

\author{Matthew Parno}
\orcid{0000-0002-9419-2693}
\affiliation{
\institution{Dartmouth College}
\department{Department of Mathematics}
\city{Hanover}
\state{NH}
\country{USA}
}
\email{matthew.d.parno@dartmouth.edu}

\author{Youssef Marzouk}
\orcid{**}
\affiliation{
\institution{Massachusetts Institute of Technology}
\department{Department of Aeronautics and Astronautics}
\city{Cambridge}
\state{MA}
\country{USA}
}
\email{ymarz@mit.edu}

\author{Omar Ghattas}
\orcid{}
\affiliation{
\institution{The University of Texas at Austin}
\department{Oden Institute for Computational Engineering \& Sciences}
\department{Department of Mechanical Engineering}
\department{Department of Geological Sciences}
\city{Austin}
\state{TX}
\country{USA}
}
\email{omar@oden.utexas.edu}

\author{Noemi Petra}
\orcid{0000-0002-9491-0034}
\affiliation{
\institution{University of California, Merced}
\department{Applied Mathematics, School of Natural Sciences}
\city{Merced}
\state{CA}
\country{USA}
}
\email{npetra@ucmerced.edu}



\def\biblio{}

\begin{abstract}
Bayesian inference provides a systematic framework for integration of
data with mathematical models to quantify the uncertainty in the
solution of the inverse problem. However, the solution of Bayesian inverse
problems governed by complex forward models described by partial
differential equations (PDEs) remains prohibitive with black-box
Markov chain Monte Carlo (MCMC) methods. We present hIPPYlib-MUQ, an
extensible and scalable software framework that contains
implementations of state-of-the art algorithms aimed to overcome the
challenges of high-dimensional, PDE-constrained Bayesian inverse
problems. These algorithms accelerate MCMC sampling by exploiting the
geometry and intrinsic low-dimensionality of parameter space via
derivative information and low rank approximation. The software
integrates two complementary open-source software packages, hIPPYlib
and MUQ. hIPPYlib solves PDE-constrained inverse problems using
automatically-generated adjoint-based derivatives, but it lacks full
Bayesian capabilities. MUQ provides a spectrum of powerful Bayesian
inversion models and algorithms, but expects forward models to come
equipped with gradients and Hessians to permit large-scale
solution. By combining these two complementary libraries, we created a
robust, scalable, and efficient software framework that realizes the
benefits of each and allows us to tackle complex large-scale Bayesian
inverse problems across a broad spectrum of scientific and engineering
disciplines. To illustrate the capabilities of hIPPYlib-MUQ, we
present a comparison of a number of MCMC methods available in the
integrated software on several high-dimensional Bayesian inverse
problems. These include problems characterized by both linear and
nonlinear PDEs, various noise models, and different
parameter dimensions. The results demonstrate that large ($\sim
50\times$) speedups over conventional black box and gradient-based
MCMC algorithms can be obtained by exploiting Hessian information
(from the log-posterior), underscoring the power of the integrated
hIPPYlib-MUQ framework.

\end{abstract}

\maketitle

\section{Introduction}
With the rapid explosion of observational and experimental data, a prominent
challenge is how to derive knowledge and insight from this data to make better
predictions and high-consequence decisions.  This question arises in all areas
of science, engineering, technology, and medicine, and in many cases, there are
mathematical models available that represent the underlying physical systems of
which the data is observed or measured.  These models are often subject to
considerable uncertainties stemming from unknown or uncertain input model
parameters (e.g., coefficient fields, constitutive laws, source terms,
geometries, initial and/or boundary conditions) as well as from noisy and
limited observations.  The goal is to infer these unknown model parameters from
observations of model outputs through corresponding partial differential
equation (PDE) models, and to quantify the uncertainty in the solution of such
inverse problems.

Bayesian inversion provides a systematic framework for integration of
data with complex PDE-based models to quantify uncertainties in model
parameter inference~\citep{KaipioSomersalo05, Tarantola05}.  In the
Bayesian framework, noisy data and, possibly uncertain, mathematical
models are integrated together with a prior information, yielding a
posterior probability distribution of the model parameters.  The
Markov chain Monte Carlo (MCMC) method is a common way to explore the
posterior distribution by use of sampling techniques.  However,
Bayesian inversion with complex forward models via conventional MCMC
methods faces several computational challenges.  First, characterizing
the posterior distribution of the model parameters or subsequent
predictions often requires repeated evaluations of expensive-to-solve
large-scale PDE models.  Second, the posterior distribution often has
a complex structure stemming from the nonlinear mapping from model
parameter to observed quantities.
Third, the parameters often are fields, which, after discretization,  lead to
very high-dimensional posteriors.  These difficulties make the solution of
Bayesian inverse problems with complex large-scale PDE forward models
computationally intractable.

Extensive research efforts have been devoted to overcome the prohibitiveness of
Bayesian inverse problems governed by large-scale PDEs.
With rapid progress in high-performance computing,
and advances in scalable PDE solvers, repeated
evaluations of forward PDE models for different input
parameters~\citep{petsc-web-page, trilinos-website} are becoming tractable. Furthermore, structure-exploiting MCMC
methods have effectively facilitated the exploration of complex posterior
distributions~\citep{CotterRobertsStuartEtAl12,
  BeskosGirolamiLanEtAl17,
  PetraMartinStadlerEtAl14,
  Bui-ThanhBursteddeGhattasEtAl12_gbfinalist}. Finally, dimension
reduction methods have proven to significantly reduce the computational cost of
MCMC simulations~\citep{CuiMarzoukWillcox16, ZahmCuiLawEtAl22}.  %
Applying and combining these advanced techniques can be extremely challenging. Therefore, a computational tool that will assist the computational and scientific community to apply, extend and tailor these methods will be very beneficial.

In this paper, we present a software framework to tackle large-scale Bayesian
inverse problems with PDE-based forward models, which has applications across a wide
range of science and engineering fields.  The software integrates two
open-source software packages, an Inverse Problems Python library
(hIPPYlib)~\citep{VillaPetraGhattas21} and the MIT
Uncertainty Quantification Library (MUQ)~\citep{Parno2014}, respecting their
attractive complementary capabilities.

hIPPYlib is an extensible software framework for the solution of
deterministic and linearized Bayesian inverse problems constrained by
complex PDE models.  Based on FEniCS \citep{LoggMardalWells12} for the
finite element approximation of PDEs and on PETSc
\citep{BalayAbhyankarAdamsEtAl14} for high-performance linear algebra
operations and solvers, hIPPYlib allows users to describe (and solve)
the underlying PDE-based forward model (required by the inverse
problem solver) in a relatively straightforward way.  hIPPYlib also
contains implementations of efficient numerical methods for the
solution of deterministic and linearized Bayesian inverse
problems. These include globalized inexact Newton-conjugate
gradient~\citep{AkcelikBirosGhattasEtAl06a, BorziSchulz12},
adjoint-based computation of gradients and Hessian
actions~\citep{Troltzsch10}, randomized linear
algebra~\citep{HalkoMartinssonTropp11}, and scalable sampling from
large-scale Gaussian fields.  The state-of-the-art algorithms implemented in
hIPPYlib efficiently deliver the solution of the linearized
Bayesian inverse problem.
hIPPYlib is, however, mainly designed for 
deterministic and linearized Bayesian inverse problems, and lacks full
Bayesian inversion capabilities.

MUQ complements hIPPYlib's capabilities with more support for the formulation
and solution of Bayesian inference problems.  MUQ is a modular software
framework designed to address uncertainty quantification problems involving
complex models. The software provides an abstract modeling interface for
combining physical (e.g., PDEs) and statistical components (e.g., additive
error models, Gaussian process priors, etc.) to define Bayesian posterior
distributions in a flexible and semi-intrusive way.  MUQ also contains a suite
of powerful uncertainty quantification algorithms including Markov chain Monte
Carlo (MCMC) methods~\citep{Parno2018}, transport maps~\citep{Marzouk2016},
likelihood-informed subspaces, sparse adaptive generalized polynomial chaos
(gPC) expansions~\citep{Conrad2013}, Karhunen-Lo\'eve expansions, Gaussian
process modeling~\citep{GPML2005, Hartikainen2010}, and prediction methods
enabling global sensitivity analysis and optimal experimental design.  To
effectively apply these tools to Bayesian inverse problems, however, MUQ needs
to be equipped with the type of gradient and/or Hessian information that
hIPPYlib can provide.

By interfacing these two software libraries, we aim to create a robust, scalable,
efficient, flexible, and easy-to-use software framework that overcomes the
computational challenges inherent in complex large-scale Bayesian inverse
problems.  Representative features of the software are summarized as follows:
\begin{itemize}
  \item The software combines the benefits of the two packages, hIPPYlib and
    MUQ, to enable scalable solution of
    Bayesian inverse problems governed by large-scale PDEs.
  \item Various advanced MCMC methods are available that can exploit problem
    structure (e.g., the derivative/Hessian information of the log-posterior).
  \item The software makes use of sparsity, low-dimensionality, and geometric
    structure of the log-posterior to achieve scalable and efficient MCMC
    methods.
  \item Convergence diagnostics are implemented to assess the quality
    of MCMC samples.
\end{itemize}

In the following sections, we first briefly review the Bayesian
formulation of inverse problems governed by PDEs both in
infinite-dimensional and in finite-dimensional spaces
(Section~\ref{sec:bayesian_framework}).  We then describe MCMC methods
used to characterize the posterior
(Section~\ref{sec:geometry_aware_mcmc}) and summarize convergence
diagnostics available in the software
(Section~\ref{sec:mcmc_diagnostics}).  Next, we present the design of
hIPPYlib-MUQ (Section~\ref{sec:software_framework}).  Finally, we
present numerous benchmark problems and a step-by-step implementation
guide to illustrate the key aspect of the present software
(Section~\ref{sec:numerical_illustration}). Section~\ref{sec:conclusions}
provides concluding remarks.

\section{The Bayesian inference framework}
\label{sec:bayesian_framework}

In this section, we present a brief discussion of the Bayesian
inference approach to solving inverse problems governed by PDEs. We
begin by providing an overview of the framework for
infinite-dimensional Bayesian inverse problems
following~\citet{Stuart10, Bui-ThanhGhattasMartinEtAl13a,
  PetraMartinStadlerEtAl14}.  Then we present a
brief discussion of the finite-dimensional approximations of the prior
and the posterior distributions; a lengthier discussion can be found
in~\citet{Bui-ThanhGhattasMartinEtAl13a}.  Finally, we present the
Laplace approximation to the posterior distribution, which requires
the solution of a PDE-constrained optimization problem for the
computation of the \emph{maximum a posteriori} (MAP).

\subsection{Infinite-dimensional Bayesian inverse problems}
\label{sec:infinite_dim_ip}

The objective of the inverse problem is to determine an unknown input
parameter field $\ipar$ that would give rise to given observational
(or experimental) noisy data $\obs$ by means of a (physics-based)
mathematical model.  In other words, given $\obs \in \mathbb{R}^q$, we
seek to infer $\ipar \in \iparspace$ such that
\begin{equation}\label{eq:noisemodel}
  \obs \approx \ff(\ipar),
\end{equation}
where $\ff: \iparspace \to \mathbb{R}^q$ is the {\it
  parameter-to-observable} map that predicts observations from a given
parameter $m$ through a forward mathematical model, and $\iparspace$
is an infinite-dimensional Hilbert space of functions defined on a
domain $\D \subset \mathbb{R}^d$.  Note that the evaluation of this
map involves solving the forward PDE model given $m$, followed by
extracting the observations from the solution of the forward
problem. In what follows, we assume that the forward equation residual
is continuously Fr\'echet-differentiable and its Jacobian a continuous
linear operator with continuous inverse~\citep{GhattasWillcox21}.

In the Bayesian approach, the inverse problem is framed as a
statistical inference problem.  The uncertain parameter $\ipar$ and
the observational data $\obs$ are deemed as random variables and the
solution is a conditional probability distribution $\mupost (\ipar |
\obs)$ that represents level of confidence in the estimation of the
parameter given the data.  The approach combines a {\it prior model}
reflecting our prior knowledge or assumptions about the parameters
before data are acquired, and a {\it likelihood model} representing
the probability that a given set of parameters might give rise to the
observed data.

Using the Radon-Nikodym derivative~\citep{Williams1991} of the posterior measure
$\mupost (\ipar)$ with respect to the prior measure $\muprior (\ipar)$, Bayes' theorem in
infinite dimensions is stated as
\begin{equation} \label{eq:bayes-abstract}
   \frac{d\mupost}{d\muprior} \propto \like(\obs | \ipar),
\end{equation}
where $\like$ denotes the likelihood function.  For detailed conditions under
which the posterior measure is well-defined, we refer the reader to~\citet{Stuart10}.

Without loss of generality, the probability density function of the likelihood
can be expressed as
\begin{equation}\label{eq:likelihood}
  \like(\obs | \ipar) \propto \exp\Big\{-\Phi(\ipar; \obs)\Big\}.
\end{equation}
The negative log-likelihood function $\Phi(\ipar; \obs)$ has different forms
depending on how one models the noise that stems from measurement uncertainties
and/or modeling errors; for example, in the case of an additive Gaussian noise
model $\obs = \ff(\ipar) + \vec{\eta}$ with a Gaussian noise random variable
$\vec{\eta} \in \mathbb{R}^q$ with zero mean and covariance matrix $\ncov \in
\mathbb{R}^{q \times q}$,
it has the
form $\Phi(\ipar; \obs) = \half \| \ff(\ipar) - \vec{d} \|^2_{\ncov^{-1}}$, where $\|
\cdot \|_{\ncov^{-1}}$ denotes the $L^2$ norm weighted by the inverse noise
covariance $\ncov^{-1}$.

We take the prior to be a Gaussian measure\footnote{The use of a
  Gaussian measure (often with Mat\'ern covariance covariance kernel)
  is a common choice in infinite-dimensional inverse problems, where
  the inversion parameter is a spatially varying field. In
  hIPPYlib-MUQ, non-Gaussian prior probability distributions can be
  handled in a case-by-case basis by introducing a nonlinear
  \emph{invertible} change of variables to map the desired prior
  distribution to a Gaussian one.}, i.e., $\muprior =
\GM{\iparpr}{\Cprior}$, and assume that samples from the prior
distribution are square-integrable functions in the domain $\D$, i.e.,
belong to $L^2(\D)$.  The prior covariance operator $\Cprior$ is
constructed to be a trace-class operator to guarantee bounded variance
of samples from the prior distribution and well-posedness of the
Bayesian inverse problem~\citep{Stuart10,
  Bui-ThanhGhattasMartinEtAl13a, VillaPetraGhattas21} for detailed
explanation.  Specifically, we take the prior covariance operator as
the inverse of the $v$th power of a Laplacian-like operator, namely
$\Cprior := \Acal^{-v} = (-\gamma \Delta + \delta I)^{-v}; \ v >
\frac{d}{2}$. Here $\gamma$ and $\delta > 0$ control the correlation
length and the pointwise variance of the prior
operator~\citep{VillaPetraGhattas21, LindgrenRueLindstroem11}. These
choices of prior ensure that $\Cprior$ is a trace-class operator,
guaranteeing bounded pointwise variance and a well-posed
infinite-dimensional Bayesian inverse problem~\citep{Stuart10,
  Bui-ThanhGhattasMartinEtAl13a}.

\subsection{Discretization of the Bayesian formulation}

Here, we briefly present the finite-dimensional approximation of the Bayesian
formulation described in the previous section. %
We consider a finite-dimensional subspace $\iparspace_h$ of
$\iparspace$, defined by the span of a set of basis functions
$\LRc{\phi_j}_{j=1}^n$.  For example, in hIPPYlib-MUQ, these basis
functions are globally continuous piecewise polynomials on each
element of a mesh discretization of the domain
$\D$~\citep{BeckerCareyOden81,StrangFix88}. These are the natural
choice for the discretization of the elliptic operator $\Acal$ used in
the definition of the prior covariance.  The parameter field $\ipar$
is then approximated as $\ipar \approx \ipar_h =
\sum_{j=1}^nm_j\phi_j$, and, in what follows, $\bs m
=\LRp{m_1,\hdots,m_n}^T\in \R^n$ denotes the vector of the finite
element coefficients of $\ipar_h$.

The finite-dimensional approximation of the prior measure $\muprior$ is now
specified by the density
\begin{equation}
  \prior(\vec{m})
  \propto\exp\Big(
    -\frac{1}{2}\|\vec{m}-\mpr\|^2_{\prcov^{-1}}
  \Big),
\end{equation}
where $\mpr \in \mathbb{R}^n$ and $\prcov \in \mathbb{R}^{n \times n}$ are the
mean vector and the covariance matrix that arise upon discretization of
$\iparpr$ and $\Cprior$, respectively.
We refer the reader to~\citet{Bui-ThanhGhattasMartinEtAl13a,
VillaPetraGhattas21} for the explicit expression of the prior covariance matrix
$\prcov$.

Then the Bayes' theorem for the density of the finite-dimensional approximation
of the posterior measure $\mupost$ is given by
\begin{equation}\label{eq:posterior}
  \post(\vec{m}) \defeq \post(\vec{m} | \obs) \propto \like (\obs | \vec{m})
  \prior(\vec{m}).
\end{equation}
The finite-dimensional posterior probability density function can be expressed
explicitly as
\begin{equation}\label{eq:posterior2}
  \post(\vec{m})
  \propto\exp\Big(
  -\Phi(\vec{m}; \obs)
  -\frac{1}{2}\|\vec{m}-\mpr\|^2_{\prcov^{-1}}
  \Big).
\end{equation}
Here, evaluating $\Phi(\vec{m}; \obs)$ requires constructing a finite dimensional discretization of the parameter-to-observable map, $\FF(\vec{m})$.

\subsection{The Laplace approximation of the  posterior distribution}

In general, the posterior probability distribution~\eqref{eq:posterior2} is not
Gaussian due to the nonlinearity of the parameter-to-observable map.
In this section, we discuss the solution to the so-called \emph{linearized} Bayesian inverse problem
by use of the Laplace approximation. %
The Laplace approximation amounts to constructing a Gaussian distribution
centered at the \emph{maximum a posteriori} (MAP) point.  The MAP point
represents the most probable value of the parameter vector over the posterior,
 i.e., %
\begin{equation}\label{eq:objfunction-bayesian}
  \vec{\map} := \argmin_{\vec{m}} (- \log \post(\vec{m})) = \argmin_{\vec{m}}
  \biggl[ \Phi(\vec{m}; \obs)
  + \frac{1}{2}\|\vec{m}-\mpr\|^2_{\prcov^{-1}}\biggr].
\end{equation}
The covariance matrix of the Laplace approximation is the inverse of the
Hessian of
the negative log-posterior evaluated at the MAP point.
Then under the Laplace approximation, the solution of the linearized Bayesian inverse problem is given by
\begin{equation}\label{eq:lapl_discrete}
  \hat{\pi}_{\text{post}}(\vec{\ipar}) \sim \mathcal{N}(\vec{\map},\postcov),
\end{equation}
with
\begin{equation}
  \postcov
  = \mat{H}^{-1} (\vec{\map})
  = \left(\Hmisfit(\vec{\map}) + \prcov^{-1}
  \right)^{-1},
  \label{eq:postcov_discrete}
\end{equation}
where $\mat{H}(\vec{\map})$ and $\Hmisfit(\vec{\map})$ denote the Hessian
matrices of, respectively, the negative log-posterior and the negative
log-likelihood evaluated at the MAP point; see~\citet{VillaPetraGhattas21} for
a derivation of this Hessian using the adjoint-based method.

The quality of the Gaussian approximation of the posterior depends on
the degree of nonlinearity in the parameter-to-observable map, the
noise covariance matrix, and the number of
observations~\citep{Gelman2004, Tarantola05,
  Bui-ThanhGhattasMartinEtAl13a, IsaacPetraStadlerEtAl15a,
  EvansSwartz00, Press03, Stigler86, TierneyKadane86, Wong01}. When
the parameter-to-observable map is linear and the additive noise and
prior models are both Gaussian, the Laplace approximation is exact.
Even if the parameter-to-observable map is significantly nonlinear, the Laplace
approximation is a crucial ingredient to achieve scalable, efficient, and
accurate posterior sampling with MCMC methods, as we will discuss in the
following section.

Note that the Laplace approximation involves the Hessian of the negative
log-likelihood (the data misfit part of the Hessian), which may not be explicitly constructed
when the parameter dimension is large. %
However, the data typically provide only limited information about the
parameter field, and thus the eigenspectrum of the Hessian matrix often decays
very rapidly.
We exploit this compact nature of the Hessian to overcome its prohibitive
computational cost, and construct a low-rank approximation of the data misfit
Hessian matrix using a matrix-free method (such as the randomized subspace
iteration~\citep{HalkoMartinssonTropp11}).

Concretely, we construct a low-rank approximation of the data misfit Hessian,
i.e., $\Hmisfit \approx \prcov^{-1} \mat{V}_r \boldsymbol{\Lambda}_r \mat{V}_r^T
\prcov^{-1}$, where $\boldsymbol{\Lambda}_r = \text{diag}(\lambda_1, \ldots,
\lambda_r) \in \mathbb{R}^{r \times r}$ and $\mat{V}_r = [\vec{v}_1, \ldots,
\vec{v_r}] \in \mathbb{R}^{n \times r}$ contain the $r$ largest eigenvalues and
corresponding eigenvectors, respectively, of the generalized symmetric
eigenvalue problems
\begin{equation}\label{eq:eigenproblem_1}
  \Hmisfit \vec{v}_i = \lambda_i \prcov^{-1} \vec{v}_i; \quad i = 1, \ldots, n.
\end{equation}
Note that the eigenvectors $\vec{v}_i$ are orthonormal with respect to
$\prcov^{-1}$, that is $\vec{v}_i^T \prcov^{-1} \vec{v}_j = \delta_{ij}$, where
$\delta_{ij}$ is the Kronecker delta.
With this row-rank approximation and using the Sherman-Morrison-Woodbury
formula~\citep{GolubVan96}, we obtain, for the inverse of the Hessian
in~\eqref{eq:postcov_discrete},
\begin{equation}
  \mat{H}^{-1} = \left(\Hmisfit +  \prcov^{-1}\right)^{-1} = \prcov -
  \mat{V}_r \mat{D}_r \mat{V}_r^T + \mathcal{O}\left(\sum_{i=r+1}^n \frac{
  \lambda_i }{1 + \lambda_i}\right),
  \label{eq:sm}
\end{equation}
where $\mat{D}_r =  \diag(\lambda_1/(\lambda_1+1), \dots,
\lambda_r/(\lambda_r+1)) \in \R^{r\times r}$.
We refer the reader to~\citet{VillaPetraGhattas21} for a detailed description
of how the randomized algorithm~\citep{HalkoMartinssonTropp11} proceeds to
construct the low-rank approximation of the Hessian, including the associated
computational complexity.

We can see from the last remainder term in~\eqref{eq:sm} that to
obtain an accurate low-rank approximation of $\mat{H}^{-1}$, we must
keep eigenvectors corresponding to eigenvalues that are greater than
1.  This approximation is used to efficiently perform various
operations related to the Hessian, for example, applying the square-root
inverse of the Hessian to a vector, which is needed to draw samples
from the Gaussian approximation discussed in this section. We remark
that the efficiency and scalability of our approach is based on the
low rankness of the data misfit Hessian. We refer the reader
to~\citet{GhattasWillcox21}, where this argument is made via model
problems where low-rankness can be analytically shown, and for
references to more complex problems where it can be shown empirically.

\section{MCMC Sampling}
\label{sec:geometry_aware_mcmc}
As mentioned above, when the parameter-to-observable map is nonlinear,
the Laplace approximation may be a poor approximation of the posterior.
In this case, one needs to apply a sampling-based method to explore the full
posterior. In this section, we focus on several advanced Markov chain Monte Carlo (MCMC) methods
available in the present software.  We outline the general structure of MCMC
methods with a brief discussion of their key features.  We then present various
diagnostics to assess the convergence of MCMC simulations.

\subsection{Markov chain Monte Carlo}

MCMC provides a flexible framework for exploring the posterior
distribution.  It generates samples from the posterior distribution that
can be employed in Monte Carlo approximations of posterior expectations.  For example, the posterior expectation of a quantity of
interest $\qoi(\ipar)$ can be approximated by
\begin{equation}
  \int \qoi(\ipar)\, d\mupost \approx \frac{1}{N}\sum_{i=1}^N
  \qoi\left(\ipar_i\right),
  \label{eq:mcest}
\end{equation}
where each $\ipar_i\sim \mupost$ is a sample of the posterior distribution.

MCMC techniques construct ergodic Markov chains where the posterior
distribution is the unique stationary distribution of the
chain~\citep{RobertCasella05}.
Asymptotically, the states of the Markov chain are therefore exact samples of
the posterior distribution and can be used in \eqref{eq:mcest}.
Markov chains are defined in terms of a
transition kernel, which is a position dependent probability distribution
$\tk(\cdot | \vec{\ipar}_i)$ over state $\vec{\ipar}_{i+1}$ in the chain given
the previous state $\vec{\ipar}_i$, i.e. $\vec{\ipar}_{i+1} \sim \tk(\cdot |
\vec{\ipar}_i)$.   Note that chains of finite length must be employed in
practice and the statistical accuracy of the Monte Carlo estimator is therefore
highly dependent on the ability of the transition kernel to efficiently explore
the parameter space.

There are several frameworks for constructing transition kernels that
are appropriate for MCMC, including the well known Metropolis-Hastings
(MH) rule~\citep{MetropolisRosenbluthRosenbluthEtAl53,Hastings70},
Gibbs sampler (e.g.,~\citep{CasellaGeorge92}), and delayed rejection
(DR)~\citep{Mira2001}.  MUQ provides implementations of these
frameworks, as well as the generalized Metropolis-Hastings (gMH)
kernel~\citep{calderhead2014general} and multilevel MCMC framework
of~\citep{dodwell2019multilevel}.  Most of these frameworks start by
drawing samples from one or more proposal distributions $q_1(\cdot |
\vec{\ipar}_i),\, \ldots,\, q_K(\cdot | \vec{\ipar}_i)$ that are easy
to sample from (e.g., Gaussian) and then ``correct'' the proposed
samples to obtain exact, but correlated, posterior samples.  In the MH
and DR kernels, corrections take the form of accepting or rejecting
the proposed point.  In the gMH kernel, the correction involves
analytically sampling a finite state Markov chain over multiple
proposed points.  Intuitively, proposal distributions that capture the
shape of the posterior, either locally around $m_i$ or globally over
the parameter space, tend to require fewer ``corrections'' and yield
more efficient algorithms.

\paragraph{Proposal Distributions}
Let $q(\cdot | \vec{\ipar}_i)$ denote a proposal distribution that is
``parameterized" by the current state of the chain $\vec{\ipar}_i$.   We
require that the proposal distribution is easily sampled and that its density
can be efficiently evaluated.  The MH rule
\citep{MetropolisRosenbluthRosenbluthEtAl53, Hastings70} defines a transition
kernel $\tk_{MH}(\cdot | \vec{\ipar}_i)$ through a two step process:  first
draw a random sample $\vec{\ipar}^\prime \sim q(\cdot | \vec{\ipar}_i)$ from
the proposal distribution, and then accept the proposed sample
$\vec{\ipar}^\prime$ as the next step in the chain $\vec{\ipar}_{i+1}$ with a
probability $\alpha = \min \{1, \gamma\}$ where $\gamma = \frac{
\post(\vec{\ipar}^\prime)}{\post(\vec{\ipar}_i)} \frac{q(\vec{\ipar}_i |
\vec{\ipar}^\prime)}{q( \vec{\ipar}\prime | \vec{\ipar}_i)}$.
If rejected, set $\vec{\ipar}_{i+1}=\vec{\ipar}_i$.   Under mild technical
conditions on the proposal distribution (see e.g., \citet{roberts2004general}),
the MH rule defines a Markov chain that is ergodic and has $\mupost$ as a
stationary distribution, thus enabling states in the chain to be used in Monte
Carlo estimators.  Note that the detailed balance condition (see e.g.,
\citet{owen2013monte}) is commonly employed to verify that a Markov chain has
$\mupost$ as a stationary distribution, but this condition alone is not
sufficient to guarantee that the chain will converge to the stationary
distribution.  See \citet{roberts2004general} for a detailed discussion of MH
convergence and convergence rates.

While the MH rule will yield a valid MCMC kernel for a
large class of proposal distributions, the dependence of the proposal on the
previous state, combined with possible rejection of the proposed state, results
in inter-sample correlations in the Markov chain.   Because of these
correlations, the error of the Monte Carlo approximation in \eqref{eq:mcest}
will be larger when using MCMC than in the classic Monte Carlo setting with
independent samples.   Markov chains with large correlations will result in
larger estimator variance.
To reduce correlation in the Markov chain, we seek
proposal distributions that can take large steps with a high probability of
acceptance.  From the acceptance probability in the MH rule
we see that this can occur when the proposal density $q(\vec{\ipar} |
\vec{\ipar}_i)$ is a good approximation to $\post(\vec{\ipar})$, so
that $\gamma$ is close to 1.

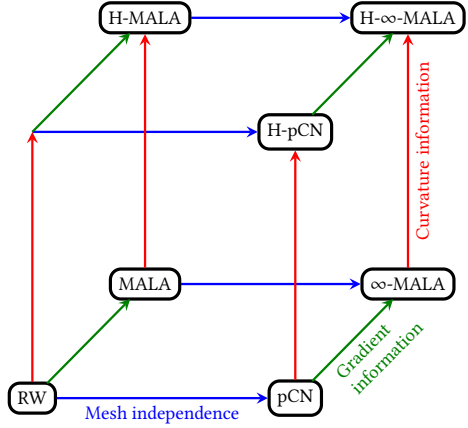
\begin{figure}[tb]
  \begin{tikzpicture}[->, >=stealth, line width=0.2ex,
  every node/.style={font=\fontsize{9}{10.2}\selectfont, scale=0.8},
  mybox/.style={very thick, rounded corners, inner sep=4}]

  \node [] (0) at (0, 0) {};
  \node [] (1) at (3.5, 0) {};
  \node [] (2) at (0, 3.5) {};
  \node [] (3) at (3.5, 3.5) {};
  \node [] (4) at (1.5, 1.5) {};
  \node [] (5) at (5, 1.5) {};
  \node [] (6) at (1.5, 5) {};
  \node [] (7) at (5, 5) {};

  \node [mybox, draw=black, name path=random_walker] (random_walker) at (0) {RW};
  \node [mybox, draw=black, name path=pcn] (pcn) at (1) {pCN};
  \node [mybox, draw=black, name path=hpcn] (hpcn) at (3) {H-pCN};
  \node [mybox, draw=black, name path=mala] (mala) at (4) {MALA};
  \node [mybox, draw=black, name path=infmala] (infmala) at (5) {$\infty$-MALA};
  \node [mybox, draw=black, name path=hmala] (hmala) at (6) {H-MALA};
  \node [mybox, draw=black, name path=hinfmala] (hinfmala) at (7) {H-$\infty$-MALA};

  \draw [color=blue] (random_walker.east) to node[below] {Mesh independence} (pcn.west);
  \draw [color=blue] (mala.east) to (infmala.west);
  \draw [color=blue] (2.center) to (hpcn.west);
  \draw [color=blue] (hmala.east) to (hinfmala.west);

  \draw [draw=none, name path=line1] (0.center) to (4.center);
  \node [draw=none, name intersections={of=line1 and random_walker, by={int1}}] {};
  \node [draw=none, name intersections={of=line1 and mala, by={int2}}] {};
  \draw [color=green!50!black] (int1.center) to (int2.center);

  \draw [draw=none, name path=line2] (1.center) to (5.center);
  \node [draw=none, name intersections={of=line2 and pcn, by={int3}}] {};
  \node [draw=none, name intersections={of=line2 and infmala, by={int4}}] {};
  \draw [color=green!50!black] (int3.center) to node[sloped, below, align=left,
    xshift=0.2cm] {Gradient \\ information} (int4.center);

  \draw [draw=none, name path=line3] (3.center) to (7.center);
  \node [draw=none, name intersections={of=line3 and hpcn, by={int5}}] {};
  \node [draw=none, name intersections={of=line3 and hinfmala, by={int6}}] {};
  \draw [color=green!50!black] (int5.center) to (int6.center);

  \draw [draw=none, name path=line4] (2.center) to (6.center);
  \node [draw=none, name intersections={of=line4 and hmala, by={int7}}] {};
  \draw [color=green!50!black] (2.center) to (int7.center);

  \draw [color=red] (random_walker.north) to (2.center);
  \draw [color=red] (pcn.north) to (hpcn.south);
  \draw [color=red] (infmala.north) to node[sloped, below] {Curvature information} (hinfmala.south);
  \draw [color=red] (mala.north) to (hmala.south);


\end{tikzpicture}
  \caption{The relationship of various MCMC proposal distributions with respect
  to mesh-refinement independence (blue arrow), gradient awareness (green
  arrow), and curvature awareness (red arrow). The abbreviations stand for the
  following MCMC proposals: RW for random walk, pCN for preconditioned
  Crank-Nicolson, MALA for Metropolis-adjusted Langevin algorithm, H-pCN for
  curvature-informed pCN, H-MALA for curvature-informed MALA, $\infty$-MALA for
  infinite-dimensional MALA, and H-$\infty$-MALA for curvature-informed
  infinite-dimensional MALA.}
\label{fig:MCMC_prop}
\end{figure}

We now turn to describing specific proposal distributions used in hIPPYlib-MUQ.
First, we begin by describing common proposal mechanisms that exploit gradient
and curvature information to accelerate sampling in finite-dimensional spaces.
These algorithms comprise the left face of the cube in Figure
\ref{fig:MCMC_prop}.  We then show how these ideas can be extended to construct
proposals with performance that is independent of mesh-refinement (i.e., independent of dimension), thus
``lifting'' the derivative-accelerated proposals to an infinite-dimensional
setting.   This ``lifting'' operation transforms proposals on the left face in
Figure \ref{fig:MCMC_prop} to their dimension-independent analogs on the right
face of the proposal cube.

\paragraph{Exploiting Gradient and Curvature Information}
Perhaps the simplest and most common, but not generally efficient, proposal distribution takes the form of
a Gaussian distribution centered at the current state in the chain,
\begin{equation}
  q_{\text{RW}}(\vec{\ipar} | \vec{\ipar}_i) = \mathcal{N}\left(\vec{\ipar}_i,
  \boldsymbol{\Gamma}_{\text{prop}} \right),
  \label{eq:rwprop}
\end{equation}
where $\boldsymbol{\Gamma}_{\text{prop}} \in \mathbb{R}^{n
\times n}$ is a user defined covariance matrix.  When used
with the MH rule, this \emph{random walk} (RW) proposal yields an
MCMC algorithm that is commonly called the random walk Metropolis
algorithm.   The adaptive Metropolis (AM) algorithm employs a variant of this
proposal where the covariance $\boldsymbol{\Gamma}_{\text{prop}}$ is adapted based on
previous samples~\citep{Haario2001}.  A proposal covariance $\boldsymbol{\Gamma}_{\text{prop}}$ that
matches posterior covariance increases efficiency, but the random walk proposal
is still a poor approximation of the posterior density.%

A slightly more efficient proposal can be obtained through a one-step
Euler-Maruyama discretization of the Langevin stochastic differential equations
\citep{RobertsStramer03}.  The resulting Langevin proposal takes the form
\begin{equation}
 q_{\text{MALA}}(\vec{\ipar} | \vec{\ipar}_i) = \mathcal{N}\left( \vec{\ipar}_i + \tau
   \boldsymbol{\Gamma}_{\text{prop}} \nabla \log \post(\vec{\ipar}_i),\, 2\tau
 \boldsymbol{\Gamma}_{\text{prop}} \right),
 \label{eq:malaprop}
\end{equation}
where $\tau$ is the step size parameter.
MH samplers with this proposal are called \emph{Metropolis-adjusted Langevin
algorithms} (MALA).      
Like the AM algorithm, adapting the covariance of the MALA proposal can also
improve performance~\citep{Atchade06, MarshallRoberts12}.  It is also common 
to use an approximation of the posterior covariance, such as the inverse of 
the log-posterior Hessian, to help the MALA proposal capture the posterior 
correlation.  In this work for example, we employ a low rank-based approximation 
of the log-posterior Hessian at the MAP point (c.f. eq. \eqref{eq:sm})
\begin{equation}
  q_{\text{H-MALA}}(\vec{\ipar} | \vec{\ipar}_i) = \mathcal{N}\left(\vec{\ipar}_i + \tau
  \mat{H}^{-1} \nabla \log \post(\vec{\ipar}_i), 2\tau  \mat{H}^{-1}\right).
\end{equation}
This metric is similar to the one used by \citet{martin2012stochastic} and is equivalent to the preconditioned MALA proposal in
\eqref{eq:malaprop} using the covariance of the Laplace approximation in
\eqref{eq:postcov_discrete}. 

Both \eqref{eq:rwprop} and \eqref{eq:malaprop} use a covariance that is
constant across the parameter space.   Allowing this covariance to adapt to the
local correlation structure of the posterior density enables higher order
approximations to be obtained, resulting in more efficient MCMC algorithms.
In~\citet{girolami2011riemann}, a differential geometric viewpoint was employed to define a family of proposal mechanisms on a Riemannian manifold.
Adapting the MALA proposal in \eqref{eq:malaprop} to this manifold
setting and ignoring the manifold's curvature, results in
\begin{equation}
  q_{\text{sMMALA}}(\vec{\ipar} | \vec{\ipar}_i) = \mathcal{N}\left(\vec{\ipar}_i + \tau
  \mathbf{G}^{-1}(\vec{\ipar}_i) \nabla \log \post(\vec{\ipar}_i), 2\tau
\mathbf{G}^{-1}(\vec{\ipar}_i)\right),
\end{equation}
where $\mathbf{G}(\vec{\ipar})$ is a position-dependent metric tensor.  This is known as the simplified
Manifold MALA (sMMALA) proposal. \citet{girolami2011riemann} defined the metric tensor $\mathbf{G}(\vec{\ipar})$ using the expected Fisher information metric, which provides
a positive definite approximation of the posterior covariance at the point
$\vec{\ipar}$.    

Hamiltonian Monte Carlo techniques, including the No-U-Turn Sampler
(NUTS)~\citep{hoffman2014no}, are another important class of MCMC
proposals.  These techniques approximately solve a Hamiltonian system
to take large jumps in parameter space.  While efficient in many
scenarios (see e.g.,~\citet{Neal10}), especially with purely
statistical models, we have found that solving the Hamiltonian system
involves an intractable number of posterior gradient evaluations on
our PDE-based problems of interest.  The transport map MCMC algorithms
of~\citet{Parno2018} are also not considered here because of the
challenge of building high-dimensional transformations.

\paragraph{Dimension-Independent Proposal Distributions}
For finite-dimensional parameters, the random walk and MALA proposals
defined above can be used with the MH rule for MCMC.  However, their
performance is not discretization invariant.  As the discretization of
the function $\ipar$ is refined, the performance of the samplers on
the finite-dimensional posterior $\post(\vec{\ipar})$ will worsen.  As
the dimension increases, the difference between the largest two
eigenvalues of the MCMC transition kernel (i.e. the spectral gap),
goes to zero and the mixing times of the chains grows indefinitely;
see~\citet{hairer2014spectral,CotterRobertsStuartEtAl12} for details.
Some modifications to the proposals are necessary to obtain
``dimension-independent'' performance.  The works
of~\citet{CotterRobertsStuartEtAl12}, \citet{BeskosGirolamiLanEtAl17},
and \citet{bardsley2020scalable}, for example, modify existing
finite-dimensional proposals to ensure the algorithm performance is
independent of mesh refinement.

The dimension-independent analog of the RW proposal is the preconditioned
Crank-Nicolson (pCN) proposal introduced in~\citet{CotterRobertsStuartEtAl12}.   It
takes the form
\begin{equation}
  q_{\text{pCN}}(\vec{\ipar} | \vec{\ipar}_i) = \mathcal{N}\left( \mpr +
\sqrt{1-\beta^2}(\vec{\ipar}_i-\mpr),\, \beta^2 \prcov \right).
 \label{eq:pcnprop}
\end{equation}
Notice that when $\beta=1$, the pCN proposal is equal to the prior
distribution.  The MALA proposal was also adapted
in~\citet{CotterRobertsStuartEtAl12} to obtain the
infinite-dimensional MALA ($\infty$-MALA) proposal
\begin{equation}
  q_{\text{MALA}}^{\infty}(\vec{\ipar} | \vec{\ipar}_i) = \mathcal{N}\left(\sqrt{1 -
    \beta^2}\vec{\ipar}_i + \beta \frac{\sqrt{h}}{2} \left(\mpr - \prcov \nabla
  \Phi(\vec{\ipar}_i)\right), \, \beta^2 \prcov\right),
\label{eq:infmala}
\end{equation}
where $\beta = 4\sqrt{h} / (4 + h)$ and $h$ is a parameter that can be
tuned.  While the pCN and $\infty$-MALA proposals result in
discretization-invariant Metropolis-Hastings algorithms, they suffer
from the same deficiencies as their finite-dimensional RW and MALA
analogs, i.e., they do not capture the posterior geometry.

Several efforts have worked to minimize this deficiency, see for
example~\citet{BeskosGirolamiLanEtAl17, RudolphSprungk16,
  PinskiSimpsonStuartEtAl15, PetraMartinStadlerEtAl14}.  We consider a
generalization of the pCN proposal described
in~\citet{PinskiSimpsonStuartEtAl15}.  It incorporates the MAP point
and the posterior curvature information at that point into the pCN
proposal, which is denoted by H-pCN and takes the form
\begin{equation}
  q_{\text{H-pCN}}(\vec{\ipar} | \vec{\ipar}_i) = \mathcal{N}\left( \vec{\map} +
  \sqrt{1-\beta^2}(\vec{\ipar}_i-\vec{\map}),\, \beta^2 \mat{H}^{-1} \right).
 \label{eq:hpcn}
\end{equation}
Another method that can exploit the posterior geometry is an extension of
the $\infty$-MALA proposal discussed in~\citet{BeskosGirolamiLanEtAl17}:
\begin{equation}
  q_{\text{sMMALA}}^{\infty}(\vec{\ipar} | \vec{\ipar}_i) = \mathcal{N}\left(
 \mu^\prime(\vec{\ipar}_i),\, \Gamma^\prime(\vec{\ipar}_i)\right),
 \label{eq:infsmmala}
\end{equation}
where
\begin{eqnarray}
 \mu^\prime(\vec{\ipar}_i) &=&
 \sqrt{1 - \beta^2}\vec{\ipar}_i + \beta \frac{\sqrt{h}}{2} \left(
   \vec{\ipar}_i -
   \mathbf{G}^{-1} \prcov^{-1} (\vec{\ipar}_i-\mpr) - \mathbf{G}^{-1} \nabla
   \Phi(\vec{\ipar}_i)
 \right)\\
 \Gamma^\prime(\vec{\ipar}_i)  &=& \beta^2 \mathbf{G}^{-1}(\vec{\ipar}_i).
\end{eqnarray}
This $\infty$-sMMALA proposal simplifies to $\infty$-MALA when
$\mathbf{G}^{-1}(\vec{\ipar}_i)=\prcov$.
When $\mathbf{G}(\vec{\ipar}_i)$ is
the Laplace approximation Hessian from \eqref{eq:postcov_discrete}, the
$\infty$-sMMALA proposal simplifies to
\begin{equation}
  q_{\text{H-MALA}}^\infty(\vec{\ipar} | \vec{\ipar}_i) = \mathcal{N}\left(
 \sqrt{1 - \beta^2}\vec{\ipar}_i + \beta \frac{\sqrt{h}}{2} \left(
   \vec{\ipar}_i -
   \mat{H}^{-1} \prcov^{-1} (\vec{\ipar}_i-\mpr) - \mat{H}^{-1} \nabla
   \Phi(\vec{\ipar}_i)
 \right),
 \beta^2 \mat{H}^{-1} \right),
 \label{eq:hinfmala}
\end{equation}
which we denote by H-$\infty$-MALA.

\paragraph{Alternative Transition Kernels}
The proposal distributions above are classically considered in the
context of a Metropolis-Hasting kernel.  However, there are
alternative transition kernels that also result in ergodic Markov
chains.  Here we consider transition kernels constructed from the
delayed rejection approach of~\citet{Mira2001} as well as
Metropolis-within-Gibbs kernels, which repeatedly use the
Metropolis-Hastings rule on different conditional slices of the
posterior distribution to construct the Markov chain.  In particular,
we consider the family of dimension-independent likelihood-informed
(DILI) approaches~\citep{CuiLawMarzouk16, cui2021data}, which define a
Metropolis-within-Gibbs sampler that inherits dimension-independent
properties from an appropriate dimension-independent proposal. By
dimension-independence here we mean that the acceptance rate and mixing
properties will not deteriorate when the dimension of the problem
increases.

The delayed rejection kernel allows multiple proposals to be attempted in each
step of the Markov chain.  This can be advantageous when using multiple
proposals with complementary properties.  For example, it is possible to start
with a proposal that attempts to make large ambitious jumps across the
parameter space but may have low acceptance probability while falling back on a
more conservative proposal that takes smaller steps with a larger probability
of acceptance.  Similarly, it is possible to start with a proposal that is more
computationally efficient (e.g., does not require gradient information) but
less likely to be accepted, while employing a more expensive proposal mechanism
in a second stage to ensure the chain explores the space.  In either case, if
the first proposed move is rejected by the Metropolis-Hastings rule, another
more
expensive proposal that is more likely to be accepted can be tried with an
adjusted acceptance probability.  More than two stages can also be employed.  

DILI divides the parameter space into a finite-dimensional subspace, which can
be explored with standard proposal mechanisms, and a complementary infinite-dimensional
space that can be explored with a dimension-independent approach,
such as those described above.  The resulting transition kernel is more
complicated than the Metropolis-Hastings rule, but inherits the
dimension-independent properties of the complementary space proposal.  The
likelihood-informed subspace is computed using the generalized eigenvalue
problem~\eqref{eq:eigenproblem_1}.  If the eigenvalue is larger than one, it
indicates that the likelihood function dominates the prior density in that
direction.  The same low rank structure used to approximate the posterior
Hessian can therefore be used to decompose the parameter space into a
likelihood-informed subspace (LIS) spanned by the columns of $\mat{V}_r$ and an
orthogonal complementary space (CS).   %
Within each subspace, a standard Metropolis-Hastings kernel is employed.   As
long as the kernel in the CS uses a dimension-independent proposal (typically
pCN), then the DILI sampler will remain dimension-independent.   Unlike the
original implementation described in~\citet{CuiLawMarzouk16}, the MUQ
implementation does not use a whitening transform and thus avoids computing any
symmetric decomposition of the prior covariance.   In general, the Hessian used
in \eqref{eq:eigenproblem_1} can be adapted to capture more correlation
structure.  However, we did not find this necessary in the numerical
experiments below.

\begin{figure}
  \centering
  \begin{tikzpicture}
\tikzmath{\Ly=1; \Lx=1.618*\Ly;}


\draw[thick] (-0.1*\Lx,-0.2*\Ly) rectangle (3.2*\Lx, 1.5*\Ly);
\node[anchor=north west] at (-0.1*\Lx, 1.5*\Ly) {\footnotesize Metropolis-Hastings Algorithms};

\draw [fill=chain!40] (0,\Ly) -- (0,0) -- (\Lx,0) -- (\Lx,0.325*\Ly) -- (1.125*\Lx, 0.5*\Ly) -- (\Lx,0.675*\Ly) -- (\Lx, \Ly) -- cycle;
\node[anchor=center] at (0.5*\Lx, 0.5*\Ly) {\scriptsize  Chain};

\begin{scope}[xshift=1.05*\Lx cm]
\draw [fill=kernel!40] (0,\Ly) -- (0,0.675*\Ly) -- (0.125*\Lx,0.5*\Ly) -- (0,0.326*\Ly) -- (0,0) -- (\Lx,0) -- (\Lx, 0.325*\Ly) -- (1.125*\Lx, 0.325*\Ly) -- (1.125*\Lx, 0.675*\Ly) -- (\Lx, 0.675*\Ly) -- (\Lx, \Ly) -- cycle;
\node[anchor=center] at (0.55*\Lx, 0.5*\Ly) {\scriptsize MH Kernel};
\end{scope}

\begin{scope}[xshift=2.1*\Lx cm]
\draw [fill=proposal!40] (0,\Ly) -- (0,0.675*\Ly) -- (0.125*\Lx,0.675*\Ly) -- (0.125*\Lx,0.325*\Ly) -- (0,0.325*\Ly) -- (0,0) -- (\Lx,0) -- (\Lx, \Ly) -- cycle;
\node[anchor=center] at (0.55*\Lx, 0.5*\Ly) {\scriptsize Proposal};
\end{scope}

\begin{scope}[xshift=4.5*\Lx cm, yshift=-\Ly cm]

\draw[thick] (-0.1*\Lx,-1.2*\Ly) rectangle (3.2*\Lx, 2.5*\Ly);
\node[anchor=north west] at (-0.1*\Lx, 2.5*\Ly) {\footnotesize Delayed Rejection Algorithms};

\draw [fill=chain!40] (0,\Ly) -- (0,0) -- (\Lx,0) -- (\Lx,0.325*\Ly) -- (1.125*\Lx, 0.5*\Ly) -- (\Lx,0.675*\Ly) -- (\Lx, \Ly) -- cycle;
\node[anchor=center] at (0.5*\Lx, 0.5*\Ly) {\scriptsize  Chain};

\begin{scope}[xshift=1.05*\Lx cm]
\draw [fill=kernel!40] (0,2*\Ly) -- (0,0.675*\Ly) -- (0.125*\Lx,0.5*\Ly) -- (0,0.326*\Ly) -- (0,-\Ly) -- (\Lx,-\Ly) -- (\Lx, 0.325*\Ly-\Ly) -- (1.125*\Lx, 0.325*\Ly-\Ly) -- (1.125*\Lx, 0.675*\Ly-\Ly) -- (\Lx, 0.675*\Ly-\Ly) -- (\Lx, 0.325*\Ly) -- (1.125*\Lx, 0.325*\Ly) -- (1.125*\Lx, 0.675*\Ly) -- (\Lx, 0.675*\Ly) -- (\Lx, 1.325*\Ly) -- (1.125*\Lx, 1.325*\Ly) -- (1.125*\Lx, 1.675*\Ly) -- (\Lx, 1.675*\Ly) -- (\Lx, 2*\Ly) -- cycle;
\node[anchor=center] at (0.55*\Lx, 0.5*\Ly) {\scriptsize DR Kernel};
\end{scope}

\begin{scope}[xshift=2.1*\Lx cm,yshift=\Ly cm]
\draw [fill=proposal!40] (0,\Ly) -- (0,0.675*\Ly) -- (0.125*\Lx,0.675*\Ly) -- (0.125*\Lx,0.325*\Ly) -- (0,0.325*\Ly) -- (0,0) -- (\Lx,0) -- (\Lx, \Ly) -- cycle;
\node[anchor=center] at (0.55*\Lx, 0.5*\Ly) {\scriptsize Proposal 1};
\end{scope}

\begin{scope}[xshift=2.1*\Lx cm]
\node[anchor=center] at (0.55*\Lx, 0.55*\Ly) {$\vdots$};
\end{scope}

\begin{scope}[xshift=2.1*\Lx cm,yshift=-\Ly cm]
\draw [fill=proposal!40] (0,\Ly) -- (0,0.675*\Ly) -- (0.125*\Lx,0.675*\Ly) -- (0.125*\Lx,0.325*\Ly) -- (0,0.325*\Ly) -- (0,0) -- (\Lx,0) -- (\Lx, \Ly) -- cycle;
\node[anchor=center] at (0.55*\Lx, 0.5*\Ly) {\scriptsize Proposal $J$};
\end{scope}

\end{scope}

\begin{scope}[yshift=-2.5*\Ly cm]

\draw[thick] (-0.1*\Lx,-0.75*\Ly) rectangle (4.25*\Lx, 2.05*\Ly);
\node[anchor=north west] at (-0.1*\Lx, 2.05*\Ly) {\footnotesize DILI Algorithms};

\draw [fill=chain!40] (0,\Ly) -- (0,0) -- (\Lx,0) -- (\Lx,0.325*\Ly) -- (1.125*\Lx, 0.5*\Ly) -- (\Lx,0.675*\Ly) -- (\Lx, \Ly) -- cycle;
\node[anchor=center] at (0.5*\Lx, 0.5*\Ly) {\scriptsize  Chain};

\begin{scope}[xshift=1.05*\Lx cm]
\draw [fill=kernel!40] (0,1.55*\Ly) -- (0,0.675*\Ly) -- (0.125*\Lx,0.5*\Ly) -- (0,0.326*\Ly) -- (0,-0.55*\Ly) -- (\Lx,-0.55*\Ly) -- (\Lx, 0.325*\Ly-0.55*\Ly) -- (1.125*\Lx, -0.05*\Ly) -- (\Lx, 0.675*\Ly-0.55*\Ly) -- (\Lx, 0.325*\Ly+0.55*\Ly) -- (1.125*\Lx, 1.05*\Ly) -- (\Lx, 0.675*\Ly+0.55*\Ly) -- (\Lx, 1.55*\Ly) -- cycle;
\node[anchor=center] at (0.55*\Lx, 0.5*\Ly) {\scriptsize DILI Kernel};
\end{scope}

\begin{scope}[xshift=2.1*\Lx cm, yshift=0.55*\Ly cm]
\draw [fill=kernel!40] (0,\Ly) -- (0,0.675*\Ly) -- (0.125*\Lx,0.5*\Ly) -- (0,0.326*\Ly) -- (0,0) -- (\Lx,0) -- (\Lx, 0.325*\Ly) -- (1.125*\Lx, 0.325*\Ly) -- (1.125*\Lx, 0.675*\Ly) -- (\Lx, 0.675*\Ly) -- (\Lx, \Ly) -- cycle;
\node[anchor=center] at (0.55*\Lx, 0.5*\Ly) {\scriptsize MH Kernel};
\end{scope}

\begin{scope}[xshift=3.15*\Lx cm, yshift=0.55*\Ly cm]
\draw [fill=proposal!40] (0,\Ly) -- (0,0.675*\Ly) -- (0.125*\Lx,0.675*\Ly) -- (0.125*\Lx,0.325*\Ly) -- (0,0.325*\Ly) -- (0,0) -- (\Lx,0) -- (\Lx, \Ly) -- cycle;
\node[anchor=center] at (0.55*\Lx, 0.5*\Ly) {\scriptsize LIS Proposal};
\end{scope}

\begin{scope}[xshift=2.1*\Lx cm, yshift=-0.55*\Ly cm]
\draw [fill=kernel!40] (0,\Ly) -- (0,0.675*\Ly) -- (0.125*\Lx,0.5*\Ly) -- (0,0.326*\Ly) -- (0,0) -- (\Lx,0) -- (\Lx, 0.325*\Ly) -- (1.125*\Lx, 0.325*\Ly) -- (1.125*\Lx, 0.675*\Ly) -- (\Lx, 0.675*\Ly) -- (\Lx, \Ly) -- cycle;
\node[anchor=center] at (0.55*\Lx, 0.5*\Ly) {\scriptsize MH Kernel};
\end{scope}

\begin{scope}[xshift=3.15*\Lx cm, yshift=-0.55*\Ly cm]
\draw [fill=proposal!40] (0,\Ly) -- (0,0.675*\Ly) -- (0.125*\Lx,0.675*\Ly) -- (0.125*\Lx,0.325*\Ly) -- (0,0.325*\Ly) -- (0,0) -- (\Lx,0) -- (\Lx, \Ly) -- cycle;
\node[anchor=center] at (0.55*\Lx, 0.5*\Ly) {\scriptsize CS Proposal };
\end{scope}

%

\end{scope}

\end{tikzpicture}
  \caption{The flexible framework of hIPPYlib-MUQ allows many different
  combinations of transition kernels and proposal distributions to be employed.
Note that each kernel can
interact with any proposal distribution, which enables many different MCMC
algorithms to be constructed from the same basic components. }
  \label{fig:mcmc_puzzle}
\end{figure}
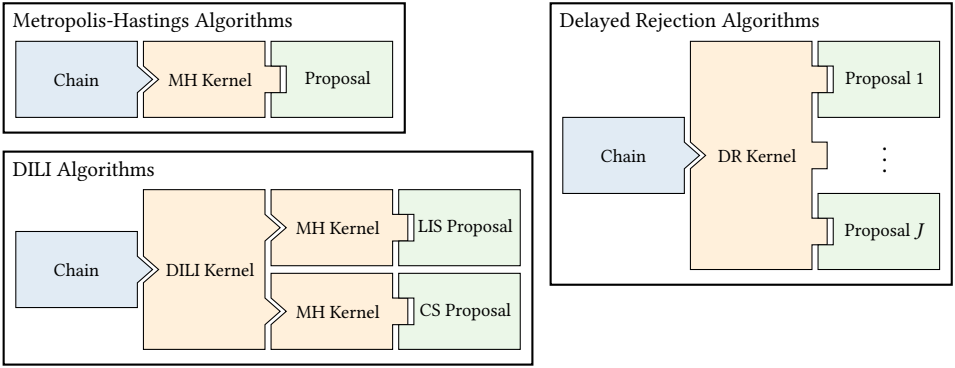

\paragraph{Assembling an MCMC Algorithm}
It is possible to combine nearly any of the proposals and kernels described
above, resulting in myriad possible MCMC algorithms.    As suggested in Figure
\ref{fig:mcmc_puzzle}, there are three fundamental building blocks to an MCMC
algorithm.   The chain keeps track of previous points and allows computing
Monte Carlo estimates.   The kernel defines a mechanism for sampling the next state
$\vec{\ipar}_{i+1}$ given the value of the current state $\vec{\ipar}_i$ and
one or more proposal distributions.   The proposal defines a position specific probability distribution that can be easily sampled and has a density that can be efficiently evaluated.
We mimic these abstract interfaces in our software design to define and test a large
number of kernel-proposal combinations.

\subsection{MCMC diagnostics}
\label{sec:mcmc_diagnostics}
Two questions naturally arise when analyzing a length $I$ Markov chain
$[\mathbf{m}_1,\ldots,\mathbf{m}_{I}]$ produced by MCMC.  First, has the chain
converged to the stationary distribution?  Second, what is the statistical
efficiency of the chain, that is, how many independent samples does the chain
have that actually contribute to the accuracy of Monte Carlo estimators?
Most theoretical guarantees are asymptotic, and it is important to
quantitatively answer these questions when employing finite-length MCMC chains.
Based on these considerations, this section describes the diagnostics
implemented in hIPPYlib-MUQ to check the convergence and statistical efficiency
of high-dimensional MCMC chains.

\subsubsection{Assessing Convergence}
To assess convergence, we compute two different asymptotically unbiased
estimators of the posterior covariance: one that is an overestimate for finite
$I$ and one that is an underestimate for finite $I$.  As the ratio of these two
estimates approaches one, we can be confident that the MCMC chain has converged
(see e.g.,~\citet{Brooks1998,Gelman2004, Vehtari2020}).

The estimates are based on running $J$ independent chains starting from
randomly chosen points that are more disperse than the posterior distribution
$\mupost$, where we define a ``disperse'' distribution as one that has a larger
covariance than $\mupost$.  Each chain has the same length $I$.

Letting $\mathbf{m}_{ij}$ be the $i$th MCMC sample in chain $j$, we define the
within-sequence covariance matrix $\mathbf{W}$ and the between-sequence
covariance matrix $\mathbf{B}$ as
\begin{alignat}{3}
  \mathbf{W} &= \frac{1}{J(I-1)} \sum_{j=1}^J \sum_{i=1}^I
  (\mathbf{m}_{ij} - \bar{\mathbf{m}}_{.j})
  (\mathbf{m}_{ij} - \bar{\mathbf{m}}_{.j})^T;
  \quad
  &&\bar{\mathbf{m}}_{.j} &&= \frac{1}{I} \sum_{i=1}^I \mathbf{m}_{ij},
  \label{eq:variance_within}\\
  \mathbf{B} &= \frac{I}{J-1} \sum_{j=1}^J
  (\bar{\mathbf{m}}_{.j} - \bar{\mathbf{m}}_{..})
  (\bar{\mathbf{m}}_{.j} - \bar{\mathbf{m}}_{..})^T;
  &&\bar{\mathbf{m}}_{..} &&= \frac{1}{J} \sum_{j=1}^J \bar{\mathbf{m}}_{.j}.
  \label{eq:variancce_across}
\end{alignat}
As pointed out in~\citet{Brooks1998}, $\mathbf{W}$ and $\mathbf{B}$ can be
combined to produce an estimate $\widehat{\mathbf{V}}$ of the posterior covariance that
takes the form
\begin{equation}
  \widehat{\mathbf{V}} = \frac{I-1}{I} \mathbf{W} + \frac{J+1}{JI}\mathbf{B}.\label{eq:vhat}
\end{equation}
The overdispersion of the initial points in each chain causes
$\widehat{\mathbf{V}}$ to
overestimate the posterior covariance for finite $I$.  On the other hand, the
average within-chain covariance $\mathbf{W}$ will tend to underestimate the covariance
because the chains have not explored the entire parameter space.   Comparing
$\mathbf{W}$ and $\widehat{\mathbf{V}}$ thus provides a way of assessing convergence.

The $\hat{R}$ statistic of~\citet{Gelman2004} and~\citet{Vehtari2020} is a
common way of comparing $\mathbf{W}$ and $\widehat{\mathbf{V}}$.  It uses the ratio of the
diagonal component of $\widehat{\mathbf{V}}$ and $\mathbf{W}$ to construct a componentwise
convergence diagnostic.  For high dimensional problems,
it is more natural to consider a multivariate convergence diagnostic.
We will therefore employ the multivariate potential scale reduction
factor (MPSRF) of~\citet{Brooks1998}, which is a natural extension of
the componentwise $\hat{R}$ statistic.  The MPSRF is defined by
\begin{equation}
\begin{aligned}
  \text{MPSRF} &= \sqrt{\max_a \frac{a^T \widehat{\mathbf{V}} a}{a^T \mathbf{W} a}} %
    & =  \sqrt{\frac{I-1}{I} + \frac{J+1}{JI} \lambda_{\text{max}}} \, ,
\end{aligned}
\end{equation}
where $\lambda_{\text{max}}$ is the largest eigenvalue satisfying the
generalized eigenvalue problem $\mathbf{B} \vec{v} = \lambda \mathbf{W} \vec{v}$.

Note that $\text{MPSRF}\geq 1$ when $\lambda_{\text{max}}>1$, which occurs when the chains have overdispersed starting points that cause the inter-chain variance $\mathbf{B}$ to be larger than the within-chain variance $\mathbf{W}$. When the MPSRF approaches 1,
the variance within each sequence approaches the variance across sequences,
thus indicating that each individual chain has converged to the
target distribution.  The literature contains several recommendations for 
values of MPSRF that indicate convergence; for example, \citet{gelman1992inference}
suggest the commonly used value $\text{MPSRF}<1.1$ while \citet{Vehtari2020} 
argues for the more conservative threshold $\text{MPSRF}<1.01$.

\subsubsection{Statistical Efficiency}
The samples in an MCMC chain are generally correlated, which increases the
variance of Monte Carlo estimators constructed with MCMC samples. For a
quantity of interest $\qoi(\vec{\ipar})$, the effective sample size (ESS) of a
Markov chain is defined as the number of independent samples of the posterior
that would be needed to estimate $\mathbb{E}[\qoi]$ with the same statistical
accuracy as an estimate from the Markov chain.  The ESS is therefore a measure
of how much information is contained in the MCMC chain.  In this work, it is
commonly assumed that the ESS is derived for estimators of the posterior mean,
i.e., $\mathbb{E}[\qoi]=\mathbb{E}[\vec{\ipar}]$, and here we derive the ESS under
this common assumption.  Other ESS variants, like those described by \citet{Vehtari2020}, are 
more suitable for problems involving tail probabilities, but the implementation 
of these methods in hIPPYlib-MUQ is left to future work.

There are several ways of estimating the ESS. For instance, spectral approaches use the
integrated autocorrelation of the MCMC chain to estimate the effective sample
size (see e.g., \citet{Gelman2004,wolff2004monte}).  Other common methods use
the statistics of small sample batches (see e.g.,
\citet{flegal2010batch,vats2019multivariate}).  MUQ provides implementations of
both spectral and batch methods.   Here we focus on the spectral formulation of
ESS however, because it gives additional insight into the structure of MCMC
chains.  The ESS for component $k$ of $\vec{\ipar}$ is defined by

\begin{equation}
    \label{eqn:ess}
    \text{ESS}_k = \frac{JI}{1 + 2 \sum_{t=1}^\infty \rho_{kt}},
\end{equation}
where $\rho_{kt}$ is the autocorrelation function of component $k$ in the MCMC chain at lag $t$.
Here, the autocorrelation function $\rho_{kt}$ is estimated by the following
formula~\citep{Gelman2004}:
\begin{equation}
    \label{eqn:acf}
  \rho_{kt} \approx \hat{\rho}_{kt} = 1 - \frac{v_{kt}}{2 \hat{V}_{kk}},
\end{equation}
where $\hat{V}_{kk}$ is the $k$th diagonal component of the posterior
covariance estimate defined in \eqref{eq:vhat} and $v_{kt}$ is the variogram
defined by
\begin{equation}
    v_{kt} = \frac{1}{J(I-t)} \sum_{j=1}^J \sum_{i=t+1}^I
    (m_{ij,k} - m_{(i-t)j,k})^2.
\end{equation}
In practice, $\hat{\rho}_{kt}$ is noisy for large values of $t$ and we truncate
the summation \eqref{eqn:ess} at some lag $t'$. Following common practice, we
choose $t'\ge 0$ to be the lag for which the sum successive autocorrelation
estimates $\hat{\rho}_{2t'} + \hat{\rho}_{2t' + 1}$ is
negative~\citep{Gelman2004}.

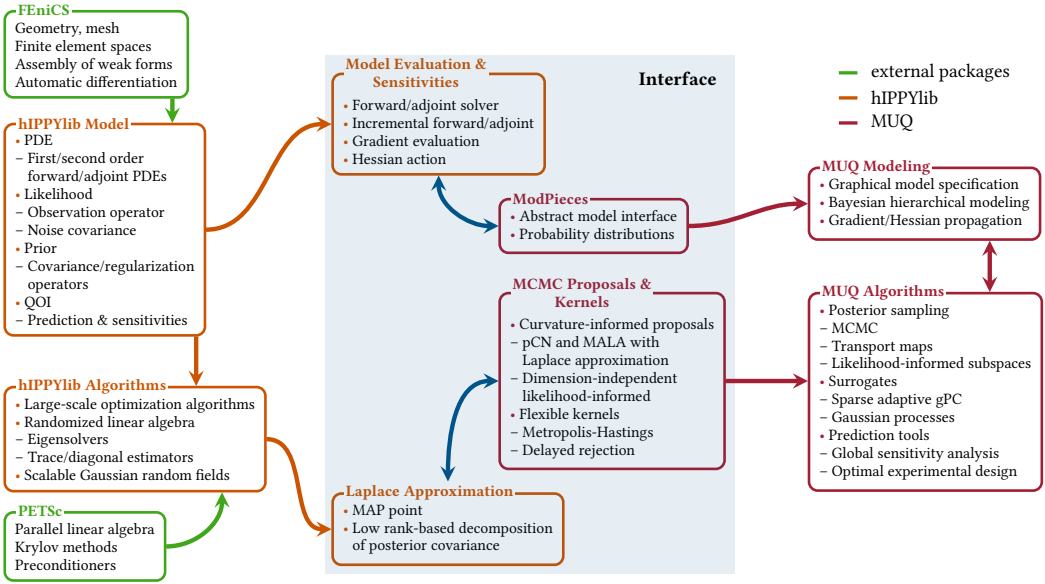
\begin{figure}[ht]
  \centering
  \resizebox{\textwidth}{!}{\begin{tikzpicture}[->, >=stealth, line width=0.2em,
  every node/.style={font=\fontsize{9}{10.2}\selectfont, scale=0.8}]
  \def\mybox#1#2#3#4#5{
    \node[draw=#1, very thick, rounded corners,
    anchor=north west, inner sep=6, align=left]
    (#2) at #3 {#4};
    \node[#1, inner sep=1, fill=white,
    anchor=west, right=0.5em, align=center] at (#2.north west) {#5};
  }

  \mybox{darkgreen}{fenics}{(0,9)}
  {Geometry, mesh\\
    Finite element spaces\\
    Assembly of weak forms\\
  Automatic differentiation}
  {\textbf{FEniCS}}

  \mybox{utorange}{hlmodel}{(0,7.2)}
  {\textcolor{utorange}{\textbullet} PDE\\
    -- First/second order\\
    \phantom{--} forward/adjoint PDEs\\
    \textcolor{utorange}{\textbullet} Likelihood\\
    -- Observation operator\\
    -- Noise covariance\\
    \textcolor{utorange}{\textbullet} Prior\\
    -- Covariance/regularization\\
    \phantom{--} operators\\
    \textcolor{utorange}{\textbullet} QOI\\
    -- Prediction \& sensitivities
  }
  {\textbf{hIPPYlib Model}}

  \mybox{utorange}{hlalgorithms}{(0,3)}
  {\textcolor{utorange}{\textbullet} Large-scale optimization algorithms\\
    \textcolor{utorange}{\textbullet} Randomized linear algebra\\
    -- Eigensolvers\\
    -- Trace/diagonal estimators\\
  \textcolor{utorange}{\textbullet} Scalable Gaussian random fields}
  {\textbf{hIPPYlib Algorithms}}

  \mybox{darkgreen}{petsc}{(0,1)}
  {Parallel linear algebra\\ Krylov methods\\ Preconditioners}
  {\textbf{PETSc}}

  \mybox{utorange}{hloutputs1}{(5.3,8)}
  {\phantom{\tiny phantom}\\
    \textcolor{utorange}{\textbullet} Forward/adjoint solver\\
    \textcolor{utorange}{\textbullet} Incremental forward/adjoint\\
    \textcolor{utorange}{\textbullet} Gradient evaluation\\
  \textcolor{utorange}{\textbullet} Hessian action}
  {\textbf{Model Evaluation \&}\\ \textbf{Sensitivities}}

  \mybox{utorange}{hloutputs2}{(5.3,1.3)}
  {
    \textcolor{utorange}{\textbullet} MAP point\\
    \textcolor{utorange}{\textbullet} Low rank-based decomposition\\
  \phantom{\textcolor{utorange}{\textbullet}} of posterior covariance}
  {\textbf{Laplace Approximation}}

  \mybox{MitRed}{muqmodpieces}{(8,6)}
  {\textcolor{MitRed}{\textbullet} Abstract model interface\\
  \textcolor{MitRed}{\textbullet} Probability distributions}{\textbf{ModPieces}}

  \mybox{MitRed}{muqmcmc}{(8,4.5)}
  {\phantom{\tiny phantom}\\
    \textcolor{MitRed}{\textbullet} Curvature-informed proposals\\
    -- pCN and MALA with\\
    \phantom{--} Laplace approximation\\
    -- Dimension-independent\\
    \phantom{--} likelihood-informed\\
    \textcolor{MitRed}{\textbullet} Flexible kernels\\
    -- Metropolis-Hastings\\
    -- Delayed rejection
  }
  {\textbf{MCMC Proposals \&}\\ \textbf{Kernels}}

  \mybox{MitRed}{muqmodeling}{(13,6.5)}
  {\textcolor{MitRed}{\textbullet} Graphical model specification\\
    \textcolor{MitRed}{\textbullet} Bayesian hierarchical modeling\\
  \textcolor{MitRed}{\textbullet} Gradient/Hessian propagation}
  {\textbf{MUQ Modeling}}

  \mybox{MitRed}{muqalgorithms}{(13,4.5)}
  {\textcolor{MitRed}{\textbullet} Posterior sampling\\
    -- MCMC\\
    -- Transport maps\\
    -- Likelihood-informed subspaces\\
    \textcolor{MitRed}{\textbullet}
    Surrogates\\
    -- Sparse adaptive gPC\\
    -- Gaussian processes\\
    \textcolor{MitRed}{\textbullet}
    Prediction tools\\
    -- Global sensitivity analysis\\
    -- Optimal experimental design
  }
  {\textbf{MUQ Algorithms}}

  \fill[opacity=0.1, color=lakeblue]
  ([xshift=-8em, yshift=-4.7em] muqmcmc.south west) rectangle
  ([xshift=+2.1em, yshift=+6.5em] muqmodpieces.north east);
  \node[above] at ([xshift=-.5em, yshift=+4.8em] muqmodpieces.north east) {\large \textbf{Interface}};

  \path ([xshift=+3.5em] fenics.south) edge[darkgreen]
  ([xshift=+3.5em] fenics.south |- hlmodel.north)
  (petsc.east) edge[out=0, in=-90, darkgreen] ([xshift=+4em] hlalgorithms.south)
  ([xshift=+4.2em] hlmodel.south) edge[utorange]
  ([xshift=+4.2em] hlmodel.south |- hlalgorithms.north)
  (hlmodel.east) edge[out=0, in=180, utorange] (hloutputs1.west)
  (hlalgorithms.east) edge[out=0, in=180, utorange] (hloutputs2.west)
  (hloutputs1.south) edge[<->, out=-90, in=180, lakeblue] (muqmodpieces.west)
  ([yshift=+0.5em] hloutputs2.north) edge[<->, out=90, in=180, lakeblue] (muqmcmc.west)
  (muqmodpieces.east) edge[out=0, in=180, MitRed]
  (muqmodeling.west)
  (muqmcmc.east) edge[MitRed] (muqmcmc.east -| muqalgorithms.west)
  ([xshift=+3em]muqmodeling.south) edge[<->, MitRed]
  ([xshift=+3em]muqmodeling.south |- muqalgorithms.north);

  \draw[-,color=darkgreen, line width=1.5pt] (13.5,8) -- (13.8, 8);
  \node [color=black, anchor=west] at (13.9,8) {\large external packages};
  \draw[-,color=utorange, line width=1.5pt] (13.5,7.6) -- (13.8, 7.6);
  \node [color=black, anchor=west] at (13.9,7.6) {\large hIPPYlib};
  \draw[-,color=MitRed, line width=1.5pt] (13.5,7.2) -- (13.8, 7.2);
  \node [color=black, anchor=west] at (13.9,7.2) {\large MUQ};

\end{tikzpicture}}
  \caption{Description of the functionalities of hIPPYlib and MUQ and their
  interface. Orange and red boxes represent hIPPYlib and MUQ functionalities,
  respectively. Green boxes indicate external software libraries, FEniCS and PETSc, that provide parallel implementation of finite element discretizations and solvers. Arrows represent one-way or reciprocal interactions.}
  \label{fig:tikz/integration_diagram}
\end{figure}

\section{Software framework}
\label{sec:software_framework}

hIPPYlib-MUQ is a Python interface that
integrates two open source software libraries, hIPPYlib and MUQ, into a unique software
framework, allowing the user to implement state-of-the-art Bayesian inversion
algorithms in a seamless way.
We outline in Figure~\ref{fig:tikz/integration_diagram} the main
functionalities of hIPPYlib and MUQ and the integration of their complementary
components.

hIPPYlib, built on FEniCS and PETSc for the discretization and
solution of PDEs, provides Python implementations of scalable
adjoint-based algorithms for solving large-scale deterministic and
linearized Bayesian inverse problems governed by PDEs.  hIPPYlib model
component provides a collection of libraries by which users can
describe the forward PDE, the prior model, and the likelihood model in
the FEniCS form language \citep{LoggMardalWells12}.  hIPPYlib
algorithms component incorporates optimization algorithms, randomized
linear algebra, and scalable sampling of Gaussian fields that are
required to efficiently solve the deterministic and linearized
Bayesian inverse problems.  The hiPPYlib user manual can be found at
\citet{Villa2020}, which includes details on the software
installation, documentation, and tutorials.

MUQ is an easy-to-use software framework for defining and solving uncertainty
quantification problems.
MUQ modeling tools allow users to easily and flexibly construct complicated
models, including Bayesian hierarchical models, in a semi-intrusive way that
enables efficient gradient and Hessian evaluations.
MUQ also implements a variety of advanced uncertainty quantification
techniques including MCMC sampling methods, surrogate models (e.g., Gaussian
processes), and prediction tools (e.g., global sensitivity analysis).
We refer the reader to \href{http://muq.mit.edu}{muq.mit.edu}
for a detailed description of the software with installation instructions and step-by-step tutorials. 

We note that a significant synergistic effect can be obtained by making use of
complementary aspects of these two software libraries: hIPPYlib outputs such as
the gradient and Hessian evaluations and the Laplace approximation, and MUQ's
advanced MCMC sampling modules and flexible modeling capabilities.
In the integrated framework, hIPPYlib is used to define
the forward model, the prior, and the likelihood, to compute the maximum a
posteriori (MAP) point, and to construct a Laplace approximation of
the posterior distribution based on approximations of the posterior covariance
as a low-rank update of the prior~\citep{Bui-ThanhGhattasMartinEtAl13a}.  MUQ
is employed to exploit advanced MCMC methods to fully characterize the
posterior distribution in non-Gaussian/nonlinear settings.  hIPPYlib-MUQ offers
a set of wrappers that encapsulate the functionality of hIPPYlib in a way that
various features of hIPPYlib can be accessed by MUQ.  A key aspect of
hIPPYlib-MUQ is that it enables the use of curvature-informed MCMC methods,
which is crucial for efficient and scalable exploration of the posterior
distribution for large-scale Bayesian inverse problems. 

In the context of hIPPYlib-MUQ, hIPPYlib provides tools for (i) automatically
implementing adjoint gradients and Hessian actions, (ii) efficiently sampling
Gaussian Markov Random fields (GMRF), and (iii) constructing Laplace
approximations with low-rank Hessians.
On the other hand, MUQ provides two important capabilities: (i) an
abstract graphical modeling framework that provides an interface for
implementing model components (e.g., prior distributions or forward models)
and enables multiple components of a model or inverse problem to be easily
composed and differentiated, and (ii) a suite of MCMC algorithms, including
curvature-informed and discretization-invariant methods.   
The adjoint techniques of hIPPYlib enable MUQ to efficiently compute gradients and Hessian
actions of a graphical model with PDE-based components.  The efficient GMRF
sampling and low-rank Laplace approximations accelerate MUQ's
discretization-invariant MCMC techniques, which use these hIPPYlib tools within
the MCMC proposal.    The details of our object orientated approach for
seamlessly blending these MUQ and hIPPYlib tools are provided below.

\begin{figure}[h]
  \centering
  \resizebox{0.5\textwidth}{!}{\begin{tikzpicture}[->, >=stealth, line width=0.2em,
  every node/.style={font=\fontsize{9}{10.2}\selectfont, scale=0.8},
  mybox/.style={very thick, rounded corners, inner sep=4}]

  %
  %
  \node [] (LaplaceGaussian)
    {{\color{lakeblue} \textbf{LaplaceGaussian}}};
  \node [mybox, draw=utorange, below=0.1cm of LaplaceGaussian] (LaplacianPrior)
    {
      {\color{utorange} \textbf{LaplacianPrior}}
    };
  \node [mybox, draw=lakeblue, inner sep=8, fit={(LaplaceGaussian)
    (LaplacianPrior)}] (LaplacianPriorBox) {};

  \node [below=0.2cm of LaplacianPriorBox] (BiLaplaceGaussian)
    {{\color{lakeblue} \textbf{BiLaplaceGaussian}}};
  \node [mybox, draw=utorange, below=0.1cm of BiLaplaceGaussian] (BiLaplacianPrior)
    {
      {\color{utorange} \textbf{BiLaplacianPrior}}
    };
  \node [mybox, draw=lakeblue, inner sep=8, fit={(BiLaplaceGaussian)
    (BiLaplacianPrior)}] (BiLaplacianPriorBox) {};

  \node [below=0.3cm of BiLaplacianPriorBox] (LAPosteriorGaussian)
    {{\color{lakeblue} \textbf{LAPosteriorGaussian}}};
  \node [mybox, draw=utorange, below=0.1cm of LAPosteriorGaussian, inner
    sep=3.5] (GaussianLRPosterior)
    {
      {\color{utorange} \textbf{GaussianLRPosterior}}
    };
  \node [mybox, draw=lakeblue, inner sep=11, fit={(LAPosteriorGaussian)
    (GaussianLRPosterior)}] (LAPosteriorGaussianBox) {};

  \node [below=0.4cm of LAPosteriorGaussianBox] (Param2LogLikelihood)
    {{\color{lakeblue} \textbf{Param2LogLikelihood}}};
  \node [mybox, draw=utorange, below=0.1cm of Param2LogLikelihood] (Model)
    {
      {\color{utorange} \textbf{Model}}
    };
  \node [mybox, draw=lakeblue, inner sep=8, fit={(Param2LogLikelihood)
    (Model)}] (Param2LogLikelihoodBox) {};

  %
  %
  \node [mybox, draw=MitRed, MitRed, right= of BiLaplacianPriorBox, inner sep=5] (GaussianBase)
    {\textbf{GaussianBase}};

  \node [mybox, draw=MitRed, MitRed, below= of GaussianBase, inner sep=5] (Density)
    {\textbf{Density}};

  \node [mybox, draw=MitRed, MitRed, right= of GaussianBase, inner sep=5] (Distribution)
    {\textbf{Distribution}};

  \node [mybox, draw=MitRed, MitRed, inner sep=5] (ModPiece) at
    (Param2LogLikelihoodBox -| Distribution) {\textbf{ModPiece}};

  \path [very thick] (LaplacianPriorBox.east) edge [dashed,out=0,in=180]
    ([yshift=0.1cm] GaussianBase.west)
    (BiLaplacianPriorBox.east) edge [dashed,out=0,in=180] (GaussianBase.west)
    (LAPosteriorGaussianBox.east) edge [dashed,out=0,in=180] ([yshift=-0.1cm] GaussianBase.west);

  \path [very thick] (Param2LogLikelihoodBox.east) edge [dashed,out=0,in=180]
    (ModPiece.west);

  \path [very thick] (GaussianBase.east) edge [dashed,out=0,in=180]
    (Distribution.west)
    (Density.east) edge [dashed,out=0,in=180] ([yshift=-0.1cm] Distribution.west)
    (Density.east) edge [dashed,out=0,in=180] ([yshift=0.1cm] ModPiece.west);
\end{tikzpicture}}
  \caption{Class hierarchy for hIPPYlib-MUQ framework. Classes of hIPPYlib, MUQ, and
  the interface are colored in orange, red, and blue, respectively. Dashed arrows
  represent inheritance relationship between two classes: the arrowhead
  attaches to the super-class and the other attaches to the sub-class.
  }
  \label{fig:tikz/class_hierarchy}
\end{figure}
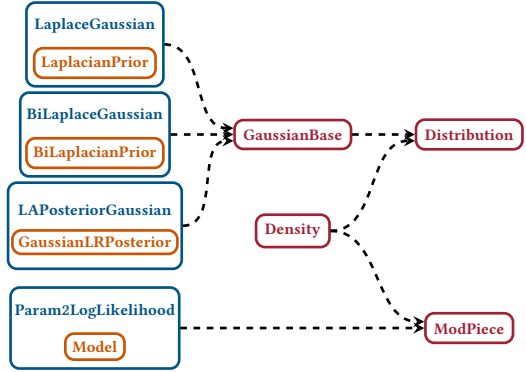

\begin{figure}[h]
  \centering
  \resizebox{\textwidth}{!}{\begin{tikzpicture}[->, >=stealth, line width=0.2em,
  every node/.style={font=\fontsize{9}{10.2}\selectfont, scale=0.8, align=center},
  mybox/.style={very thick, rounded corners, inner sep=6}]

  \node [align=left,draw,thick] (code_snippet)
  {
    \begin{lstlisting}[language=Python,flexiblecolumns=true,basicstyle=\sffamily]
# Example code snippet
import muq.Modeling as mm
import hippylib2muq as hm

# ... Use hIPPYlib to define prior and model variables

# Convert hiPPYlib components to MUQ components
prior_density = hm.BiLaplaceGaussian(prior).AsDensity()
likelihood = hm.Param2LogLikelihood(model)

# Add all of the components to the graph
graph = mm.WorkGraph()
graph.AddNode(mm.IdentityOperator(dim), 'Parameter')
graph.AddNode(prior_density, 'Prior')
graph.AddNode(likelihood, 'Likelihood')
graph.AddNode(mm.DensityProduct(2), 'Posterior')

# Define right branch:  Parameter->Prior->Posterior
graph.AddEdge('Parameter', 0, 'Prior', 0)
graph.AddEdge('Prior', 0, 'Posterior', 0)

# Define left branch:  Parameter->Likelihood->Posterior
graph.AddEdge('Parameter', 0, 'Likelihood', 0)
graph.AddEdge('Likelihood', 0, 'Posterior', 1)
    \end{lstlisting}
  };

  %
  %
  \coordinate [left=4.0cm of code_snippet] (model);

  \node [mybox, draw=black, above=1.0cm of model.north]
    (parameter) {Parameter\\ ({\color{MitRed} IdentityOperator})};
  \node [mybox, draw=black, left=0.5cm of model] (likelihood) {Likelihood\\
    ({\color{lakeblue} Param2LogLikelihood})};
  \node [mybox, draw=black, right=0.5cm of model] (prior) {Prior\\
    ({\color{lakeblue} e.g., LaplaceGaussian})};
  \node [mybox, draw=black, below=1.0cm of model.south] (posterior)
    {Posterior\\ ({\color{MitRed} DensityProduct})};

  \node [mybox, draw=black, above=0.5cm of parameter] (input) {Input};
  \node [mybox, draw=black, below=0.5cm of posterior] (output) {Output};

  \path [thick] (parameter.south) edge[out=-30, in=150] (prior.north)
    (parameter.south) edge[out=-150, in=30] (likelihood.north)
    (prior.south) edge[out=-150, in=30] (posterior.north)
    (likelihood.south) edge[out=-30, in=150] (posterior.north)
    (input.south) edge (parameter.north)
    (posterior.south) edge (output.north);

\end{tikzpicture}}
  \caption{Graphical description of Bayesian posterior modeling using hIPPYlib-MUQ
  software framework (left) and an example code snippet (right). In the left
  figure, class names of MUQ and the interface are colored in red and blue, respectively.
  MUQ's \texttt{WorkGraph} class provides a way to combine all the Bayesian
  posterior model components by its member functions \texttt{AddNode} and
  \texttt{AddEdge}.
  MUQ's \texttt{IdentityOperator} class identifies input parameters and
  the input argument \texttt{dim} represents the parameter dimension.
  MUQ's \texttt{DensityProduct} class defines the product of prior and likelihood
  densities and the input argument 2 means the number of input densities.
  The second and forth arguments of the function \texttt{AddEdge} mean the
  output index of the first argument and the input index of the third argument,
  respectively.
  For example, \texttt{graph.AddEdge('Parameter', 0, 'Prior', 0)} means that
  the output of 'Parameter' indexed 0 is connected to the input of 'Prior' indexed 0.
}
  \label{fig:tikz/graph_model}
\end{figure}
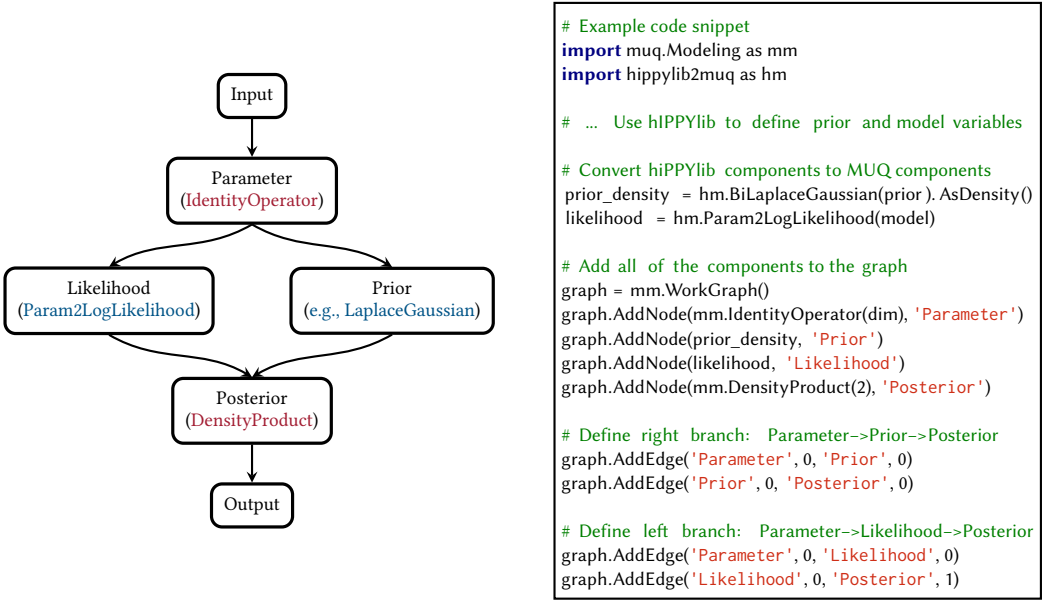

Figure~\ref{fig:tikz/class_hierarchy} gives an overview of the Python
classes implemented by the hIPPYlib-MUQ interface. %
Inherited from MUQ classes, the interface classes wrap the hIPPYlib
functionalities needed to achieve curvature-informed MCMC sampling methods.
These include:
\begin{enumerate}

  \item Prior Gaussian interface (\texttt{LaplaceGaussian} and
    \texttt{BiLaplaceGaussian}) to enable the use of hIPPYlib prior models
    (\texttt{LaplacianPrior} and \texttt{BiLaplacianPrior}) in MUQ probability
    distribution models (\texttt{GaussianBase}).
    The interface allows MUQ to use Gaussian prior with covariance as inverse
    of the $v$-th power of a Laplacian-like operator ($v = 1$ for
    \texttt{LaplaceGaussian} and $v=2$ for \texttt{BiLaplacianPrior}) and
    scalable sampling techniques for Gaussian random fields.

  \item Likelihood interface (\texttt{Param2LogLikelihood}) to
    incorporate hIPPYlib likelihood functions (\texttt{Model}) into MUQ
    Bayesian models (\texttt{ModPiece}). 
    The interface enables the MUQ model to efficiently perform the model evaluation (the
    parameter-to-observable map) including forward PDE solves and the adjoint-based computation of gradient
    and Hessian actions.

  \item Laplace approximation interface (\texttt{LAPosteriorGaussian}) to
    let MUQ MCMC modules get access to the Laplace approximation of the posterior distribution
    generated by hIPPYlib (\texttt{GaussianLRPosterior}).
    This interface provides the MAP point and/or the low-rank approximation of
    the Hessian at the MAP point for MUQ MCMC proposals, which leads to
    efficient and scalable sampling of the posterior.

\end{enumerate}
These interface classes can then be used to form a Bayesian posterior model
governed by PDEs using MUQ's graphical modeling interface (\texttt{WorkGraph}) as
shown in Figure~\ref{fig:tikz/graph_model}, as well as to construct MCMC
proposals.

hIPPYlib-MUQ also implements the MCMC convergence diagnostics described in
Section \ref{sec:mcmc_diagnostics}. These include the
potential scale reduction factor and its extension to multivariate parameter
cases~\citep{Brooks1998}, the autocorrelation function, and the effective sample
size.   A detailed description of all classes and functionalities of hIPPYlib-MUQ can also be found at
  \href{https://hippylib2muq.readthedocs.io/en/latest/modules.html}
  {https://hippylib2muq.readthedocs.io/en/latest/modules.html.}

\section{Numerical Illustration}
\label{sec:numerical_illustration}
The objective of this section is to showcase applications of the integrated
software framework discussed in previous sections via a step-by-step
implementation procedure.  We focus on comparing the performance of several
MCMC methods available in the software framework.  For the illustration we
first revisit the model problem considered in~\citet{VillaPetraGhattas21}, an
inverse problem of reconstructing the log-diffusion coefficient field in a
two-dimensional elliptic partial differential equation. We then consider a
nonlinear $p$-Poisson problem in three dimensions for which the parameter field
in a Robin boundary condition is inferred. In this section, we summarize the
Bayesian formulation of the example problems and present numerical results
obtained using the proposed software framework, hIPPYlib-MUQ version 0.2.0 that
builds on hIPPYlib version 3.1.0 with FEniCS version 2019 and MUQ version
0.3.0; see
\url{https://hippylib2muq.readthedocs.io/en/latest/installation.html}
for the software installation instruction and its dependencies. The
accompanying Jupyter notebook provides a detailed description of the
hIPPYlib-MUQ implementations; see
\url{https://hippylib2muq.readthedocs.io/en/latest/tutorial.html}.
Additional examples, including some with hierarchical models, can also be found
at \href{http://muq.mit.edu}{muq.mit.edu}.

\subsection{Coefficient field inversion in a two-dimensional Poisson linear PDE}

We first consider the coefficient field inversion in a Poisson partial
differential equation given pointwise noisy state measurements.  We begin by
describing the forward model setup and quantity of interest, followed by the
definition of the prior and the likelihood distributions.  Next, we present the
Laplace approximation of the posterior and apply several MCMC methods to
characterize the posterior distribution, as well as the predictive posterior
distribution of a scalar quantity of interest. Finally, the scalability of the
proposed methods with respect to the parameter dimension is then assessed in a
mesh refinement study. 

\subsubsection{Forward model}
\label{sec:fwd_model_2d}
Let $\Omega \in \mathbb{R}^d (d = 2, 3)$ be an open bounded domain with
boundary $\partial \Omega = \partial \Omega_D \cup \partial \Omega_N, \partial
\Omega_D \cap \partial \Omega_N
= \emptyset$.
Given a realization of the uncertain parameter field $m$, the state variable
$u$ is governed by

\begin{alignat}{2}
  \label{eqn:forward_strong}
  -\nabla \cdot \left( e^m \nabla u\right) &= f &\quad&\text{in } \Omega
  \nonumber,\\
  u &= g &&\text{on } \partial \Omega_D, \\
  e^m \nabla u \cdot \mathbf{n} &= h &&\text{on } \partial \Omega_N,
  \nonumber
\end{alignat}
where $f$ is a volume source term, $g$ and $h$ are the prescribed Dirichlet and
Neumann boundary data, respectively, and $\mathbf{n}$ is the outward unit normal
vector.

The weak form of \eqref{eqn:forward_strong} reads as follows: Find $u
\in \mathcal{V}_g$ such that
\begin{equation}
  \label{eqn:forward_weak}
  \langle e^m \nabla u, \nabla p \rangle = \langle f,p \rangle +
  \langle h, p \rangle_{\partial \Omega_N}
  \quad \forall p \in \mathcal{V}_0,
\end{equation}
where
\begin{align}
  \mathcal{V}_g &= \left\{v \in H^1(\Omega)|
  v = g \ \text{on} \ \partial \Omega_D \right\},
  \nonumber \\
  \mathcal{V}_0 &= \left\{v \in H^1(\Omega)|
  v = 0 \ \text{on} \ \partial \Omega_D \right\}.
\end{align}
Above, we denote the $L^2$-inner product over $\Omega$ by $\langle \cdot, \cdot
\rangle$ and that over $\partial \Omega_N$ by $\langle \cdot, \cdot
\rangle_{\partial \Omega_N}$.

As a quantity of interest, the logarithm of the normal flux through the bottom boundary
$\partial \Omega_b \subset \partial \Omega_D$ is considered. Specifically, we
define the quantity of interest $\qoi (m)$ as
\begin{equation}
  \label{eqn:qoi}
  \qoi(m) = \ln \left\{ - \int_{\partial \Omega_b} e^m \nabla u \cdot \mathbf{n} \, ds
  \right\}.
\end{equation}

In this example, we consider a unit square domain in $\mathbb{R}^2$ with
no source term ($f = 0$), no normal flux ($h = 0$) on the
left and right boundaries, and the Dirichlet condition imposed on the top
boundary ($g = 1$) and the bottom boundary ($g = 0$).

For the spatial discretization, we use quadratic finite elements for the state
variable (also for the adjoint variable) and linear finite elements for the
parameter variable.  
The computational domain is discretized using a regular mesh with 2,048
triangular elements. This leads to 4,225 and 1,089 degrees of freedom for the
state and parameter variables, respectively.  In the scalability results
presented in Section \ref{sec:mcmc_scalability_2d}, the mesh is then refined
with up to four levels of uniform refinement leading to 263,169 and 66,049
degrees of freedom for the state and parameter variables, respectively, on the
finest level.

\subsubsection{Prior model}
\label{sec:prior_model_2d}

As discussed in Section \ref{sec:bayesian_framework}, we choose the prior to be a
Gaussian distribution $\GM{\iparpr}{\Cprior}$ with $\Cprior = \Acal^{-2}$ where
$\Acal$ is a Laplacian-like operator given as
\begin{equation}
  \Acal m =
  \begin{cases}
    - \gamma \nabla \cdot (\Theta \nabla m) + \delta m
            & \text{in } \Omega,\\
            \Theta \nabla m \cdot \mathbf{n} + \beta m
            & \text{on } \partial \Omega.
  \end{cases}
\end{equation}
Here, $\beta \propto \sqrt{\gamma \delta}$ is the optimal Robin coefficient
introduced to alleviate undesirable boundary effects~\citep{DaonStadler18}, and
an anisotropic tensor $\Theta$ is of the form
\begin{equation}
  \Theta =
  \begin{bmatrix}
    \theta_1 \sin^2\alpha + \theta_2 \cos^2\alpha & (\theta_1 - \theta_2) \sin\alpha \cos\alpha
    \\ (\theta_1 - \theta_2) \sin\alpha \cos\alpha & \theta_1 \cos^2\alpha + \theta_2\sin^2\alpha
\end{bmatrix}.
\end{equation}

\begin{figure}[t]
  \centering
  \resizebox{\textwidth}{!}{\input{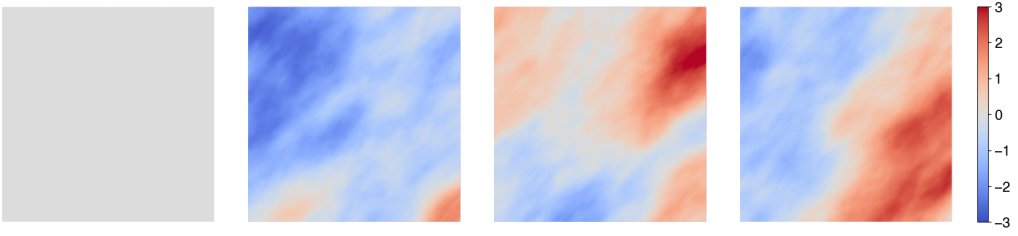}}
  \caption{Prior mean (leftmost) and three sample fields drawn from the prior
  distribution for the Poisson problem.}
  \label{fig:ex1_prior_samples}
\end{figure}
For this example we take $\gamma = 0.1, \delta = 0.5, \beta = \sqrt{\gamma \delta} / 1.42,
\theta_1 = 2.0, \theta_2 = 0.5$ and $\alpha = \pi / 4$.
Figure~\ref{fig:ex1_prior_samples} shows the prior mean $\iparpr$ and three
samples from the prior distribution.

\subsubsection{Observations with noise and the likelihood}
\label{sec:likelihood_model_2d}

\begin{figure}[t]
  \centering
  \resizebox{0.8\textwidth}{!}{\input{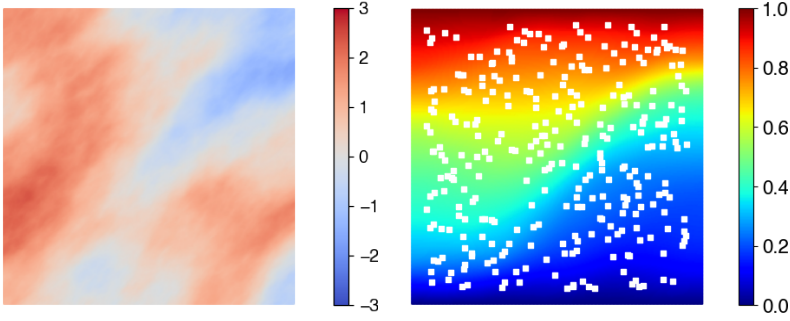}}
  \caption{True parameter field (left) and the corresponding state field
  (right) for the Poisson problem. The locations of the observation points are
  marked as white squares in
  the right figure.}
  \label{fig:ex1_mtrue_observation}
\end{figure}

We generate state observations at $l$ random locations uniformly distributed
over $[0.05, 0.95]^2$ by solving the forward problem on the
finest mesh with the true parameter field $\ipart$ (here a sample from the
prior is used) and then adding a random Gaussian noise to the resulting state
values; see Figure~\ref{fig:ex1_mtrue_observation}. The number of observations
$l$ is set to 300 for this example.
The vector of synthetic observations is given by
\begin{equation}
    \obs = \mathcal{B} u + \vec{\eta},
\end{equation}
where $\mathcal{B}$ is a linear observation operator, restricting the state
solution to the $l$ observation points.
The additive noise vector $\vec{\eta}$ has mutually independent components that
are normally distributed with zero mean and standard deviation $\sigma = 0.005$.
The likelihood function is then given by
\begin{equation}
    \pi_{\rm like}( \obs \,| \, m ) \propto \exp\left(
    -\frac{1}{2}\|\mathcal{B}\,u(m) - \obs \|^2_{\ncov^{-1}}\right),
\end{equation}
where $\ncov = \sigma^2 \mathbf{I}$.

\subsubsection{Laplace approximation of the posterior}
\label{sec:laplace_approx_2d}

We next construct the Laplace approximation of the posterior, a Gaussian
distribution $\mulaplace \sim \GM{\iparmap}{\mathcal{H}(\iparmap)^{-1}}$ with mean equal to the MAP point and covariance given by
 the Hessian of the negative
log-posterior evaluated at the MAP point.
The MAP point is obtained by minimizing the negative log-posterior,
i.e.,
\begin{equation}
  \underset{m \in \mathcal{M}}{\min} \; \mathcal{J}(m) \;:=\; \frac{1}{2} \|
  \mathcal{B}\,u(m)- \obs \|^2_{ {\bf \Gamma}_{\text{noise}}^{-1}}
  +\frac{1}{2} \| m -\iparpr \|^2_{\Cprior^{-1}}.
\end{equation}

We employ the inexact Newton-CG algorithm implemented in \texttt{hIPPYlib} to
solve the above optimization problem.
We refer the reader to~\citet{VillaPetraGhattas21} for a detailed description of
the algorithm and the expressions for the gradient and Hessian actions of the
negative log-posterior~$\mathcal{J}(m)$.

As pointed out in Section \ref{sec:bayesian_framework}, explicitly
computing the Hessian is prohibitive for large-scale problems, as this entails solving two forward-like PDEs as many times as the number of parameters.
To make the operations with the Hessian scalable with respect to the parameter
dimension, we invoke a low-rank approximation of the data misfit part of the
Hessian, retaining only $r$ eigenvectors that are the most significantly
informed directions from the data \citep{VillaPetraGhattas21}.

\begin{figure}[t]
  \centering
  \resizebox{\textwidth}{!}{\input{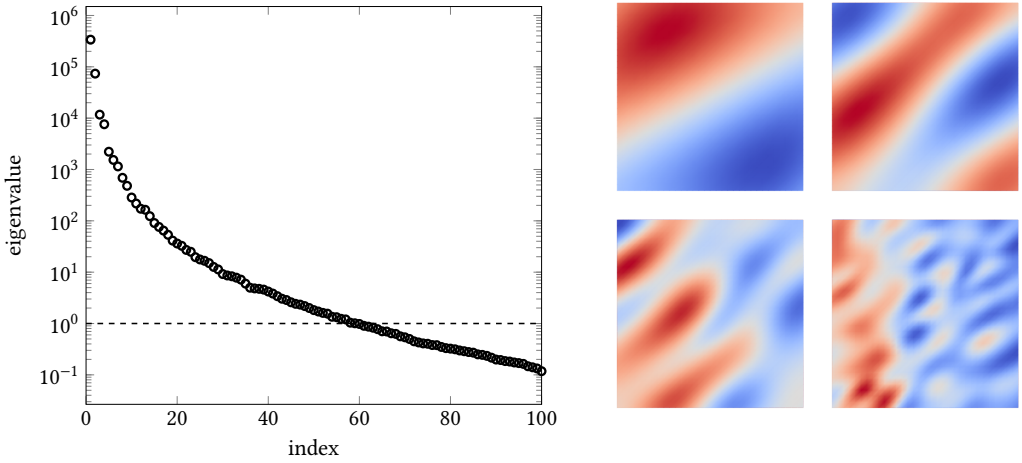}}
  \caption{Logarithmic plot of the $r=100$ dominant eigenvalues of the
  prior-preconditioned data misfit Hessian and the eigenvectors corresponding
  to the $1st, 4th, 16th$, and $64th$ largest eigenvalues for the Poisson
  problem.}%
  \label{fig:ex1_misfit_eig}
\end{figure}

\begin{figure}[t]
  \centering
  \resizebox{\textwidth}{!}{\input{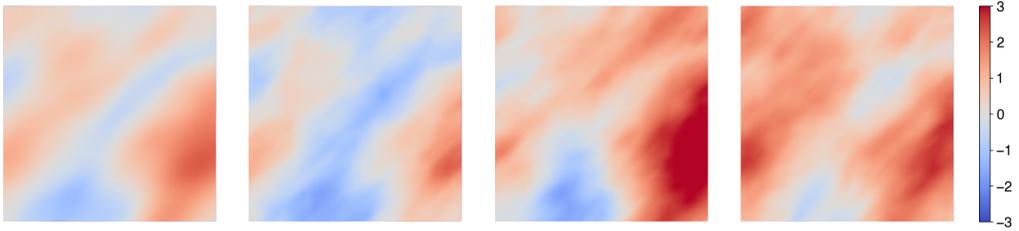}}
  \caption{The MAP point (leftmost) and three sample fields drawn from the
  Laplace approximation of the posterior distribution for the Poisson problem.}
  \label{fig:ex1_la_samples}
\end{figure}

Figure~\ref{fig:ex1_misfit_eig} shows the eigenspectrum of the
prior-preconditioned data misfit Hessian.
The double pass randomized algorithm provided by $\texttt{hIPPYlib}$
with an oversampling factor of 20 is used to accurately compute the dominant
eigenpairs.
We see that eigenvalues are smaller than 1 after around the $60th$ eigenvalue,
indicating that keeping 60 eigenpairs is sufficient for the low-rank
approximation.
Figure~\ref{fig:ex1_misfit_eig} also shows four eigenvectors, which, as expected, illustrate
that eigenvectors corresponding to smaller eigenvalues display more
fluctuations.

In Figure~\ref{fig:ex1_la_samples}, we depict the MAP point and three samples
drawn from the Laplace approximation of the posterior.

\subsubsection{Exploring the posterior using MCMC methods}
\label{sec:explore_post_2d}

In this section, we implement the advanced MCMC algorithms discussed in
the Section \ref{sec:geometry_aware_mcmc} to explore the posterior and compare their
performance.

In particular, we consider pCN, MALA, $\infty$-MALA, DR, DILI, and their
Hessian-informed counterparts.
For each method, we simulate 20 independent MCMC chains, each with 25,000
samples (after discarding 2,500 samples as burn-in), and hence draw a total of
500,000 samples from the posterior.  A sample from the Laplace approximation of
the posterior is chosen as a starting point for each chain.

For checking the convergence and statistical efficiency of MCMC chains, we
consider the subspace spanned by the $r$ dominant eigenvectors of the
generalized eigensystem in \eqref{eq:eigenproblem_1}, instead of all components
of the parameter vector $\vec{\ipar}$.
Specifically, we compute the MPSRF, autocorrelation time, and ESS with respect to a
coefficient vector $\mathbf{c} \in \mathbb{R}^r$ defined by
\begin{equation}
  \mathbf{c} = \mathbf{V}_r^T \prcov^{-1} \mathbf{m}.
\end{equation}

\begin{table}[t]
    \centering
    \caption{Comparison of the performance of several MCMC methods for the
      Poisson problem: pCN, MALA, $\infty$-MALA, DR, DILI, and their
      Hessian-informed versions.
      Acceptance rate (AR), multivariate potential scale reduction
      factor (MPSRF), and effective sample sample size (ESS)
      are reported for convergence diagnostics.
      MPSRF and ESS are computed with respect to the projection of parameter
      samples along the first 25 dominant eigenvectors of the
      prior-preconditioned data misfit Hessian at the MAP point.
      Two values of AR are listed in DR and DILI-MAP, which are for the first
      and the second proposal moves, respectively. We also provide the number
      of forward and/or adjoint PDE solves per effective sample (NPS/ES) for
      sampling efficiency. We use 20 MCMC chains, each with 25,000 iterations
      (500,000 samples in total).
      The numbers in parentheses in each method name represent the parameter
      values used ($\beta$ for pCN, $\tau$ for MALA, $h$ for
      $\infty$-MALA, and $\beta$ and $\tau$ for and DILI). The numbers in parentheses
      of the minimum ESS and the maximum ESS indicate the corresponding eigenvector
      index.}
    \label{tab:convergence_diagnostics}
    \begin{adjustbox}{max width=\textwidth}
      \begin{tabular}{l|cccccc}
      {\bf Method} & {\bf AR (\%)} & {\bf MPSRF} & {\bf Min. ESS (index)} & {\bf Max. ESS (index)} & {\bf Avg. ESS} & {\bf NPS/ES} \\
        \hline
      pCN (5.0E-3) &24 &2.629 &25 (24) &225 (8) &84 &5,952 \\
      MALA (6.0E-6) &48 &2.642 &26 (22) &874 (5) &148 &10,135 \\
      $\infty$-MALA (1.0E-5) &57 &2.943 &25 (23) &1,102 (5) &160 &9,375 \\
      \hdashline
      H-pCN (4.0E-1) &27 &1.192 &64 (1) &3,598 (15) &2,314 &216 \\
      H-MALA (6.0E-2) &60 &1.014 &545 (1) &8,868 (19) &6,459 &232 \\
      H-$\infty$-MALA (1.0E-1) &71 &1.016 &582 (1) &8,417 (18) &5,905 &254 \\
      \hdashline
      DR (H-pCN (1.0E0), H-MALA (6.0E-2)) &(4, 61) &1.013 &641 (1) &12,522 (17) &9,222 &215 \\
      DR (H-pCN (1.0E0), H-$\infty$-MALA (2.0E-1)) &(4, 48) &1.011 &613 (1) &12,812 (17) &9,141 &213 \\
      \hdashline
      DILI-PRIOR (0.8, 0.1) &(60, 33) &1.064 &314 (1) &4,667 (13) &3,216 &548 \\
      DILI-LA (0.8, 0.1) &(83, 36) &1.017 &562 (1) &10,882 (17) &7,192 &245 \\
      DILI-MAP (0.8, 0.1) &(77, 22) &1.006 &1,675 (1) &10,271 (20) &8,692 &202 %
    \end{tabular}
    \end{adjustbox}
\end{table}

Table~\ref{tab:convergence_diagnostics} shows the convergence diagnostics and
computational efficiency of the MCMC samples.
MPSRF and ESS are computed with respect to the projection of parameter samples
along the first 25 dominant eigenvectors of the prior-preconditioned data
misfit Hessian at the MAP point.
Table~\ref{tab:convergence_diagnostics} reports the
mininum, maximum, and average ESS over all the 25 projections.

The last column in Table~\ref{tab:convergence_diagnostics} represents the
number of forward and/or adjoint PDE solves required to draw a single
independent sample (average ESS is used).
This quantity can be used to measure the sampling efficiency and rank the
methods in terms of computational efficiency.
Under this metric, DILI-MAP is
the most efficient method and requires only 202 PDE solves for an effective sample.
DR (213 NPS/ES for H-$\infty$-MALA and 215 NPS/ES for H-MALA) and H-pCN (216 NPS/ES) are close seconds.

\begin{figure}[t]
    \centering
    \resizebox{0.8\textwidth}{!}{\begin{tikzpicture}
  \begin{axis}[
    scale only axis,
    xlabel={lag},
    ylabel={autocorrelation},
    xmin=0, xmax=500,
    ymin=-0.2, ymax=1.1,
    legend style={nodes=right, anchor=north west, draw=none},
    legend pos=outer north east,
    every axis plot/.append style={very thick}
    ]
    \definecolor{olive}{rgb}{0.5,0.5,0}
    \definecolor{darkgoldenrod}{rgb}{0.72,0.53,0.4}
    \definecolor{orange}{rgb}{1,0.65,0}
    \definecolor{purple}{rgb}{0.5,0,0.5}
    \definecolor{violet}{rgb}{0.93,0.51,0.93}

    \addplot[black, dashed,
    domain=\pgfkeysvalueof{/pgfplots/xmin}:\pgfkeysvalueof{/pgfplots/xmax},
    samples=2, forget plot] {0.0};
    \addplot[olive] table[x=lag,y=pcn] {data/ex1_acf.dat};
    \addlegendentry{pCN};
    \addplot[darkgoldenrod] table[x=lag,y=mala] {data/ex1_acf.dat};
    \addlegendentry{MALA};
    \addplot[orange] table[x=lag,y=infmala] {data/ex1_acf.dat};
    \addlegendentry{$\infty$-MALA};

    \addplot[blue] table[x=lag,y=hpcn] {data/ex1_acf.dat};
    \addlegendentry{H-pCN};
    \addplot[green!60!black] table[x=lag,y=hmala] {data/ex1_acf.dat};
    \addlegendentry{H-MALA};
    \addplot[cyan] table[x=lag,y=infhmala] {data/ex1_acf.dat};
    \addlegendentry{H-$\infty$-MALA};

    \addplot[darkgray] table[x=lag,y=dr] {data/ex1_acf.dat};
    \addlegendentry{DR (H-pCN, H-MALA)};
    \addplot[gray] table[x=lag,y=dr_inf] {data/ex1_acf.dat};
    \addlegendentry{DR (H-pCN, $\infty$-H-MALA)};

    \addplot[purple] table[x=lag,y=dili_prior] {data/ex1_acf.dat};
    \addlegendentry{DILI-PRIOR};
    \addplot[violet] table[x=lag,y=dili_la] {data/ex1_acf.dat};
    \addlegendentry{DILI-LA};
    \addplot[red] table[x=lag,y=dili_map] {data/ex1_acf.dat};
    \addlegendentry{DILI-MAP};
  \end{axis}
\end{tikzpicture}}
    \caption{Autocorrelation function estimate~\eqref{eqn:acf}
    of the quantity of interest $\qoi$~\eqref{eqn:qoi}
    for several MCMC methods.
    Note that the autocorrelation function plots for pCN, MALA, and
    $\infty$-MALA appear unchanged and overlap.
    }%
    \label{fig:ex1_mcmc_acf}
\end{figure}
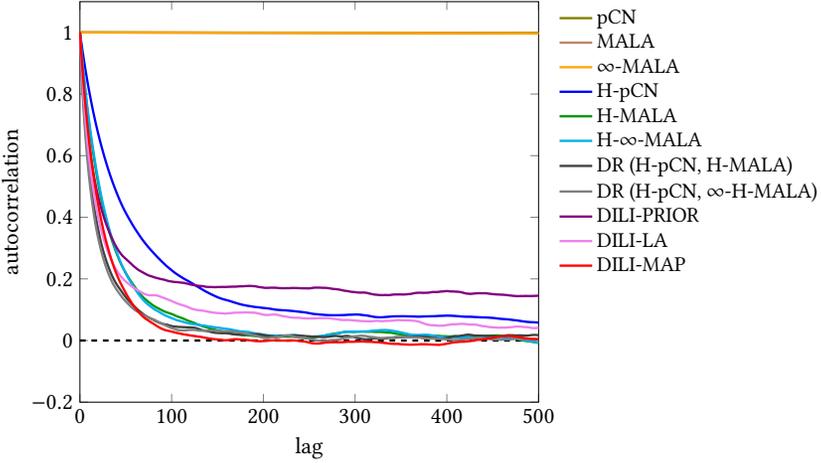

\begin{figure}[t]
  \centering
  \resizebox{\textwidth}{!}{\input{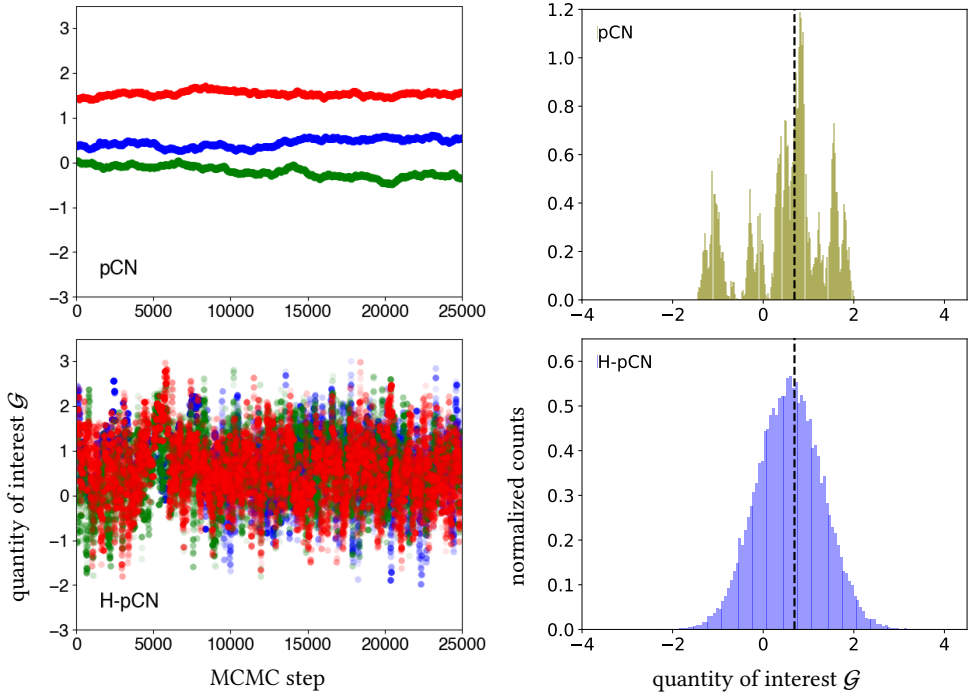}}
  \caption{Left: Trace plots of the quantity of interest $\qoi$~\eqref{eqn:qoi}
    from three MCMC chains (out of 20 independent chains) of pCN and H-pCN simulations;
    different colors (here blue, green and red) represent the trace of each
    chain.
    Right: Probability density function estimate of the quantity of interest
    $\qoi$~\eqref{eqn:qoi} computed from pCN and H-pCN samples; all the 500,000
    samples, 20 chains with 25,000 samples each, are pulled together in the
    histogram; the number of counts is normalized so that the plot represents a
    probability density function; the black dashed line represents the quantity
    of interest computed from the true parameter field $m_\text{true}$.
    }
  \label{fig:ex1_mcmc_trace_hist}
\end{figure}

We next assess the convergence of MCMC samples of the quantity of interest
in \eqref{eqn:qoi} to the predictive posterior distribution of
$\qoi(m)$: the autocorrelation function estimates of the quantity of interest
~\eqref{eqn:qoi} are shown in Figure~\ref{fig:ex1_mcmc_acf} (here, we use
the formula~\eqref{eqn:acf} to account for the use of multiple chains), 
and trace plots from three independent MCMC chains and histograms of all the MCMC samples for pCN and H-pCN
are depicted in Figure~\ref{fig:ex1_mcmc_trace_hist}.

\begin{figure}[t]
  \centering
  \resizebox{\textwidth}{!}{\input{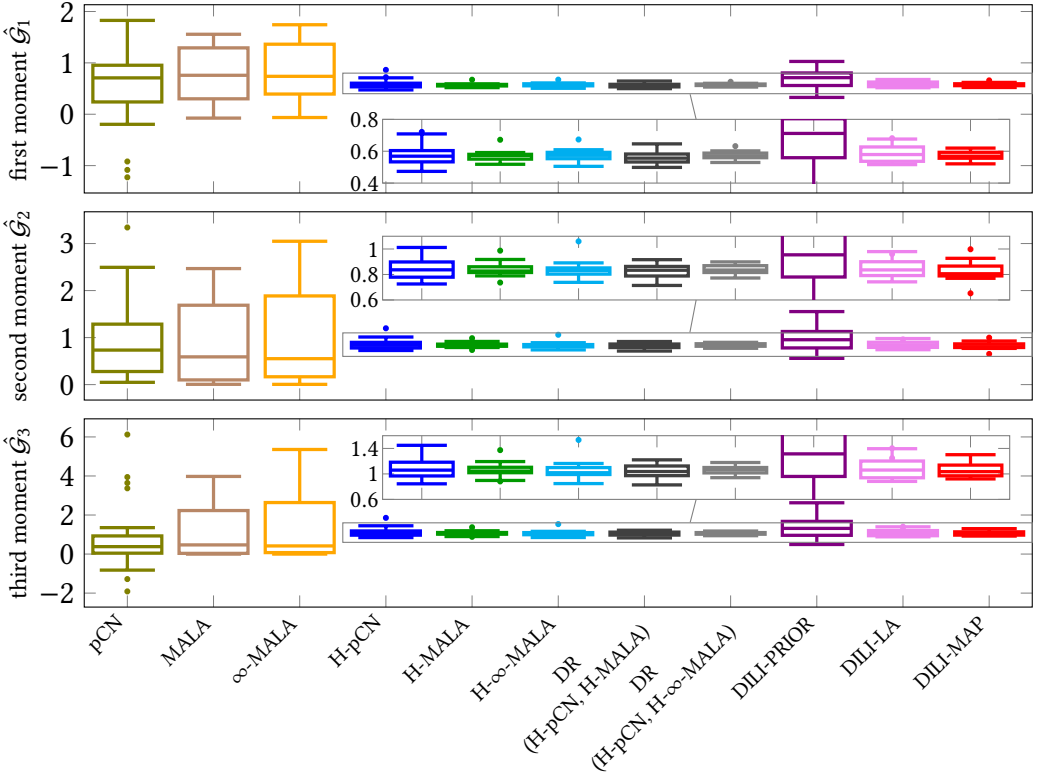}}
  \caption{Whisker plots of first, second, and third moment estimates
    ($\hat{\qoi}_k, k=1,2,3$) of the quantity of
    interest~\eqref{eqn:moment_qoi} computed by using several MCMC
    methods. The central mark is the median; lower and upper quartiles
    represent 25th and 75th percentiles, respectively. Whiskers extend
    to the extreme data points that fall within the distance from the
    lower or upper quartiles to 1.5 times the interquartile range (the
    distance between the upper and lower quartiles); all the other
    data points are plotted as outliers. The number of data points for
    each method is 20, the number of independent MCMC chains.}
  \label{fig:ex1_moment_matching}
\end{figure}

Lastly, we compare estimates of moments of the quantity of interest for the
different sampling strategies.  For each MCMC chain, the $k$th ($k = 1,2,3$)
moment of the quantity of interest computed from parameter samples
$\mathbf{m}_i$ ($i=1, 2, \ldots, N; N = 25,000$) is computed as
\begin{equation}
    \label{eqn:moment_qoi}
    \hat{\qoi}_k = \frac{1}{N} \sum_{i=1}^{N} \qoi^k (\mathbf{m}_i).
\end{equation}
The results are reported in Figure~\ref{fig:ex1_moment_matching} as
box-and-whisker plots.

From the results presented in this section, we draw the following
conclusions:
\begin{itemize}
  \item The Hessian information at the MAP point plays an important role in
    enhancing the sampling performance of the MCMC methods. In fact MCMC chains
    without the Hessian information did not converge over the entire length of
    the chain and were localized around the starting point. The convergence was
    achieved or nearly achieved only when the MCMC proposal exploited the Laplace approximation of
    the posterior that incorporates the Hessian information.
  \item DILI-MAP shows the best sampling efficiency in terms of the number of
    forward and/or adjoint PDE solves per effective sample. Note that the
    parameter value used in the MCMC methods (e.g., $\beta$ and/or $\tau$) was
    not the optimal and a different result may be obtained with different
    parameter values.
\end{itemize}
We further study the performance of MCMC methods under different problem
settings to provide more insight into the practical use of the hIPPYlib-MUQ framework.

\subsubsection{Scalability of Hessian-informed pCN}
\label{sec:mcmc_scalability_2d}

Here we investigate the effect of mesh resolution on the sampling performance.
A curvature aware MCMC method, the H-pCN is selected with $\beta=0.4$ for the
test.  The dimensions of the parameter and the state variables from a coarse
mesh (mesh 1) to the finest mesh (mesh~4) are (1,089, 4,225), (4,225, 16,641),
(16,641, 66,049), and (66,049, 263,169), respectively.

\begin{figure}[t]
  \centering
  \resizebox{\textwidth}{!}{\input{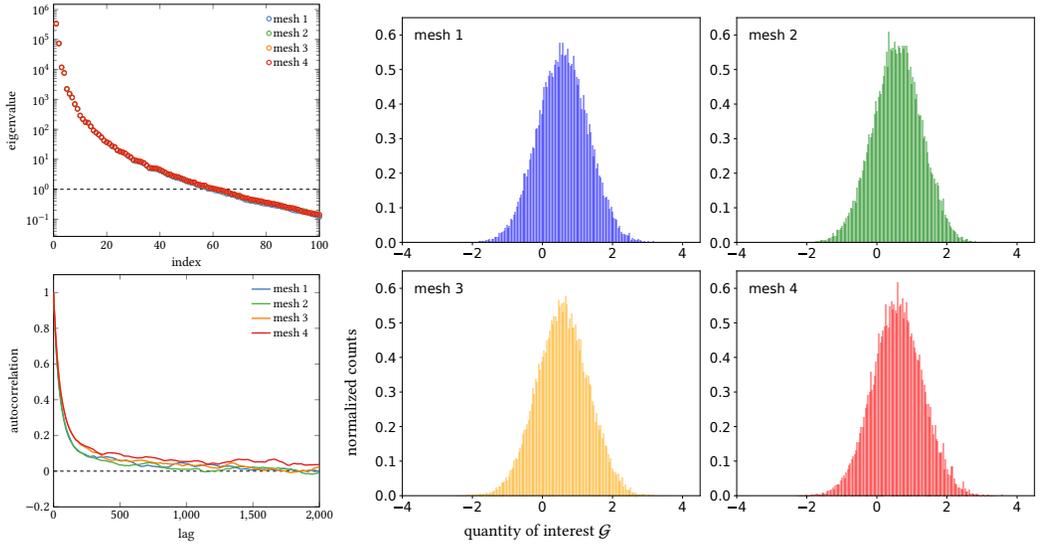}}
  \caption{Left:
  Logarithmic plot of the $r=100$ dominant eigenvalues of the
  prior-preconditioned data misfit Hessian (top), and
  autocorrelation function estimate~\eqref{eqn:acf} of the quantity of interest
  $\qoi$~\eqref{eqn:qoi} (bottom).
  Right: Probability density function
  estimate of the quantity of interest $\qoi$~\eqref{eqn:qoi}; all the samples,
  20 chains with 25,000 samples each, so 500,000 in total, are pulled together
  in the histogram; the number of counts is normalized so that the plot
  represents a probability density function.
  We consider four different meshes which are increasingly refined from the
  coarsest (mesh 1) to the finest (mesh~4).
  We use the H-pCN method ($\beta=0.4$) to draw samples.
  }
  \label{fig:ex1_eig_acf_hist_4meshes}
\end{figure}

We follow the same problem setting as before, and use the same
synthetic observations (obtained from the true parameter field generated from the
finest mesh) for all levels.
The top-left part of Figure~\ref{fig:ex1_eig_acf_hist_4meshes} shows the $r=100$ dominant eigenvalues
of the prior-preconditioned data misfit Hessian.  One observes that the
eigenspectrum is virtually independent of mesh refinement.

\begin{table}[t]
    \centering
    \caption{Acceptance rate (AR), multivariate potential scale reduction
      factor (MPSRF) and effective sample size (ESS) of the posterior samples
      generated by using the H-pCN method with different dimensions.
      We use $\beta = 0.4$ for the H-pCN method and draw in total 500,000
      samples (20 MCMC chains, each with 25,000 iterations).
      MPSRF and ESS are computed with respect to the projection of samples along the first
      25 dominant eigenvectors of the prior-preconditioned data misfit Hessian at
      the MAP point.
      The numbers in parentheses of the minimum ESS and the maximum ESS indicate
    the corresponding eigenvector index.}
    \label{tab:ex1_diagnostic_4meshes}
    \begin{adjustbox}{max width=\textwidth}
    \begin{tabular}{lccccc}
      {\bf Dimension (state, parameter)} & {\bf AR (\%)} & {\bf MPSRF} & {\bf Min. ESS (index)} & {\bf Max.
      ESS (index)} & {\bf Avg. ESS}\\
      \hline
      (4,225, 1,089) &27 &1.192 &64 (1) &3,598 (15) &2,314 \\
      (16,641, 4,225) &24 &1.333 &63 (1) &3,221 (18) &1,830\\
      (66,049, 16,641) &23 &1.075 &209 (1) &3,073 (11) &1,940\\
      (263,169, 66,049) &22 &1.117 &102 (2) &3,276 (15) &1,767 \\
    \end{tabular}
    \end{adjustbox}
\end{table}

To assess the convergence of the MCMC methods, in Table~\ref{tab:ex1_diagnostic_4meshes}
we report the acceptance rate, MPSRF, and ESS of the posterior samples. The MPSRF and ESS are computed
with respect to the projection of parameter samples along the first 25 dominant eigenvectors of the
prior-preconditioned data misfit Hessian at the MAP point. 
We present the autocorrelation function estimates~\eqref{eqn:acf} in the
bottom-left part of Figure~\ref{fig:ex1_eig_acf_hist_4meshes}, and
show histograms for the quantity of interest $\qoi$~\eqref{eqn:qoi} in the right
part of Figure~\ref{fig:ex1_eig_acf_hist_4meshes}.
The results show that %
the convergence of samples is almost independent with respect to the MPSRF and the
autocorrelation function.

\subsection{Coefficient field inversion in a Robin boundary condition for a three-dimensional $p$-Poisson nonlinear PDE}

So far, we restricted our attention to the additive Gaussian noise
model. While this additive noise model is the most commonly used
model, it would be inappropriate in some cases such as speckle noise
found in synthetic aperture radar (SAR) images. In this example, we
consider a different noise model where the noise is proportional to
the value of the observations.

The example forward model is a nonlinear PDE in three space dimensions for
which we seek to infer an unknown coefficient field in a Robin boundary
condition. %
Specifically, the forward governing equations are given by
\begin{alignat}{2}
    \label{eqn:ex2_forward_strong}
    -\nabla \cdot \left( |\nabla u|^{p-2}_\epsilon \nabla u\right) &= f &\quad
    &\text{in } \Omega, \nonumber \\
    |\nabla u|^{p-2}_\epsilon \nabla u \cdot \mathbf{n} + e^m u
    &= 0 &&\text{on } \partial \Omega_R, \\
    |\nabla u|^{p-2}_\epsilon \nabla u \cdot \mathbf{n}
    &= 0
    &&\text{on }
    \partial \Omega \setminus \partial \Omega_R,
    \nonumber
\end{alignat}
with $1 \le p \le \infty$.
Note that the $p$-Laplacian, $\nabla \cdot \left( |\nabla
u|^{p-2} \nabla u\right)$, is singular when $p<2$ and degenerates when $p>2$ at
points $\nabla u =0$~\citep{lindqvist2017notes, Brown2010}, so a regularization
term $\epsilon$ (here we take $\epsilon = 1.0 \times 10^{-8}$) is introduced in the above equation as $|\nabla u|_\epsilon =
\sqrt{|\nabla u|^2 + \epsilon}$.
The $p$-Laplacian is a nonlinear counterpart of the Laplacian operator, and
appears in many nonlinear diffusion problems (e.g., non-Newtonian fluids),
where a nonlinear diffusion is modeled as a power law type.

We assume $p=3$ and consider a thin brick domain $\Omega = [0, 1]^2
\times [0, 0.05]$ with a volume source term ($f = 1$) and a mixed
boundary condition, e.g., we impose a Robin boundary condition on the
bottom boundary surface and no normal flux on the remaining boundary
surfaces.

The problem domain $\Omega$ is discretized using a regular tetrahedral grid and using linear finite
elements for all the state, adjoint, and parameter variables. 
After discretization, the dimension is 66,564 for the state and adjoint
variables, and 16,641 for the parameter variable.

\begin{figure}[t]
    \centering
    \resizebox{\textwidth}{!}{\input{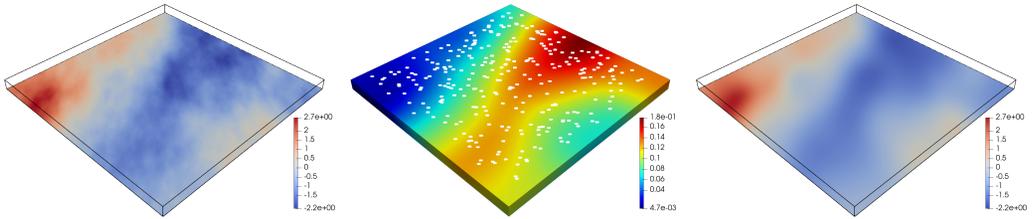}}
    \caption{Left: True parameter field on the bottom surface.
    Middle: Corresponding state field and
    $l=300$ observation points (white square marks) on the top surface. 
    Right:
    The MAP point on the bottom surface.}
    \label{fig:ex2_mtrue_observation}
\end{figure}

For the prior, we use the same Gaussian distribution as the one employed in the
previous example. 
We create synthetic state observations at $l=300$ random
locations uniformly distributed on the top boundary surface by first solving the forward
problem with the true parameter field $m_{\textrm{true}}$ obtained from a sample of 
the prior, and then multiplying a Gamma-distributed noise.
Specifically, the vector of synthetic observations is of the form
\begin{equation}
  \obs = \vec{\eta} \odot \mathcal{B} u,
  \label{eq:ex2_multiplicative_noise}
\end{equation}
where $\odot$ denotes component-wise multiplication, and each
component of $\vec{\eta}$ is independently and identically
Gamma-distributed with shape $\kappa$ and scale $\nu$, i.e., $\eta_i =
\text{Gamma}(\kappa, \nu), i=1, \ldots, q$.
We take $\kappa = 1 / \nu$ and in this case the negative log-likelihood function has the form
\begin{equation}
  \Phi(\ipar; \obs) = \kappa \sum_{i=1}^{q} \left( \log (\mathcal{B} u)_i +
  \frac{d_i}{(\mathcal{B} u)_i} \right),
  \label{eq:ex2_negative_loglikelihood}
\end{equation}
where the subscript $i$ means $i$th component of the corresponding vector.
In this example, we set $\kappa = 10^{4}$.

Figure~\ref{fig:ex2_mtrue_observation}
illustrates the true parameter field on the bottom boundary, the locations of
the observations on the top surface, and the MAP point obtained by solving the
optimization problem of minimizing the negative log-posterior.
The Laplace approximation of the posterior is then constructed based on the
low-rank factorization of the data misfit Hessian at the MAP point.
The spectrum of the prior-preconditioned data misfit Hessian indicates that
the number of dominant eigenvalues (larger than 1) is about 55.

\subsubsection{MCMC results for characterizing the posterior}

Here we present MCMC simulation results for the uncertain boundary
coefficient vector.  In this example, we consider the H-pCN method
with $\beta=0.2$ and run 40 independent MCMC chains, each with 25,000
iterations after discarding 2,500 samples as burn-in (1,000,000
samples are generated in total).  For each MCMC run, a sample from the
Laplace approximation of the posterior is taken as the starting point.

\begin{table}[t]
  \centering
  \caption{Convergence diagnostics for the $p$-Poisson problem:
    acceptance rate (AR), multivariate potential scale reduction
    factor (MPSRF), and effective sample sample size (ESS) of the projection of
    the parameter samples along the first 25 eigenvectors of the
    prior-preconditioned data misfit Hessian at the MAP point. We use H-pCN
    method ($\beta = 0.2$) with 40 chains, each with 25,000 iterations (1,000,000
  samples in total). The numbers in parentheses of the minimum ESS and the
  maximum ESS indicate the corresponding eigenvector index.}
  \label{tab:ex2_diagnostic}
  \begin{tabular}{ccccc}
    {\bf AR (\%)} & {\bf MPSRF} & {\bf Min. ESS (index)} & {\bf Max. ESS (index)} & {\bf Avg. ESS} \\
    \hline
    23 & 1.041 & 243 (0) & 2,891 (22) & 1,586
  \end{tabular}
\end{table}

\begin{figure}[t]
  \centering
  \resizebox{0.8\textwidth}{!}{\input{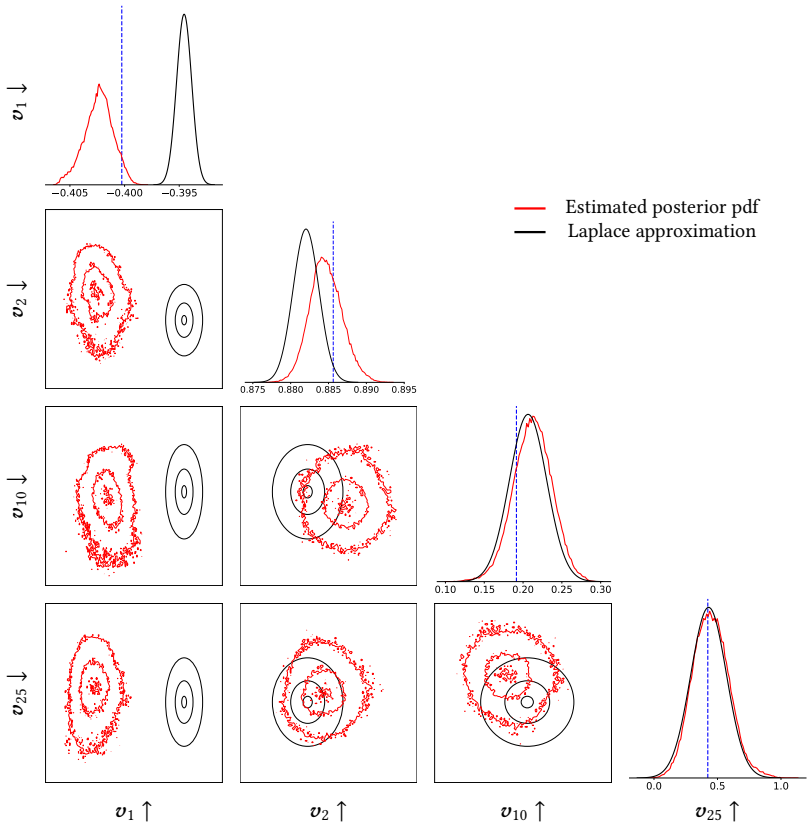}}
  \caption{Marginal distributions of the MCMC samples (red) and the
    Laplace approximation (black), both for the projection of
    parameter samples onto the eigenvectors $\vec{v}_1, \vec{v}_2,
    \vec{v}_{10}, \text{and}, \vec{v}_{25}$.  In the one-dimensional
    marginal plots, the blue dashed lines represent the projection of
    the true parameter vector onto the eigenvectors.  For the
    two-dimensional marginal plots, we only show three contours that
    represent 5\%, 50\%, and 95\% of the distribution,
    respectively. We note the difference in scaling and hence the
    appearance of a much larger difference between the estimated and
    the Laplace approximation of the posterior in the first column.}
  \label{fig:ex2_marginals}
\end{figure}

Table~\ref{tab:ex2_diagnostic} lists convergence diagnostics of the MCMC
simulation.  The parameter samples are projected onto the first 25
eigenvectors of the prior-preconditioned data misfit Hessian at the MAP point,
and the MPSRF and the ESS are evaluated based on this projection.

Figure~\ref{fig:ex2_marginals} shows marginal distributions of the
posterior MCMC samples and of the Laplace approximation. As before,
these marginals are computed with respect to the projection of the
parameter samples onto the eigenvectors.  We observe that there is a
clear difference between the marginal distributions of the MCMC
samples and those of the Laplace approximation, especially for the
eigenvector corresponding to larger eigenvalues.

The numerical studies have been carried out on the Multi-Environment
Computer for Exploration and Discovery (MERCED) cluster at UC Merced.
While there are some performance variations depending on the cluster
CPU node that a job is assigned to, approximately, a single PDE solve
required about 0.07 seconds for the first example, and about 4 seconds
for the second three-dimensional example (both are the elapsed time).

\section{Conclusion}
\label{sec:conclusions}

We have presented a robust and scalable software framework for the
solution of large-scale Bayesian inverse problems governed by
PDEs. The software integrates two complementary open-source software
libraries, hIPPYlib and MUQ, resulting in a unique software framework
that addresses the prohibitive nature of Bayesian solution of inverse
problems governed by PDEs. The main objectives of the proposed
software framework are to
\begin{enumerate}
\item provide to domain scientists a suite of sophisticated and
  computationally efficient MCMC methods that exploit Bayesian inverse
  problem structure; and
\item allow researchers to easily implement new methods and compare
  against the state of the art.
\end{enumerate}

The integration of the two libraries allows advanced MCMC methods to
exploit the geometry and intrinsic low-dimensionality of parameter
space, leading to efficient and scalable exploration of the posterior
distribution. In particular, the Laplace approximation of the
posterior is employed to generate high-quality MCMC proposals. This
approximation is based on the inverse of the Hessian of the
log-posterior, made tractable via low-rank approximation of the
Hessian of the log-likelihood. Numerical experiments on linear and
nonlinear PDE-based Bayesian inverse problems illustrate the ability
of Laplace-based proposals to accelerate MCMC sampling by factors of
$\sim 50\times$.

Despite the fast and dimension-independent convergence of these
advanced structure-exploiting MCMC methods, many Bayesian inverse
problems governed by expensive-to-solve PDEs remain out of reach. For
example, the results of section~\ref{sec:explore_post_2d} for the
Poisson coefficient inverse problem indicate that $O(10^6)$ PDE solves
may still be required even with the most efficient MCMC methods. In
such cases, hIPPYlib-MUQ can be used as a prototyping environment to
study new methods that further exploit problem structure, for example
through the use of various reduced models
(e.g.,~\citet{CuiMarzoukWillcox16}) or via advanced Hessian
approximations that go beyond low rank~\citet{AlgerRaoMeyersEtAl19,
  AmbartsumyanBoukaramBui-ThanhEtAl20}.

We also remark on limitations of our software framework. Since we rely
on FEniCS for the finite element approximation of PDEs, hIPPYlib-MUQ
inherits all the limitations and challenges that come with this
software. However, the FEniCS, hIPPYlib and MUQ developers and user
communities offer a rich support base that the users of hIPPYlib-MUQ
can build on. Support for alternative finite element implementations is
the subject of future work. Another current limitation of hIPPYlib-MUQ
is the lack of parallel implementations of MCMC methods. Therefore,
the goal for future versions of hIPPYlib-MUQ is to incorporate
multilevel parallelism.  This will include both the parallel PDE
solvers available now and additional parallel MCMC chains, will allow
solutions of even more complex PDE-based Bayesian inverse problems
with higher-dimensional parameter spaces.

\paragraph{Software Availability} hIPPYlib-MUQ is distributed under the GNU General Public License version 3 (GPL3). 
The hIPPYlib-MUQ project is hosted on GitHub (\url{https://github.com/hippylib/hippylib2muq}) and use Travis-CI for continuous integration. hIPPYlib-MUQ uses semantic versioning. The results presented in this work were obtained with hIPPYlib-MUQ version 0.3.0, hIPPYlib version 3.1.0, and MUQ version 0.3.5. A Docker image~\citep{Merkel14} containing the pre-installed software and examples is available at \url{https://hub.docker.com/r/ktkimyu/hippylib2muq}. 
hIPPYlib-MUQ documentation is hosted on ReadTheDocs (\url{https://hippylib2muq.readthedocs.io}). Users are encouraged to join the hIPPYlib and MUQ workspaces on Slack to connect with other users, get help, and discuss new features; see \url{https://hippylib.github.io/#slack-channel} and \url{https://mituq.bitbucket.io} for more information on how to join.

\begin{acks}
This work was supported by the U.S. National Science Foundation,
Software Infrastructure for Sustained Innovation (SI2: SSE \& SSI)
Program under grants ACI-1550593, ACI-1550547, and ACI-1550487 and the
Division of Mathematical Sciences under the CAREER grant 1654311. MP
and YM were also supported in part by Office of Naval Research MURI
grant N00014-20-1-2595. OG was also supported in part by Department of
Energy Advanced Scientific Computing Research grants DE-SC0021239 and
DE-SC0019303. The authors gratefully acknowledge computing time on the
Multi-Environment Computer for Exploration and Discovery (MERCED)
cluster at UC Merced, which was funded by National Science Foundation
Grant No. ACI-1429783. We thank the anonymous reviewers for their
careful reading of our manuscript and for their valuable comments and
suggestions, which helped us improve the quality of the manuscript.
\end{acks}

\bibliographystyle{ACM-Reference-Format}

\begin{thebibliography}{74}


\ifx \showCODEN    \undefined \def \showCODEN     #1{\unskip}     \fi
\ifx \showDOI      \undefined \def \showDOI       #1{#1}\fi
\ifx \showISBNx    \undefined \def \showISBNx     #1{\unskip}     \fi
\ifx \showISBNxiii \undefined \def \showISBNxiii  #1{\unskip}     \fi
\ifx \showISSN     \undefined \def \showISSN      #1{\unskip}     \fi
\ifx \showLCCN     \undefined \def \showLCCN      #1{\unskip}     \fi
\ifx \shownote     \undefined \def \shownote      #1{#1}          \fi
\ifx \showarticletitle \undefined \def \showarticletitle #1{#1}   \fi
\ifx \showURL      \undefined \def \showURL       {\relax}        \fi
\providecommand\bibfield[2]{#2}
\providecommand\bibinfo[2]{#2}
\providecommand\natexlab[1]{#1}
\providecommand\showeprint[2][]{arXiv:#2}

\bibitem[\protect\citeauthoryear{Ak\c{c}elik, Biros, Ghattas, Hill, Keyes, and
  van Bloeman~Waanders}{Ak\c{c}elik et~al\mbox{.}}{2006}]%
        {AkcelikBirosGhattasEtAl06a}
\bibfield{author}{\bibinfo{person}{Volkan Ak\c{c}elik}, \bibinfo{person}{George
  Biros}, \bibinfo{person}{Omar Ghattas}, \bibinfo{person}{Judith Hill},
  \bibinfo{person}{David Keyes}, {and} \bibinfo{person}{Bart van
  Bloeman~Waanders}.} \bibinfo{year}{2006}\natexlab{}.
\newblock \showarticletitle{Parallel {PDE}-constrained optimization}.
\newblock In \bibinfo{booktitle}{\emph{Parallel Processing for Scientific
  Computing}}, \bibfield{editor}{\bibinfo{person}{M.~Heroux},
  \bibinfo{person}{P.~Raghaven}, {and} \bibinfo{person}{H.~Simon}} (Eds.).
  \bibinfo{publisher}{SIAM}.
\newblock


\bibitem[\protect\citeauthoryear{Alger, Rao, Meyers, Bui-Thanh, and
  Ghattas}{Alger et~al\mbox{.}}{2019}]%
        {AlgerRaoMeyersEtAl19}
\bibfield{author}{\bibinfo{person}{N. Alger}, \bibinfo{person}{V. Rao},
  \bibinfo{person}{A. Meyers}, \bibinfo{person}{T. Bui-Thanh}, {and}
  \bibinfo{person}{O. Ghattas}.} \bibinfo{year}{2019}\natexlab{}.
\newblock \showarticletitle{Scalable matrix-free adaptive product-convolution
  approximation for locally translation-invariant operators}.
\newblock \bibinfo{journal}{\emph{SIAM Journal on Scientific Computing}}
  \bibinfo{volume}{41}, \bibinfo{number}{4} (\bibinfo{year}{2019}),
  \bibinfo{pages}{A2296--A2328}.
\newblock
\urldef\tempurl%
\url{https://arxiv.org/abs/1805.06018}
\showURL{%
\tempurl}


\bibitem[\protect\citeauthoryear{Ambartsumyan, Boukaram, Bui-Thanh, Ghattas,
  Keyes, Stadler, Turkiyyah, and Zampini}{Ambartsumyan et~al\mbox{.}}{2020}]%
        {AmbartsumyanBoukaramBui-ThanhEtAl20}
\bibfield{author}{\bibinfo{person}{Ilona Ambartsumyan}, \bibinfo{person}{Wajih
  Boukaram}, \bibinfo{person}{Tan Bui-Thanh}, \bibinfo{person}{Omar Ghattas},
  \bibinfo{person}{David Keyes}, \bibinfo{person}{Georg Stadler},
  \bibinfo{person}{George Turkiyyah}, {and} \bibinfo{person}{Stefano Zampini}.}
  \bibinfo{year}{2020}\natexlab{}.
\newblock \showarticletitle{Hierarchical Matrix Approximations of {H}essians
  Arising in Inverse Problems Governed by {PDE}s}.
\newblock \bibinfo{journal}{\emph{SIAM Journal on Scientific Computing}}
  \bibinfo{volume}{42}, \bibinfo{number}{5} (\bibinfo{year}{2020}),
  \bibinfo{pages}{A3397--A3426}.
\newblock


\bibitem[\protect\citeauthoryear{Atchad\'e}{Atchad\'e}{2006}]%
        {Atchade06}
\bibfield{author}{\bibinfo{person}{Yves~F. Atchad\'e}.}
  \bibinfo{year}{2006}\natexlab{}.
\newblock \showarticletitle{An adaptive version for the {M}etropolis adjusted
  {L}angevin algorithm with a truncated drift}.
\newblock \bibinfo{journal}{\emph{Methodology and Computing in Applied
  Probability}}  \bibinfo{volume}{8} (\bibinfo{year}{2006}),
  \bibinfo{pages}{235--254}.
\newblock


\bibitem[\protect\citeauthoryear{Balay, Abhyankar, Adams, Brown, Brune,
  Buschelman, Dalcin, Dener, Eijkhout, Gropp, Kaushik, Knepley, May, McInnes,
  Mills, Munson, Rupp, Sanan, Smith, Zampini, and Zhang}{Balay
  et~al\mbox{.}}{2018}]%
        {petsc-web-page}
\bibfield{author}{\bibinfo{person}{Satish Balay}, \bibinfo{person}{Shrirang
  Abhyankar}, \bibinfo{person}{Mark~F. Adams}, \bibinfo{person}{Jed Brown},
  \bibinfo{person}{Peter Brune}, \bibinfo{person}{Kris Buschelman},
  \bibinfo{person}{Lisandro Dalcin}, \bibinfo{person}{Alp Dener},
  \bibinfo{person}{Victorand Eijkhout}, \bibinfo{person}{William~D. Gropp},
  \bibinfo{person}{Dinesh Kaushik}, \bibinfo{person}{Matthew~G. Knepley},
  \bibinfo{person}{Dave~A. May}, \bibinfo{person}{Lois~Curfman McInnes},
  \bibinfo{person}{Richard~Tran Mills}, \bibinfo{person}{Todd Munson},
  \bibinfo{person}{Karl Rupp}, \bibinfo{person}{Patrick Sanan},
  \bibinfo{person}{Barry~F. Smith}, \bibinfo{person}{Stefano Zampini}, {and}
  \bibinfo{person}{Hong Zhang}.} \bibinfo{year}{2018}\natexlab{}.
\newblock \bibinfo{title}{{PETS}c {W}eb page}.
\newblock \bibinfo{howpublished}{\url{http://www.mcs.anl.gov/petsc}}.
\newblock
\urldef\tempurl%
\url{http://www.mcs.anl.gov/petsc}
\showURL{%
\tempurl}


\bibitem[\protect\citeauthoryear{Balay, Abhyankar, Adams, Brown, Brune,
  Buschelman, Eijkhout, Gropp, Kaushik, Knepley, McInnes, Rupp, Smith, and
  Zhang}{Balay et~al\mbox{.}}{2014}]%
        {BalayAbhyankarAdamsEtAl14}
\bibfield{author}{\bibinfo{person}{Satish Balay}, \bibinfo{person}{Shrirang
  Abhyankar}, \bibinfo{person}{Mark~F. Adams}, \bibinfo{person}{Jed Brown},
  \bibinfo{person}{Peter Brune}, \bibinfo{person}{Kris Buschelman},
  \bibinfo{person}{Victor Eijkhout}, \bibinfo{person}{William~D. Gropp},
  \bibinfo{person}{Dinesh Kaushik}, \bibinfo{person}{Matthew~G. Knepley},
  \bibinfo{person}{Lois~Curfman McInnes}, \bibinfo{person}{Karl Rupp},
  \bibinfo{person}{Barry~F. Smith}, {and} \bibinfo{person}{Hong Zhang}.}
  \bibinfo{year}{2014}\natexlab{}.
\newblock \bibinfo{title}{{PETS}c {W}eb page}.
\newblock \bibinfo{howpublished}{\url{http://www.mcs.anl.gov/petsc}}.
\newblock
\urldef\tempurl%
\url{http://www.mcs.anl.gov/petsc}
\showURL{%
\tempurl}


\bibitem[\protect\citeauthoryear{Bardsley, Cui, Marzouk, and Wang}{Bardsley
  et~al\mbox{.}}{2020}]%
        {bardsley2020scalable}
\bibfield{author}{\bibinfo{person}{Johnathan~M Bardsley},
  \bibinfo{person}{Tiangang Cui}, \bibinfo{person}{Youssef~M Marzouk}, {and}
  \bibinfo{person}{Zheng Wang}.} \bibinfo{year}{2020}\natexlab{}.
\newblock \showarticletitle{Scalable optimization-based sampling on function
  space}.
\newblock \bibinfo{journal}{\emph{SIAM Journal on Scientific Computing}}
  \bibinfo{volume}{42}, \bibinfo{number}{2} (\bibinfo{year}{2020}),
  \bibinfo{pages}{A1317--A1347}.
\newblock


\bibitem[\protect\citeauthoryear{Becker, Carey, and Oden}{Becker
  et~al\mbox{.}}{1981}]%
        {BeckerCareyOden81}
\bibfield{author}{\bibinfo{person}{E.~B. Becker}, \bibinfo{person}{G.~F.
  Carey}, {and} \bibinfo{person}{J.~T. Oden}.} \bibinfo{year}{1981}\natexlab{}.
\newblock \bibinfo{booktitle}{\emph{Finite Elements: An Introduction, Vol I}}.
\newblock \bibinfo{publisher}{Prentice Hall}, \bibinfo{address}{Englewoods
  Cliffs, New Jersey}.
\newblock


\bibitem[\protect\citeauthoryear{Beskos, Girolami, Lan, Farrell, and
  Stuart}{Beskos et~al\mbox{.}}{2017}]%
        {BeskosGirolamiLanEtAl17}
\bibfield{author}{\bibinfo{person}{Alexandros Beskos}, \bibinfo{person}{Mark
  Girolami}, \bibinfo{person}{Shiwei Lan}, \bibinfo{person}{Patrick~E Farrell},
  {and} \bibinfo{person}{Andrew~M Stuart}.} \bibinfo{year}{2017}\natexlab{}.
\newblock \showarticletitle{Geometric {MCMC} for infinite-dimensional inverse
  problems}.
\newblock \bibinfo{journal}{\emph{J. Comput. Phys.}}  \bibinfo{volume}{335}
  (\bibinfo{year}{2017}), \bibinfo{pages}{327--351}.
\newblock


\bibitem[\protect\citeauthoryear{Borz{\`\i} and Schulz}{Borz{\`\i} and
  Schulz}{2012}]%
        {BorziSchulz12}
\bibfield{author}{\bibinfo{person}{Alfio Borz{\`\i}} {and}
  \bibinfo{person}{Volker Schulz}.} \bibinfo{year}{2012}\natexlab{}.
\newblock \bibinfo{booktitle}{\emph{Computational Optimization of Systems
  Governed by Partial Differential Equations}}.
\newblock \bibinfo{publisher}{SIAM}.
\newblock


\bibitem[\protect\citeauthoryear{Brooks and Gelman}{Brooks and Gelman}{1998}]%
        {Brooks1998}
\bibfield{author}{\bibinfo{person}{Stephen~P Brooks} {and}
  \bibinfo{person}{Andrew Gelman}.} \bibinfo{year}{1998}\natexlab{}.
\newblock \showarticletitle{{General Methods for Monitoring Convergence of
  Iterative Simulations}}.
\newblock \bibinfo{journal}{\emph{Journal of Computational and Graphical
  Statistics}} \bibinfo{volume}{7}, \bibinfo{number}{4} (\bibinfo{date}{dec}
  \bibinfo{year}{1998}), \bibinfo{pages}{434--455}.
\newblock
\showISSN{1061-8600}
\urldef\tempurl%
\url{https://doi.org/10.1080/10618600.1998.10474787}
\showDOI{\tempurl}


\bibitem[\protect\citeauthoryear{Brown}{Brown}{2010}]%
        {Brown2010}
\bibfield{author}{\bibinfo{person}{Jed Brown}.}
  \bibinfo{year}{2010}\natexlab{}.
\newblock \showarticletitle{{Efficient Nonlinear Solvers for Nodal High-Order
  Finite Elements in 3D}}.
\newblock \bibinfo{journal}{\emph{Journal of Scientific Computing}}
  \bibinfo{volume}{45}, \bibinfo{number}{1} (\bibinfo{year}{2010}),
  \bibinfo{pages}{48--63}.
\newblock
\showISSN{1573-7691}
\urldef\tempurl%
\url{https://doi.org/10.1007/s10915-010-9396-8}
\showDOI{\tempurl}


\bibitem[\protect\citeauthoryear{Bui-Thanh, Burstedde, Ghattas, Martin,
  Stadler, and Wilcox}{Bui-Thanh et~al\mbox{.}}{2012}]%
        {Bui-ThanhBursteddeGhattasEtAl12_gbfinalist}
\bibfield{author}{\bibinfo{person}{Tan Bui-Thanh}, \bibinfo{person}{Carsten
  Burstedde}, \bibinfo{person}{Omar Ghattas}, \bibinfo{person}{James Martin},
  \bibinfo{person}{Georg Stadler}, {and} \bibinfo{person}{Lucas~C. Wilcox}.}
  \bibinfo{year}{2012}\natexlab{}.
\newblock \showarticletitle{{Extreme-scale UQ for Bayesian inverse problems
  governed by PDEs}}. In \bibinfo{booktitle}{\emph{SC12: Proceedings of the
  International Conference for High Performance Computing, Networking, Storage
  and Analysis}}.
\newblock
\newblock
\shownote{{G}ordon {B}ell {P}rize finalist.}


\bibitem[\protect\citeauthoryear{Bui-Thanh, Ghattas, Martin, and
  Stadler}{Bui-Thanh et~al\mbox{.}}{2013}]%
        {Bui-ThanhGhattasMartinEtAl13a}
\bibfield{author}{\bibinfo{person}{T. Bui-Thanh}, \bibinfo{person}{O. Ghattas},
  \bibinfo{person}{J. Martin}, {and} \bibinfo{person}{G. Stadler}.}
  \bibinfo{year}{2013}\natexlab{}.
\newblock \showarticletitle{A computational framework for infinite-dimensional
  Bayesian inverse problems Part {I}: The linearized case, with application to
  global seismic inversion}.
\newblock \bibinfo{journal}{\emph{SIAM Journal on Scientific Computing}}
  \bibinfo{volume}{35}, \bibinfo{number}{6} (\bibinfo{year}{2013}),
  \bibinfo{pages}{A2494--A2523}.
\newblock


\bibitem[\protect\citeauthoryear{Calderhead}{Calderhead}{2014}]%
        {calderhead2014general}
\bibfield{author}{\bibinfo{person}{Ben Calderhead}.}
  \bibinfo{year}{2014}\natexlab{}.
\newblock \showarticletitle{A general construction for parallelizing
  Metropolis- Hastings algorithms}.
\newblock \bibinfo{journal}{\emph{Proceedings of the National Academy of
  Sciences}} \bibinfo{volume}{111}, \bibinfo{number}{49}
  (\bibinfo{year}{2014}), \bibinfo{pages}{17408--17413}.
\newblock


\bibitem[\protect\citeauthoryear{Casella and George}{Casella and
  George}{1992}]%
        {CasellaGeorge92}
\bibfield{author}{\bibinfo{person}{George Casella} {and}
  \bibinfo{person}{Edward~I. George}.} \bibinfo{year}{1992}\natexlab{}.
\newblock \showarticletitle{Explaining the {G}ibbs sampler}.
\newblock \bibinfo{journal}{\emph{The American Statistician}}
  \bibinfo{volume}{46}, \bibinfo{number}{3} (\bibinfo{year}{1992}),
  \bibinfo{pages}{167--174}.
\newblock


\bibitem[\protect\citeauthoryear{Conrad and Marzouk}{Conrad and
  Marzouk}{2013}]%
        {Conrad2013}
\bibfield{author}{\bibinfo{person}{Patrick~R Conrad} {and}
  \bibinfo{person}{Youssef~M Marzouk}.} \bibinfo{year}{2013}\natexlab{}.
\newblock \showarticletitle{Adaptive {S}molyak pseudospectral approximations}.
\newblock \bibinfo{journal}{\emph{SIAM Journal on Scientific Computing}}
  \bibinfo{volume}{35}, \bibinfo{number}{6} (\bibinfo{year}{2013}),
  \bibinfo{pages}{A2643--A2670}.
\newblock
\urldef\tempurl%
\url{https://doi.org/10.1137/120890715}
\showDOI{\tempurl}


\bibitem[\protect\citeauthoryear{Cotter, Roberts, Stuart, and White}{Cotter
  et~al\mbox{.}}{2012}]%
        {CotterRobertsStuartEtAl12}
\bibfield{author}{\bibinfo{person}{S.~L. Cotter}, \bibinfo{person}{G.~O.
  Roberts}, \bibinfo{person}{A.~M. Stuart}, {and} \bibinfo{person}{D. White}.}
  \bibinfo{year}{2012}\natexlab{}.
\newblock \showarticletitle{{MCMC} methods for functions: modifying old
  algorithms to make them faster}.
\newblock  (\bibinfo{year}{2012}).
\newblock
\newblock
\shownote{submitted.}


\bibitem[\protect\citeauthoryear{Cui, Law, and Marzouk}{Cui
  et~al\mbox{.}}{2016a}]%
        {CuiLawMarzouk16}
\bibfield{author}{\bibinfo{person}{T. Cui}, \bibinfo{person}{K.J.H. Law}, {and}
  \bibinfo{person}{Y.M. Marzouk}.} \bibinfo{year}{2016}\natexlab{a}.
\newblock \showarticletitle{Dimension-independent likelihood-informed {MCMC}}.
\newblock \bibinfo{journal}{\emph{J. Comput. Phys.}}  \bibinfo{volume}{304}
  (\bibinfo{year}{2016}), \bibinfo{pages}{109--137}.
\newblock


\bibitem[\protect\citeauthoryear{Cui, Marzouk, and Willcox}{Cui
  et~al\mbox{.}}{2016b}]%
        {CuiMarzoukWillcox16}
\bibfield{author}{\bibinfo{person}{Tiangang Cui}, \bibinfo{person}{Youssef
  Marzouk}, {and} \bibinfo{person}{Karen Willcox}.}
  \bibinfo{year}{2016}\natexlab{b}.
\newblock \showarticletitle{Scalable posterior approximations for large-scale
  Bayesian inverse problems via likelihood-informed parameter and state
  reduction}.
\newblock \bibinfo{journal}{\emph{J. Comput. Phys.}}  \bibinfo{volume}{315}
  (\bibinfo{year}{2016}), \bibinfo{pages}{363--387}.
\newblock


\bibitem[\protect\citeauthoryear{Cui and Zahm}{Cui and Zahm}{2021}]%
        {cui2021data}
\bibfield{author}{\bibinfo{person}{Tiangang Cui} {and} \bibinfo{person}{Olivier
  Zahm}.} \bibinfo{year}{2021}\natexlab{}.
\newblock \showarticletitle{Data-free likelihood-informed dimension reduction
  of Bayesian inverse problems}.
\newblock \bibinfo{journal}{\emph{Inverse Problems}} \bibinfo{volume}{37},
  \bibinfo{number}{4} (\bibinfo{year}{2021}), \bibinfo{pages}{045009}.
\newblock


\bibitem[\protect\citeauthoryear{Daon and Stadler}{Daon and Stadler}{2018}]%
        {DaonStadler18}
\bibfield{author}{\bibinfo{person}{Yair Daon} {and} \bibinfo{person}{Georg
  Stadler}.} \bibinfo{year}{2018}\natexlab{}.
\newblock \showarticletitle{Mitigating the Influence of Boundary Conditions on
  Covariance Operators Derived from Elliptic {PDEs}}.
\newblock \bibinfo{journal}{\emph{Inverse Problems and Imaging}}
  \bibinfo{volume}{12}, \bibinfo{number}{5} (\bibinfo{year}{2018}),
  \bibinfo{pages}{1083--1102}.
\newblock
\showeprint{1610.05280}


\bibitem[\protect\citeauthoryear{Dodwell, Ketelsen, Scheichl, and
  Teckentrup}{Dodwell et~al\mbox{.}}{2019}]%
        {dodwell2019multilevel}
\bibfield{author}{\bibinfo{person}{Tim~J. Dodwell}, \bibinfo{person}{Christian
  Ketelsen}, \bibinfo{person}{Robert Scheichl}, {and}
  \bibinfo{person}{Aretha~L. Teckentrup}.} \bibinfo{year}{2019}\natexlab{}.
\newblock \showarticletitle{Multilevel Markov chain Monte Carlo}.
\newblock \bibinfo{journal}{\emph{SIAM Rev.}} \bibinfo{volume}{61},
  \bibinfo{number}{3} (\bibinfo{year}{2019}), \bibinfo{pages}{509--545}.
\newblock


\bibitem[\protect\citeauthoryear{Evans and Swartz}{Evans and Swartz}{2000}]%
        {EvansSwartz00}
\bibfield{author}{\bibinfo{person}{M. Evans} {and} \bibinfo{person}{T.
  Swartz}.} \bibinfo{year}{2000}\natexlab{}.
\newblock \bibinfo{booktitle}{\emph{{A}pproximating integrals via {M}onte
  {C}arlo and deterministic methods}}. Vol.~\bibinfo{volume}{20}.
\newblock \bibinfo{publisher}{OUP Oxford}.
\newblock


\bibitem[\protect\citeauthoryear{Flegal and Jones}{Flegal and Jones}{2010}]%
        {flegal2010batch}
\bibfield{author}{\bibinfo{person}{James~M Flegal} {and}
  \bibinfo{person}{Galin~L Jones}.} \bibinfo{year}{2010}\natexlab{}.
\newblock \showarticletitle{Batch means and spectral variance estimators in
  Markov chain Monte Carlo}.
\newblock \bibinfo{journal}{\emph{The Annals of Statistics}}
  \bibinfo{volume}{38}, \bibinfo{number}{2} (\bibinfo{year}{2010}),
  \bibinfo{pages}{1034--1070}.
\newblock


\bibitem[\protect\citeauthoryear{Gelman, Carlin, Stern, and Rubin}{Gelman
  et~al\mbox{.}}{2004}]%
        {Gelman2004}
\bibfield{author}{\bibinfo{person}{Andrew Gelman}, \bibinfo{person}{John~B
  Carlin}, \bibinfo{person}{Hal~S Stern}, {and} \bibinfo{person}{Donald~B
  Rubin}.} \bibinfo{year}{2004}\natexlab{}.
\newblock \bibinfo{title}{{Bayesian data analysis}}.
\newblock
\newblock


\bibitem[\protect\citeauthoryear{Gelman and Rubin}{Gelman and Rubin}{1992}]%
        {gelman1992inference}
\bibfield{author}{\bibinfo{person}{Andrew Gelman} {and}
  \bibinfo{person}{Donald~B Rubin}.} \bibinfo{year}{1992}\natexlab{}.
\newblock \showarticletitle{Inference from iterative simulation using multiple
  sequences}.
\newblock \bibinfo{journal}{\emph{Statistical science}} (\bibinfo{year}{1992}),
  \bibinfo{pages}{457--472}.
\newblock


\bibitem[\protect\citeauthoryear{Ghattas and Willcox}{Ghattas and
  Willcox}{2021}]%
        {GhattasWillcox21}
\bibfield{author}{\bibinfo{person}{O. Ghattas} {and} \bibinfo{person}{K.
  Willcox}.} \bibinfo{year}{2021}\natexlab{}.
\newblock \showarticletitle{Learning physics-based models from data:
  perspectives from inverse problems and model reduction}.
\newblock \bibinfo{journal}{\emph{Acta Numerica}}  \bibinfo{volume}{30}
  (\bibinfo{year}{2021}), \bibinfo{pages}{445--554}.
\newblock
\urldef\tempurl%
\url{https://doi.org/doi:10.1017/S0962492921000064}
\showDOI{\tempurl}


\bibitem[\protect\citeauthoryear{Girolami and Calderhead}{Girolami and
  Calderhead}{2011}]%
        {girolami2011riemann}
\bibfield{author}{\bibinfo{person}{Mark Girolami} {and} \bibinfo{person}{Ben
  Calderhead}.} \bibinfo{year}{2011}\natexlab{}.
\newblock \showarticletitle{Riemann manifold {L}angevin and {H}amiltonian
  {M}onte {C}arlo methods}.
\newblock \bibinfo{journal}{\emph{Journal of the Royal Statistical Society:
  Series B (Statistical Methodology)}} \bibinfo{volume}{73},
  \bibinfo{number}{2} (\bibinfo{year}{2011}), \bibinfo{pages}{123--214}.
\newblock


\bibitem[\protect\citeauthoryear{Golub and Van~Loan}{Golub and
  Van~Loan}{1996}]%
        {GolubVan96}
\bibfield{author}{\bibinfo{person}{Gene~H. Golub} {and}
  \bibinfo{person}{Charles~F. Van~Loan}.} \bibinfo{year}{1996}\natexlab{}.
\newblock \bibinfo{booktitle}{\emph{Matrix Computations}
  (\bibinfo{edition}{third} ed.)}.
\newblock \bibinfo{publisher}{Johns Hopkins University Press},
  \bibinfo{address}{Baltimore, MD}.
\newblock


\bibitem[\protect\citeauthoryear{Haario, Saksman, and Tamminen}{Haario
  et~al\mbox{.}}{2001}]%
        {Haario2001}
\bibfield{author}{\bibinfo{person}{Heikki Haario}, \bibinfo{person}{Eero
  Saksman}, {and} \bibinfo{person}{Johanna Tamminen}.}
  \bibinfo{year}{2001}\natexlab{}.
\newblock \showarticletitle{{An Adaptive Metropolis Algorithm}}.
\newblock \bibinfo{journal}{\emph{Bernoulli}} \bibinfo{volume}{7},
  \bibinfo{number}{2} (\bibinfo{date}{sep} \bibinfo{year}{2001}),
  \bibinfo{pages}{223--242}.
\newblock
\showISSN{13507265}
\urldef\tempurl%
\url{https://doi.org/10.2307/3318737}
\showDOI{\tempurl}


\bibitem[\protect\citeauthoryear{Hairer, Stuart, and Vollmer}{Hairer
  et~al\mbox{.}}{2014}]%
        {hairer2014spectral}
\bibfield{author}{\bibinfo{person}{Martin Hairer}, \bibinfo{person}{Andrew~M
  Stuart}, {and} \bibinfo{person}{Sebastian~J Vollmer}.}
  \bibinfo{year}{2014}\natexlab{}.
\newblock \showarticletitle{Spectral gaps for a {M}etropolis--{H}astings
  algorithm in infinite dimensions}.
\newblock \bibinfo{journal}{\emph{The Annals of Applied Probability}}
  \bibinfo{volume}{24}, \bibinfo{number}{6} (\bibinfo{year}{2014}),
  \bibinfo{pages}{2455--2490}.
\newblock


\bibitem[\protect\citeauthoryear{Halko, Martinsson, and Tropp}{Halko
  et~al\mbox{.}}{2011}]%
        {HalkoMartinssonTropp11}
\bibfield{author}{\bibinfo{person}{Nathan Halko}, \bibinfo{person}{Per~Gunnar
  Martinsson}, {and} \bibinfo{person}{Joel~A. Tropp}.}
  \bibinfo{year}{2011}\natexlab{}.
\newblock \showarticletitle{Finding structure with randomness: {P}robabilistic
  algorithms for constructing approximate matrix decompositions}.
\newblock \bibinfo{journal}{\emph{SIAM Rev.}} \bibinfo{volume}{53},
  \bibinfo{number}{2} (\bibinfo{year}{2011}), \bibinfo{pages}{217--288}.
\newblock


\bibitem[\protect\citeauthoryear{Hartikainen and S{\"a}rkk{\"a}}{Hartikainen
  and S{\"a}rkk{\"a}}{2010}]%
        {Hartikainen2010}
\bibfield{author}{\bibinfo{person}{Jouni Hartikainen} {and}
  \bibinfo{person}{Simo S{\"a}rkk{\"a}}.} \bibinfo{year}{2010}\natexlab{}.
\newblock \showarticletitle{Kalman filtering and smoothing solutions to
  temporal {G}aussian process regression models}. In
  \bibinfo{booktitle}{\emph{2010 {IEEE} international workshop on machine
  learning for signal processing}}. IEEE, \bibinfo{pages}{379--384}.
\newblock
\urldef\tempurl%
\url{https://doi.org/10.1109/MLSP.2010.5589113}
\showDOI{\tempurl}


\bibitem[\protect\citeauthoryear{Hastings}{Hastings}{1970}]%
        {Hastings70}
\bibfield{author}{\bibinfo{person}{W.~Keith Hastings}.}
  \bibinfo{year}{1970}\natexlab{}.
\newblock \showarticletitle{Monte {C}arlo sampling methods using {M}arkov
  chains and their applications}.
\newblock \bibinfo{journal}{\emph{Biometrika}} \bibinfo{volume}{57},
  \bibinfo{number}{1} (\bibinfo{year}{1970}), \bibinfo{pages}{97--109}.
\newblock


\bibitem[\protect\citeauthoryear{Hoffman and Gelman}{Hoffman and
  Gelman}{2014}]%
        {hoffman2014no}
\bibfield{author}{\bibinfo{person}{Matthew~D Hoffman} {and}
  \bibinfo{person}{Andrew Gelman}.} \bibinfo{year}{2014}\natexlab{}.
\newblock \showarticletitle{The No-U-Turn sampler: adaptively setting path
  lengths in Hamiltonian Monte Carlo.}
\newblock \bibinfo{journal}{\emph{Journal of Machine Learning Research}}
  \bibinfo{volume}{15}, \bibinfo{number}{1} (\bibinfo{year}{2014}),
  \bibinfo{pages}{1593--1623}.
\newblock


\bibitem[\protect\citeauthoryear{Isaac, Petra, Stadler, and Ghattas}{Isaac
  et~al\mbox{.}}{2015}]%
        {IsaacPetraStadlerEtAl15a}
\bibfield{author}{\bibinfo{person}{Tobin Isaac}, \bibinfo{person}{Noemi Petra},
  \bibinfo{person}{Georg Stadler}, {and} \bibinfo{person}{Omar Ghattas}.}
  \bibinfo{year}{2015}\natexlab{}.
\newblock \showarticletitle{Scalable and efficient algorithms for the
  propagation of uncertainty from data through inference to prediction for
  large-scale problems, with application to flow of the {A}ntarctic ice sheet}.
\newblock \bibinfo{journal}{\emph{J. Comput. Phys.}}  \bibinfo{volume}{296}
  (\bibinfo{date}{September} \bibinfo{year}{2015}), \bibinfo{pages}{348--368}.
\newblock
\urldef\tempurl%
\url{https://doi.org/10.1016/j.jcp.2015.04.047}
\showDOI{\tempurl}


\bibitem[\protect\citeauthoryear{Kaipio and Somersalo}{Kaipio and
  Somersalo}{2005}]%
        {KaipioSomersalo05}
\bibfield{author}{\bibinfo{person}{Jari Kaipio} {and} \bibinfo{person}{Erkki
  Somersalo}.} \bibinfo{year}{2005}\natexlab{}.
\newblock \bibinfo{booktitle}{\emph{Statistical and Computational Inverse
  Problems}}. \bibinfo{series}{Applied Mathematical Sciences},
  Vol.~\bibinfo{volume}{160}.
\newblock \bibinfo{publisher}{Springer-Verlag New York}.
\newblock
\urldef\tempurl%
\url{https://doi.org/10.1007/b138659}
\showDOI{\tempurl}


\bibitem[\protect\citeauthoryear{Lindgren, Rue, and Lindstr{\"o}m}{Lindgren
  et~al\mbox{.}}{2011}]%
        {LindgrenRueLindstroem11}
\bibfield{author}{\bibinfo{person}{Finn Lindgren}, \bibinfo{person}{H{\aa}vard
  Rue}, {and} \bibinfo{person}{Johan Lindstr{\"o}m}.}
  \bibinfo{year}{2011}\natexlab{}.
\newblock \showarticletitle{An explicit link between {G}aussian fields and
  {G}aussian {M}arkov random fields: the stochastic partial differential
  equation approach}.
\newblock \bibinfo{journal}{\emph{Journal of the Royal Statistical Society:
  Series B (Statistical Methodology)}} \bibinfo{volume}{73},
  \bibinfo{number}{4} (\bibinfo{year}{2011}), \bibinfo{pages}{423--498}.
\newblock
\showISSN{1467-9868}
\urldef\tempurl%
\url{https://doi.org/10.1111/j.1467-9868.2011.00777.x}
\showDOI{\tempurl}


\bibitem[\protect\citeauthoryear{Lindqvist}{Lindqvist}{2017}]%
        {lindqvist2017notes}
\bibfield{author}{\bibinfo{person}{Peter Lindqvist}.}
  \bibinfo{year}{2017}\natexlab{}.
\newblock \bibinfo{booktitle}{\emph{{Notes on the p-Laplace equation}}}.
\newblock Number 161. \bibinfo{publisher}{University of Jyv{\"{a}}skyl{\"{a}}}.
\newblock


\bibitem[\protect\citeauthoryear{Logg, Mardal, and Wells}{Logg
  et~al\mbox{.}}{2012}]%
        {LoggMardalWells12}
\bibfield{editor}{\bibinfo{person}{Anders Logg}, \bibinfo{person}{Kent-Andre
  Mardal}, {and} \bibinfo{person}{Garth~N. Wells}} (Eds.).
  \bibinfo{year}{2012}\natexlab{}.
\newblock \bibinfo{booktitle}{\emph{Automated Solution of Differential
  Equations by the Finite Element Method}}. \bibinfo{series}{Lecture Notes in
  Computational Science and Engineering}, Vol.~\bibinfo{volume}{84}.
\newblock \bibinfo{publisher}{Springer}.
\newblock
\urldef\tempurl%
\url{https://doi.org/10.1007/978-3-642-23099-8}
\showDOI{\tempurl}


\bibitem[\protect\citeauthoryear{Marshall and Roberts}{Marshall and
  Roberts}{2012}]%
        {MarshallRoberts12}
\bibfield{author}{\bibinfo{person}{Tristan Marshall} {and}
  \bibinfo{person}{Gareth Roberts}.} \bibinfo{year}{2012}\natexlab{}.
\newblock \showarticletitle{An Adaptive Approach to Langevin MCMC}.
\newblock \bibinfo{journal}{\emph{Statistics and Computing}}
  \bibinfo{volume}{22}, \bibinfo{number}{5} (\bibinfo{date}{Sept.}
  \bibinfo{year}{2012}), \bibinfo{pages}{1041–1057}.
\newblock
\showISSN{0960-3174}
\urldef\tempurl%
\url{https://doi.org/10.1007/s11222-011-9276-6}
\showDOI{\tempurl}


\bibitem[\protect\citeauthoryear{Martin, Wilcox, Burstedde, and Ghattas}{Martin
  et~al\mbox{.}}{2012}]%
        {martin2012stochastic}
\bibfield{author}{\bibinfo{person}{James Martin}, \bibinfo{person}{Lucas~C
  Wilcox}, \bibinfo{person}{Carsten Burstedde}, {and} \bibinfo{person}{Omar
  Ghattas}.} \bibinfo{year}{2012}\natexlab{}.
\newblock \showarticletitle{A stochastic Newton MCMC method for large-scale
  statistical inverse problems with application to seismic inversion}.
\newblock \bibinfo{journal}{\emph{SIAM Journal on Scientific Computing}}
  \bibinfo{volume}{34}, \bibinfo{number}{3} (\bibinfo{year}{2012}),
  \bibinfo{pages}{A1460--A1487}.
\newblock


\bibitem[\protect\citeauthoryear{Marzouk, Moselhy, Parno, and Spantini}{Marzouk
  et~al\mbox{.}}{2016}]%
        {Marzouk2016}
\bibfield{author}{\bibinfo{person}{Youssef Marzouk}, \bibinfo{person}{Tarek
  Moselhy}, \bibinfo{person}{Matthew Parno}, {and} \bibinfo{person}{Alessio
  Spantini}.} \bibinfo{year}{2016}\natexlab{}.
\newblock \bibinfo{booktitle}{\emph{Sampling via Measure Transport: An
  Introduction}}.
\newblock \bibinfo{publisher}{Springer International Publishing},
  \bibinfo{pages}{1--41}.
\newblock
\showISBNx{978-3-319-11259-6}
\urldef\tempurl%
\url{https://doi.org/10.1007/978-3-319-11259-6_23-1}
\showDOI{\tempurl}


\bibitem[\protect\citeauthoryear{Merkel}{Merkel}{2014}]%
        {Merkel14}
\bibfield{author}{\bibinfo{person}{Dirk Merkel}.}
  \bibinfo{year}{2014}\natexlab{}.
\newblock \showarticletitle{Docker: Lightweight Linux Containers for Consistent
  Development and Deployment}.
\newblock \bibinfo{journal}{\emph{Linux J.}} \bibinfo{volume}{2014},
  \bibinfo{number}{239}, Article \bibinfo{articleno}{2} (\bibinfo{year}{2014}).
\newblock
\showISSN{1075-3583}
\urldef\tempurl%
\url{http://dl.acm.org/citation.cfm?id=2600239.2600241}
\showURL{%
\tempurl}


\bibitem[\protect\citeauthoryear{Metropolis, Rosenbluth, Rosenbluth, Teller,
  and Teller}{Metropolis et~al\mbox{.}}{1953}]%
        {MetropolisRosenbluthRosenbluthEtAl53}
\bibfield{author}{\bibinfo{person}{Nicholas Metropolis},
  \bibinfo{person}{Arianna~W. Rosenbluth}, \bibinfo{person}{Marshall~N.
  Rosenbluth}, \bibinfo{person}{Augusta~H. Teller}, {and}
  \bibinfo{person}{Edward Teller}.} \bibinfo{year}{1953}\natexlab{}.
\newblock \showarticletitle{Equation of State Calculations by Fast Computing
  Machines}.
\newblock \bibinfo{journal}{\emph{The Journal of Chemical Physics}}
  \bibinfo{volume}{21}, \bibinfo{number}{6} (\bibinfo{year}{1953}),
  \bibinfo{pages}{1087--1092}.
\newblock
\urldef\tempurl%
\url{https://doi.org/10.1063/1.1699114}
\showDOI{\tempurl}


\bibitem[\protect\citeauthoryear{Mira et~al\mbox{.}}{Mira
  et~al\mbox{.}}{2001}]%
        {Mira2001}
\bibfield{author}{\bibinfo{person}{Antonietta Mira} {et~al\mbox{.}}}
  \bibinfo{year}{2001}\natexlab{}.
\newblock \showarticletitle{{On Metropolis-Hastings algorithms with delayed
  rejection}}.
\newblock \bibinfo{journal}{\emph{Metron}} \bibinfo{volume}{59},
  \bibinfo{number}{3-4} (\bibinfo{year}{2001}), \bibinfo{pages}{231--241}.
\newblock


\bibitem[\protect\citeauthoryear{Neal}{Neal}{2010}]%
        {Neal10}
\bibfield{author}{\bibinfo{person}{R.~M. Neal}.}
  \bibinfo{year}{2010}\natexlab{}.
\newblock \bibinfo{booktitle}{\emph{Handbook of Markov Chain Monte Carlo}}.
\newblock \bibinfo{publisher}{Chapman \& Hall / CRC Press}, Chapter MCMC using
  Hamiltonian dynamics.
\newblock


\bibitem[\protect\citeauthoryear{Owen}{Owen}{2013}]%
        {owen2013monte}
\bibfield{author}{\bibinfo{person}{Art~B Owen}.}
  \bibinfo{year}{2013}\natexlab{}.
\newblock \showarticletitle{Monte Carlo theory, methods and examples}.
\newblock  (\bibinfo{year}{2013}).
\newblock


\bibitem[\protect\citeauthoryear{Parno, Davis, Conrad, and Marzouk}{Parno
  et~al\mbox{.}}{2014}]%
        {Parno2014}
\bibfield{author}{\bibinfo{person}{Matthew Parno}, \bibinfo{person}{Andrew
  Davis}, \bibinfo{person}{Patrick Conrad}, {and} \bibinfo{person}{YM
  Marzouk}.} \bibinfo{year}{2014}\natexlab{}.
\newblock \bibinfo{title}{MIT {U}ncertainty {Q}uantification ({MUQ})
  {L}ibrary}.
\newblock
\newblock
\urldef\tempurl%
\url{https://muq.mit.edu}
\showURL{%
\tempurl}


\bibitem[\protect\citeauthoryear{Parno and Marzouk}{Parno and Marzouk}{2018}]%
        {Parno2018}
\bibfield{author}{\bibinfo{person}{Matthew~D Parno} {and}
  \bibinfo{person}{Youssef~M Marzouk}.} \bibinfo{year}{2018}\natexlab{}.
\newblock \showarticletitle{Transport map accelerated {M}arkov chain {M}onte
  {C}arlo}.
\newblock \bibinfo{journal}{\emph{SIAM/ASA Journal on Uncertainty
  Quantification}} \bibinfo{volume}{6}, \bibinfo{number}{2}
  (\bibinfo{year}{2018}), \bibinfo{pages}{645--682}.
\newblock
\urldef\tempurl%
\url{https://doi.org/10.1137/17M1134640}
\showDOI{\tempurl}


\bibitem[\protect\citeauthoryear{Petra, Martin, Stadler, and Ghattas}{Petra
  et~al\mbox{.}}{2014}]%
        {PetraMartinStadlerEtAl14}
\bibfield{author}{\bibinfo{person}{Noemi Petra}, \bibinfo{person}{James
  Martin}, \bibinfo{person}{Georg Stadler}, {and} \bibinfo{person}{Omar
  Ghattas}.} \bibinfo{year}{2014}\natexlab{}.
\newblock \showarticletitle{A computational framework for infinite-dimensional
  {B}ayesian inverse problems: {P}art {II}. {S}tochastic {N}ewton {MCMC} with
  application to ice sheet inverse problems}.
\newblock \bibinfo{journal}{\emph{SIAM Journal on Scientific Computing}}
  \bibinfo{volume}{36}, \bibinfo{number}{4} (\bibinfo{year}{2014}),
  \bibinfo{pages}{A1525--A1555}.
\newblock


\bibitem[\protect\citeauthoryear{Pinski, Simpson, Stuart, and Weber}{Pinski
  et~al\mbox{.}}{2015}]%
        {PinskiSimpsonStuartEtAl15}
\bibfield{author}{\bibinfo{person}{Frank~J Pinski}, \bibinfo{person}{Gideon
  Simpson}, \bibinfo{person}{Andrew~M Stuart}, {and} \bibinfo{person}{Hendrik
  Weber}.} \bibinfo{year}{2015}\natexlab{}.
\newblock \showarticletitle{Algorithms for {K}ullback--{L}eibler approximation
  of probability measures in infinite dimensions}.
\newblock \bibinfo{journal}{\emph{SIAM Journal on Scientific Computing}}
  \bibinfo{volume}{37}, \bibinfo{number}{6} (\bibinfo{year}{2015}),
  \bibinfo{pages}{A2733--A2757}.
\newblock


\bibitem[\protect\citeauthoryear{Press}{Press}{2003}]%
        {Press03}
\bibfield{author}{\bibinfo{person}{S.~J. Press}.}
  \bibinfo{year}{2003}\natexlab{}.
\newblock \bibinfo{booktitle}{\emph{{S}ubjective and {O}bjective {B}ayesian
  {S}tatistics: {P}rinciples, {M}ethods and {A}pplications}}.
\newblock \bibinfo{publisher}{Wiley, New York}.
\newblock


\bibitem[\protect\citeauthoryear{Rasmussen and Williams}{Rasmussen and
  Williams}{2005}]%
        {GPML2005}
\bibfield{author}{\bibinfo{person}{Carl~Edward Rasmussen} {and}
  \bibinfo{person}{Christopher K.~I. Williams}.}
  \bibinfo{year}{2005}\natexlab{}.
\newblock \bibinfo{booktitle}{\emph{Gaussian Processes for Machine Learning}}.
\newblock \bibinfo{publisher}{The MIT Press}.
\newblock
\showISBNx{026218253X}


\bibitem[\protect\citeauthoryear{Robert and Casella}{Robert and
  Casella}{2005}]%
        {RobertCasella05}
\bibfield{author}{\bibinfo{person}{Christian~P. Robert} {and}
  \bibinfo{person}{George Casella}.} \bibinfo{year}{2005}\natexlab{}.
\newblock \bibinfo{booktitle}{\emph{Monte Carlo Statistical Methods (Springer
  Texts in Statistics)}}.
\newblock \bibinfo{publisher}{Springer-Verlag New York, Inc.},
  \bibinfo{address}{Secaucus, NJ, USA}.
\newblock
\showISBNx{0387212396}


\bibitem[\protect\citeauthoryear{Roberts, Rosenthal, et~al\mbox{.}}{Roberts
  et~al\mbox{.}}{2004}]%
        {roberts2004general}
\bibfield{author}{\bibinfo{person}{Gareth~O Roberts},
  \bibinfo{person}{Jeffrey~S Rosenthal}, {et~al\mbox{.}}}
  \bibinfo{year}{2004}\natexlab{}.
\newblock \showarticletitle{{General state space Markov chains and MCMC
  algorithms}}.
\newblock \bibinfo{journal}{\emph{Probability surveys}}  \bibinfo{volume}{1}
  (\bibinfo{year}{2004}), \bibinfo{pages}{20--71}.
\newblock


\bibitem[\protect\citeauthoryear{Roberts and Stramer}{Roberts and
  Stramer}{2003}]%
        {RobertsStramer03}
\bibfield{author}{\bibinfo{person}{Gareth~O. Roberts} {and}
  \bibinfo{person}{Osnat Stramer}.} \bibinfo{year}{2003}\natexlab{}.
\newblock \showarticletitle{Langevin Diffussions and {M}etropolis-{H}astings
  Algorithms}.
\newblock \bibinfo{journal}{\emph{Methodology and Computing in Applied
  Probability}}  \bibinfo{volume}{4} (\bibinfo{year}{2003}),
  \bibinfo{pages}{337--357}.
\newblock


\bibitem[\protect\citeauthoryear{Rudolph and Sprungk}{Rudolph and
  Sprungk}{2018}]%
        {RudolphSprungk16}
\bibfield{author}{\bibinfo{person}{D. Rudolph} {and} \bibinfo{person}{B.
  Sprungk}.} \bibinfo{year}{2018}\natexlab{}.
\newblock \showarticletitle{On a Generalization of the Preconditioned
  {C}rank-{N}icolson {M}etropolis Algorithm}.
\newblock \bibinfo{journal}{\emph{Foundations of Computational Mathematics}}
  \bibinfo{volume}{18} (\bibinfo{year}{2018}), \bibinfo{pages}{309--343}.
\newblock
Issue 2.


\bibitem[\protect\citeauthoryear{{Stigler, S.~M.}}{{Stigler, S.~M.}}{1986}]%
        {Stigler86}
\bibfield{author}{\bibinfo{person}{{Stigler, S.~M.}}}
  \bibinfo{year}{1986}\natexlab{}.
\newblock \showarticletitle{{Laplace's 1774 Memoir on Inverse Probability}}.
\newblock \bibinfo{journal}{\emph{Statist. Sci.}} \bibinfo{volume}{1},
  \bibinfo{number}{3} (\bibinfo{date}{08} \bibinfo{year}{1986}),
  \bibinfo{pages}{359--363}.
\newblock
\urldef\tempurl%
\url{https://doi.org/10.1214/ss/1177013620}
\showDOI{\tempurl}


\bibitem[\protect\citeauthoryear{Strang and Fix}{Strang and Fix}{1988}]%
        {StrangFix88}
\bibfield{author}{\bibinfo{person}{G. Strang} {and} \bibinfo{person}{G.~J.
  Fix}.} \bibinfo{year}{1988}\natexlab{}.
\newblock \bibinfo{booktitle}{\emph{An Analysis of the Finite Element Method}}.
\newblock \bibinfo{publisher}{Wellesley-Cambridge Press},
  \bibinfo{address}{Wellesley, MA}.
\newblock


\bibitem[\protect\citeauthoryear{Stuart}{Stuart}{2010}]%
        {Stuart10}
\bibfield{author}{\bibinfo{person}{Andrew~M. Stuart}.}
  \bibinfo{year}{2010}\natexlab{}.
\newblock \showarticletitle{Inverse problems: {A B}ayesian perspective}.
\newblock \bibinfo{journal}{\emph{Acta Numerica}}  \bibinfo{volume}{19}
  (\bibinfo{year}{2010}), \bibinfo{pages}{451--559}.
\newblock
\urldef\tempurl%
\url{https://doi.org/10.1017/S0962492910000061}
\showDOI{\tempurl}


\bibitem[\protect\citeauthoryear{Tarantola}{Tarantola}{2005}]%
        {Tarantola05}
\bibfield{author}{\bibinfo{person}{Albert Tarantola}.}
  \bibinfo{year}{2005}\natexlab{}.
\newblock \bibinfo{booktitle}{\emph{Inverse Problem Theory and Methods for
  Model Parameter Estimation}}.
\newblock \bibinfo{publisher}{SIAM}, \bibinfo{address}{Philadelphia, PA}.
  xii+342 pages.
\newblock


\bibitem[\protect\citeauthoryear{Tierney and Kadane}{Tierney and
  Kadane}{1986}]%
        {TierneyKadane86}
\bibfield{author}{\bibinfo{person}{L. Tierney} {and} \bibinfo{person}{J.~B.
  Kadane}.} \bibinfo{year}{1986}\natexlab{}.
\newblock \showarticletitle{{Accurate Approximations for Posterior Moments and
  Marginal Densities}}.
\newblock \bibinfo{journal}{\emph{J. Amer. Statist. Assoc.}}
  \bibinfo{volume}{81}, \bibinfo{number}{393} (\bibinfo{year}{1986}),
  \bibinfo{pages}{82--86}.
\newblock
\urldef\tempurl%
\url{https://doi.org/10.1080/01621459.1986.10478240}
\showDOI{\tempurl}


\bibitem[\protect\citeauthoryear{{T}rilinos~{P}roject
  {T}eam}{{T}rilinos~{P}roject {T}eam}{2020}]%
        {trilinos-website}
\bibfield{author}{\bibinfo{person}{The {T}rilinos~{P}roject {T}eam}.}
  \bibinfo{year}{2020 (acccessed May 22, 2020)}\natexlab{}.
\newblock \bibinfo{booktitle}{\emph{The {T}rilinos {P}roject {W}ebsite}}.
\newblock
\urldef\tempurl%
\url{https://trilinos.github.io}
\showURL{%
\tempurl}


\bibitem[\protect\citeauthoryear{Tr\"oltzsch}{Tr\"oltzsch}{2010}]%
        {Troltzsch10}
\bibfield{author}{\bibinfo{person}{Fredi Tr\"oltzsch}.}
  \bibinfo{year}{2010}\natexlab{}.
\newblock \bibinfo{booktitle}{\emph{Optimal Control of Partial Differential
  Equations: Theory, Methods and Applications}}. \bibinfo{series}{Graduate
  Studies in Mathematics}, Vol.~\bibinfo{volume}{112}.
\newblock \bibinfo{publisher}{American Mathematical Society}.
\newblock


\bibitem[\protect\citeauthoryear{Vats, Flegal, and Jones}{Vats
  et~al\mbox{.}}{2019}]%
        {vats2019multivariate}
\bibfield{author}{\bibinfo{person}{Dootika Vats}, \bibinfo{person}{James~M
  Flegal}, {and} \bibinfo{person}{Galin~L Jones}.}
  \bibinfo{year}{2019}\natexlab{}.
\newblock \showarticletitle{Multivariate output analysis for Markov chain Monte
  Carlo}.
\newblock \bibinfo{journal}{\emph{Biometrika}} \bibinfo{volume}{106},
  \bibinfo{number}{2} (\bibinfo{year}{2019}), \bibinfo{pages}{321--337}.
\newblock


\bibitem[\protect\citeauthoryear{Vehtari, Gelman, Simpson, Carpenter, and
  B{\"{u}}rkner}{Vehtari et~al\mbox{.}}{2020}]%
        {Vehtari2020}
\bibfield{author}{\bibinfo{person}{Aki Vehtari}, \bibinfo{person}{Andrew
  Gelman}, \bibinfo{person}{Daniel Simpson}, \bibinfo{person}{Bob Carpenter},
  {and} \bibinfo{person}{Paul-Christian B{\"{u}}rkner}.}
  \bibinfo{year}{2020}\natexlab{}.
\newblock \showarticletitle{{Rank-Normalization, Folding, and Localization: An
  Improved $\widehat{R}$ for Assessing Convergence of MCMC (with Discussion)}}.
\newblock \bibinfo{journal}{\emph{Bayesian Analysis}} \bibinfo{volume}{16},
  \bibinfo{number}{2} (\bibinfo{year}{2020}), \bibinfo{pages}{1--26}.
\newblock
\urldef\tempurl%
\url{https://doi.org/10.1214/20-ba1221}
\showDOI{\tempurl}
\showeprint[arxiv]{arXiv:1903.08008v5}


\bibitem[\protect\citeauthoryear{Villa, Petra, and Ghattas}{Villa
  et~al\mbox{.}}{2020}]%
        {Villa2020}
\bibfield{author}{\bibinfo{person}{Umberto Villa}, \bibinfo{person}{Noemi
  Petra}, {and} \bibinfo{person}{Omar Ghattas}.}
  \bibinfo{year}{2020}\natexlab{}.
\newblock \showarticletitle{{hIPPYlib User Manual: Version 2}}.
\newblock  (\bibinfo{date}{6} \bibinfo{year}{2020}).
\newblock
\urldef\tempurl%
\url{https://doi.org/10.6084/m9.figshare.12510578.v1}
\showDOI{\tempurl}


\bibitem[\protect\citeauthoryear{Villa, Petra, and Ghattas}{Villa
  et~al\mbox{.}}{2021}]%
        {VillaPetraGhattas21}
\bibfield{author}{\bibinfo{person}{Umberto Villa}, \bibinfo{person}{Noemi
  Petra}, {and} \bibinfo{person}{Omar Ghattas}.}
  \bibinfo{year}{2021}\natexlab{}.
\newblock \showarticletitle{{h{IPPYlib}: An Extensible Software Framework for
  Large-Scale Inverse Problems Governed by PDEs: Part I: Deterministic
  Inversion and Linearized Bayesian Inference}}.
\newblock \bibinfo{journal}{\emph{ACM Trans. Math. Softw.}}
  \bibinfo{volume}{47}, \bibinfo{number}{2}, Article \bibinfo{articleno}{16}
  (\bibinfo{date}{April} \bibinfo{year}{2021}), \bibinfo{numpages}{34}~pages.
\newblock
\showISSN{0098-3500}
\urldef\tempurl%
\url{https://doi.org/10.1145/3428447}
\showDOI{\tempurl}


\bibitem[\protect\citeauthoryear{Williams}{Williams}{1991}]%
        {Williams1991}
\bibfield{author}{\bibinfo{person}{David Williams}.}
  \bibinfo{year}{1991}\natexlab{}.
\newblock \bibinfo{booktitle}{\emph{Probability with Martingales}}.
\newblock \bibinfo{publisher}{Cambridge University Press}.
\newblock


\bibitem[\protect\citeauthoryear{Wolff, Collaboration, et~al\mbox{.}}{Wolff
  et~al\mbox{.}}{2004}]%
        {wolff2004monte}
\bibfield{author}{\bibinfo{person}{Ulli Wolff}, \bibinfo{person}{Alpha
  Collaboration}, {et~al\mbox{.}}} \bibinfo{year}{2004}\natexlab{}.
\newblock \showarticletitle{Monte Carlo errors with less errors}.
\newblock \bibinfo{journal}{\emph{Computer Physics Communications}}
  \bibinfo{volume}{156}, \bibinfo{number}{2} (\bibinfo{year}{2004}),
  \bibinfo{pages}{143--153}.
\newblock


\bibitem[\protect\citeauthoryear{Wong}{Wong}{2001}]%
        {Wong01}
\bibfield{author}{\bibinfo{person}{R. Wong}.} \bibinfo{year}{2001}\natexlab{}.
\newblock \bibinfo{booktitle}{\emph{Asymptotic Approximations of Integrals}}.
\newblock \bibinfo{publisher}{Society for Industrial and Applied Mathematics}.
\newblock
\urldef\tempurl%
\url{https://doi.org/10.1137/1.9780898719260}
\showDOI{\tempurl}
\showeprint{https://epubs.siam.org/doi/pdf/10.1137/1.9780898719260}


\bibitem[\protect\citeauthoryear{Zahm, Cui, Law, Spantini, and Marzouk}{Zahm
  et~al\mbox{.}}{2022}]%
        {ZahmCuiLawEtAl22}
\bibfield{author}{\bibinfo{person}{Olivier Zahm}, \bibinfo{person}{Tiangang
  Cui}, \bibinfo{person}{Kody Law}, \bibinfo{person}{Alessio Spantini}, {and}
  \bibinfo{person}{Youssef Marzouk}.} \bibinfo{year}{2022}\natexlab{}.
\newblock \showarticletitle{Certified dimension reduction in nonlinear
  {B}ayesian inverse problems}.
\newblock \bibinfo{journal}{\emph{Math. Comp.}} \bibinfo{volume}{91},
  \bibinfo{number}{336} (\bibinfo{year}{2022}), \bibinfo{pages}{1789--1835}.
\newblock


\end{thebibliography}

\end{document}